\documentclass[11pt]{article}
\usepackage{pxfonts}
\usepackage{yfonts}
\usepackage{dsfont}
\usepackage{graphicx}
\usepackage{relsize}

\def\bel{\begin{equation}\label}
\def\eeq{\end{equation}}
\def\ds{\displaystyle}
\def\endproof{\hphantom{MM}
\hfill\llap{$\square$}\goodbreak}

\def\mt{\longrightarrow}
\def\v{\vskip 1em}

\def\ve{\varepsilon}
\def\R{\mathds R}
\def\Z{\mathds Z}
\def\C{\mathfrak{B}}

\def\Cx{\mathds C}
\def\N{{\bf N}}

\def\Re{{\bf Re}}
\def\Im {{\bf Im}}
\def\S{{\bf S}}

\def\O{{\bf O}}
\def\P{{\bf P}}
\def\Q{{\bf Q}}
\def\D{{\bf D}}
\def\W{{\bf W}}
\def\J{{\bf J}}
\def\A{{\bf A}}
\def\B{{\bf B}}
\def\H{{\bf H}}
\def\L{{\bf L}}
\def\U{{\bf U}}
\def\V{{\bf V}}

\def\N{{\bf N}}

\def\p{{\partial}}

\def\i{{\bf i}}
\def\Tilde{\widetilde}
\def\Hat{\widehat}

\def\supp{{\bf supp}}

\def\I{{\bf I}}
\def\II{{\bf II}}
\def\III{{\bf III}}

\def\M{{\bf M}}
\def\G{{\bf G}}
\def\Rec{{\bf R}}

\def\Cup{{\bigcup}}

\def\alpha{\alphaup}
\def\beta{\betaup}
\def\gamma{\gammaup}
\def\delta{\deltaup}
\def\xi{{\xiup}}
\def\eta{{\etaup}}
\def\tau{{\tauup}}
\def\rho{{\rhoup}}
\def\phi{{\phiup}}
\def\psi{{\psiup}}
\def\lambda{{\lambdaup}}
\def\omega{\omegaup}
\def\varphi{{\varphiup}}
\def\gamma{{\gammaup}}
\def\c{{\bf c}}

\def\a{{\bf a}}
\def\b{{\bf b}}

\def\p{{\bf v}}
\def\h{{\bf h}}
\def\u{{\bf u}}
\def\({\left(}
\def\){\right)}
\def\vol{{\bf vol}}

\newtheorem{lemma}{Lemma}[section]

\newtheorem{remark}{Remark}[section]

\begin{document}
 \[\begin{array}{cc}\hbox{\LARGE{\bf Fractional integration with singularity on the light-cone}}
\end{array}\]

  \[\hbox{Zipeng Wang}\]
  \begin{abstract}
 ~~In this paper, we solve a classical problem investigated by Gelfand and Shilov on a family of fractional integral operators whose kernels have singularity appeared on the light-cone in $\R^{n+1}$. As an implication, we obtain a desired sharp estimate for the $\L^p\mt\L^q$-boundedness of cone multiplier operators having a negative order.

\[\hbox{\small{\bf Contents}}\]
\[\begin{array}{lr}\ds
\hbox{1 Introduction}~~~~~~~~~~~~~~~~~~~~~~~~~~~~~~~~~~~~~~~~~~~~~~~~~~~~~~~~~~~~~~~~~~~~~~~~~~~~~~~~~~~~~~~~~~~~~~~~~~~~~~~~~~~~~1
\\ \ds
~~~~~\hbox{1.1 Sketch on the proof of Theorem One}~~~~~~~~~~~~~~~~~~~~~~~~~~~~~~~~~~~~~~~~~~~~~~~~~~~~~~~~~~~3
\\ \ds
~~~~~\hbox{1.2 A pairing formulation}~~~~~~~~~~~~~~~~~~~~~~~~~~~~~~~~~~~~~~~~~~~~~~~~~~~~~~~~~~~~~~~~~~~~~~~~~~~~~~~~~~~~6
\\ \ds
\hbox{2 Proof of Theorem Two}~~~~~~~~~~~~~~~~~~~~~~~~~~~~~~~~~~~~~~~~~~~~~~~~~~~~~~~~~~~~~~~~~~~~~~~~~~~~~~~~~~~~~~~~~~~11
\\ \ds
\hbox{3 End-point estimates regarding $\II^{\alpha}_{j}$ and $\III^{\alpha}_{j}$}~~~~~~~~~~~~~~~~~~~~~~~~~~~~~~~~~~~~~~~~~~~~~~~~~~~~~~~~13
\\ \ds~~~~~
\hbox{3.1 ${^\sharp}\II^{\alpha~\beta}_{j}$ and ${^\sharp}\III^{\alpha~\beta}_{j}$} ~~~~~~~~~~~~~~~~~~~~~~~~~~~~~~~~~~~~~~~~~~~~~~~~~~~~~~~~~~~~~~~~~~~~~~~~~~~~~~~~~~~~~~~~~~~~~~~14
\\ \ds~~~~~
\hbox{3.2 ${^\flat}\II^{\alpha~\beta}_{j}$ and ${^\flat}\III^{\alpha~\beta}_{j}$} ~~~~~~~~~~~~~~~~~~~~~~~~~~~~~~~~~~~~~~~~~~~~~~~~~~~~~~~~~~~~~~~~~~~~~~~~~~~~~~~~~~~~~~~~~~~~~~~16
\\ \ds
\hbox{4 Proof of Theorem A}~~~~~~~~~~~~~~~~~~~~~~~~~~~~~~~~~~~~~~~~~~~~~~~~~~~~~~~~~~~~~~~~~~~~~~~~~~~~~~~~~~~~~~~~~~~~~~~~18
\\ \ds
\hbox{5 Proof of Proposition One}~~~~~~~~~~~~~~~~~~~~~~~~~~~~~~~~~~~~~~~~~~~~~~~~~~~~~~~~~~~~~~~~~~~~~~~~~~~~~~~~~~~~~~~21
\\ \ds~~~~~
\hbox{5.1 Proof of Lemma One}~~~~~~~~~~~~~~~~~~~~~~~~~~~~~~~~~~~~~~~~~~~~~~~~~~~~~~~~~~~~~~~~~~~~~~~~~~~~~~~~~~~~~~22
\\ \ds~~~~~
\hbox{5.2 Size of ${^\sharp}\Hat{\P}^{\alpha~\beta}_{j}$ and ${^\sharp}\Hat{\Q}^{\alpha~\beta}_{j}$}~~~~~~~~~~~~~~~~~~~~~~~~~~~~~~~~~~~~~~~~~~~~~~~~~~~~~~~~~~~~~~~~~~~~~~~~~~~~~~~~~~~~~25
\\ \ds
\hbox{6 Convolution with a kernel having singularity on $\mathds{S}^{n-1}$}~~~~~~~~~~~~~~~~~~~~~~~~~~~~~~~~~~~~~36
\\ \ds
\hbox{7 Proof of Proposition Two}~~~~~~~~~~~~~~~~~~~~~~~~~~~~~~~~~~~~~~~~~~~~~~~~~~~~~~~~~~~~~~~~~~~~~~~~~~~~~~~~~~~~~~~41
\\ \ds~~~~~
\hbox{7.1 Proof of Lemma Two}~~~~~~~~~~~~~~~~~~~~~~~~~~~~~~~~~~~~~~~~~~~~~~~~~~~~~~~~~~~~~~~~~~~~~~~~~~~~~~~~~~~~~~43
\\ \ds
\hbox{A Some estimates regarding Bessel functions}~~~~~~~~~~~~~~~~~~~~~~~~~~~~~~~~~~~~~~~~~~~~~~~~~~~~~~~47
\\ \ds
\hbox{B Fourier transform of $\Lambda^\alpha$}~~~~~~~~~~~~~~~~~~~~~~~~~~~~~~~~~~~~~~~~~~~~~~~~~~~~~~~~~~~~~~~~~~~~~~~~~~~~~~~~~~~~~~~~~50
\\ \ds
\hbox{References}~~~~~~~~~~~~~~~~~~~~~~~~~~~~~~~~~~~~~~~~~~~~~~~~~~~~~~~~~~~~~~~~~~~~~~~~~~~~~~~~~~~~~~~~~~~~~~~~~~~~~~~~~~~~~~~~~~~54
\end{array}\]
\end{abstract}  

\section{Introduction}
 \setcounter{equation}{0}
Let
$\lambda(\alpha)={n+1\over 2}-\delta(\alpha)$ and $\delta(\alpha)=\Big({n+1\over 2n}\Big)\alpha$ for $\alpha\in\Cx$. 
 As investigated by Gelfand and Shilov \cite{Gelfand-Shilov}, $\Lambda^\alpha$ is a distribution  in $\R^{n+1}$ defined by analytic continuation from
\bel{Lambda^alpha}
\begin{array}{cc}\ds
 \Lambda^\alpha(x,t)~=~\Gamma^{-1}\left(1-\lambda(\alpha)\right)
  \left({1\over t^2-|x|^2}\right)^{\lambda(\alpha)}_+,\qquad\hbox{\small{$\Re\lambda(\alpha)<1$}}
  \end{array}
\eeq
 where $\Gamma$ is   Gamma function.

 Consider $f$ to be a Schwartz function in $\R^{n+1}$. Define 
   \bel{I_alpha}
 \I^\alpha f(x,t)~=~\iint_{\R^{n+1}} f(x-u,t-r)\Lambda^\alpha(u,r)dudr,\qquad 0<\alpha<n.
 \eeq 
 
 At $\lambda(\alpha)=1, \delta(\alpha)={n-1\over 2}$ and $\alpha={n^2-n\over n+1}$, $\I^\alpha f(x,t)$ equals a constant multiple of
 \bel{surface measure}
 \begin{array}{lr}\ds
 \int_0^\infty \left\{\int_{\mathds{S}^{n-1}} f(x-r\omega, t-r)d\omega \right\}r^{n-2}  dr~=~\int_{\R^n} f(x-u,t-|u|)|u|^{-1} du
 \end{array}
 \eeq
 where $\omega$ denotes  Lebesgue measure on the unit sphere $\mathds{S}^{n-1}$.
 
 At $\lambda(\alpha)={n-1\over 2}, \delta(\alpha)=1$ and $\alpha={2n\over n+1}$, $\I^\alpha f$ is the solution of  inhomogeneous wave equation:
 \bel{Wave}
 \begin{array}{cc}\ds
 -\partial^2_t \u(x,t)+\sum_{i=1}^n \partial^2_{x_i} \u(x,t)~=~f(x,t),\qquad (x,t)\in\R^n\times\R;
 \\ \ds
 \u(x,0)~=~\partial_t \u(x,0)~=~0.
 \end{array}
 \eeq
For $0<\Re\delta(\alpha)<1$, the Fourier transform of $\Lambda^\alpha$ agrees with
\bel{Lambda Transform}
\begin{array}{lr}
  \Hat{\Lambda}^\alpha(\xi,\tau)~=~
\pi^{\lambda(\alpha)-\delta(\alpha)-1}\Gamma\left(\delta(\alpha)\right)
 \left\{ ~{\ds \left({1\over \tau^2-|\xi|^2}\right)^{\delta(\alpha)}_-} - \sin\pi\left(\delta(\alpha)-{1\over 2}\right){\ds \left( {1\over \tau^2- |\xi|^2}\right)^{\delta(\alpha)}_+}~\right\}  
 \end{array}
\eeq
whenever $ |\xi|\neq|\tau|$.

 The $\L^p\mt\L^q$-boundedness of  $\I^\alpha$
 is essentially equivalent to the cone multiplier problem of negative orders. 
Assert
 \bel{S^alpha}
\S^{-\alpha} f(x,t)~=~\iint_{\R^{n+1}} e^{2\pi\i (x\cdot\xi+t\tau)} \Hat{f}(\xi,\tau)\Hat{\hbox{V}}^\alpha(\xi,\tau)d\xi d\tau
\eeq
of which
$\Hat{\hbox{V}}^\alpha$ is a distribution defined in $\R^{n+1}$ by analytic continuation from 
\bel{V}
\begin{array}{cc}\ds
\Hat{\hbox{V}}^\alpha(\xi,\tau)~=~\Gamma^{-1}\left(1-\delta(\alpha)\right)
 \left({1\over \tau^2- |\xi|^2}\right)^{\delta(\alpha)}_+,\qquad   \hbox{\small{$\Re\delta(\alpha)<1$}}. 
\end{array}
 \eeq 
$\diamond$ {\small Throughout,  $\C>0$ is regarded as a generic constant depending on its subindices.}
\v
{\bf Theorem One} ~~{\it Let $\I^\alpha$ defined as (\ref{Lambda^alpha})-(\ref{I_alpha}) for $0<\alpha<n$.   We have
\bel{Result One}
\begin{array}{cc}\ds
 \left\|\I^\alpha f\right\|_{\L^q(\R^{n+1})}~\leq~\C_{p~q}~\left\|f\right\|_{\L^p(\R^{n+1})},\qquad 1<p<q<\infty
\\\\ \ds
\hbox{if and only if} 
\qquad{\alpha\over n}~=~{1\over p}-{1\over q},
\qquad
 {n-1\over 2n}+\left({n+1\over 2n}\right){\alpha\over n}~<~{1\over p}~<~{n+1\over 2n}+\left({n-1\over 2n}\right){\alpha\over n}.
\end{array}
 \eeq}
 
 {\bf Theorem Two} ~~{\it Let $\S^{-\alpha}$  defined in (\ref{S^alpha})-(\ref{V}) for $0<\alpha<n$.   We have
\bel{Result Two}
\begin{array}{cc}\ds
 \left\| \S^{-\alpha} f\right\|_{\L^q(\R^{n+1})}~\leq~\C_{p~q}~\left\|f\right\|_{\L^p(\R^{n+1})},\qquad 1<p<q<\infty
\\\\ \ds
\hbox{if and only if}\qquad 
 {\alpha\over n}~=~{1\over p}-{1\over q},
\qquad
 {n-1\over 2n}+\left({n+1\over 2n}\right){\alpha\over n}~<~{1\over p}~<~{n+1\over 2n}+\left({n-1\over 2n}\right){\alpha\over n}.
 \end{array}
 \eeq
 }
 
  The homogeneity condition ${\alpha\over n}={1\over p}-{1\over q}$   as an necessity is a standard exercise of changing dilations. The range of $p$  restricted to $ {n-1\over 2n}+\left({n+1\over 2n}\right){\alpha\over n}<{1\over p}<{n+1\over 2n}+\left({n-1\over 2n}\right){\alpha\over n}$ is less obvious. 
See the counter-examples  in the paper of Oberlin \cite{Oberlin} and Lee \cite{Lee}.

\begin{figure}[h]
\centering
\includegraphics[scale=0.52]{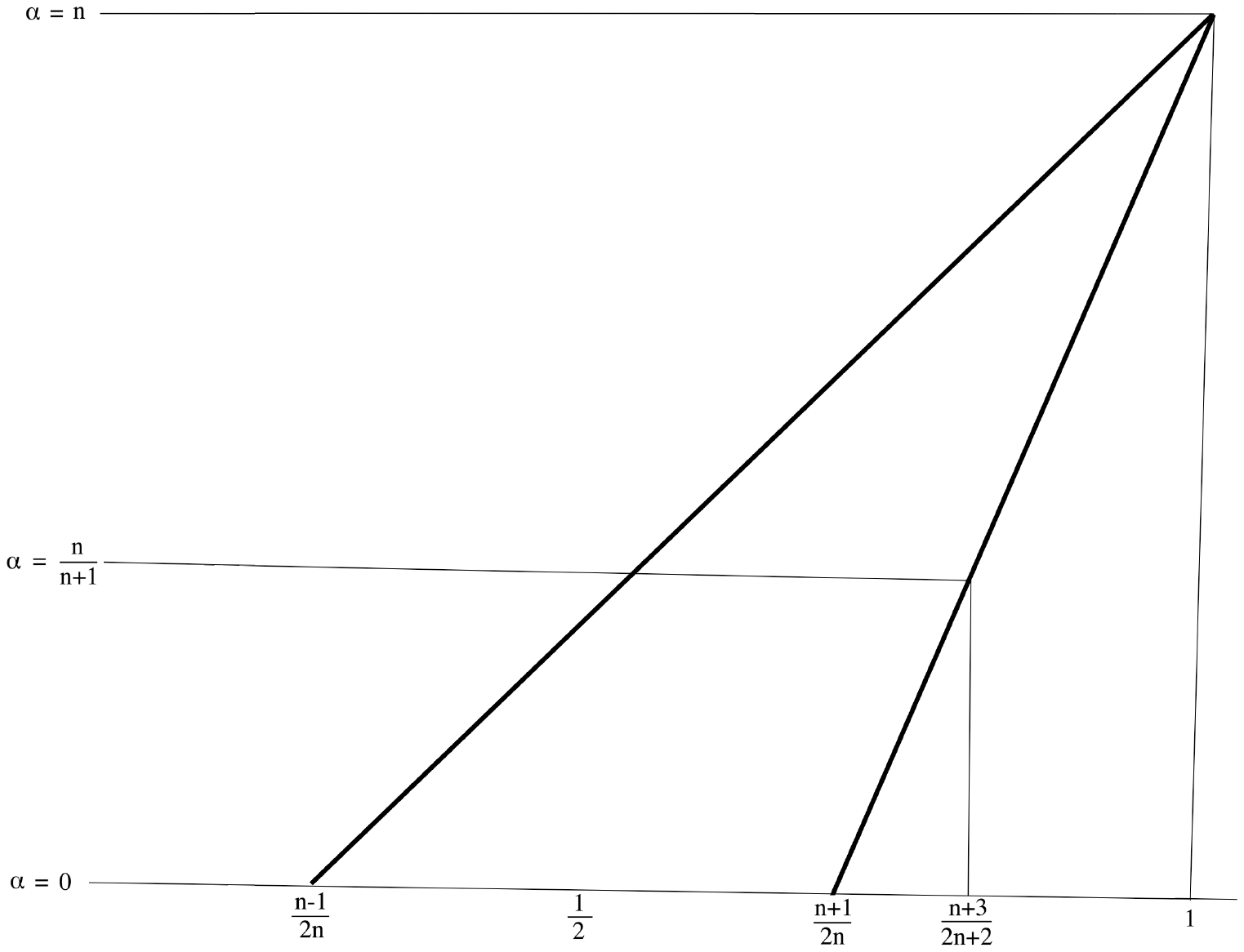}
\end{figure}

For the $\L^p\mt\L^q$-regularity of $\I^\alpha$, a number of remarkable results have been accomplished by Fefferman \cite{Fefferman''}, Ricci and Stein \cite{Ricci-Stein}, Greenleaf \cite{Greenleaf}, Littman \cite{Littman} and Oberlin \cite{Oberlin}.
\begin{remark} For ${n\over n+1}\leq \alpha<n$, {\bf Theorem One} is proved by Oberlin \cite{Oberlin} using certain mixed-norm bilinear estimates.
\end{remark}
In this paper, we give a completely new approach to prove {\bf Theorem One}. The proof of {\bf Theorem Two} is an implication.
Our framework is developed in another recent paper \cite{Wang} regarding the Bochner-Riesz summability.  

\subsection{Sketch on the proof of Theorem One} 
As investigated by Strichartz \cite{Strichartz},  $\Omega^\alpha$ is a distribution  defined in $\R^n$ by analytic continuation from
\bel{Omega^alpha}
\begin{array}{cc}\ds
\Omega^{\alpha}(x)~=~\pi^{-\lambda(\alpha)}\Gamma^{-1}\left(1-\lambda(\alpha)\right) 
 \left({1 \over1-|x|^2}\right)^{\lambda(\alpha)}_+,\qquad \hbox{\small{$\Re\lambda(\alpha)<1$}}. 
\end{array}
\eeq
Equivalently,  it can be  defined by its Fourier transform
 \bel{Omega^alpha Transform} 
\begin{array}{lr}\ds
\Hat{\Omega}^{\alpha}(\xi)~=~\left({1\over|\xi|}\right)^{{n\over 2}-\lambda(\alpha)} \J_{{n\over 2}-\lambda(\alpha)}\Big(2\pi|\xi|\Big)
\\\\ \ds~~~~~~~~~~
~=~\left({1\over|\xi|}\right)^{\delta(\alpha)-{1\over 2}} \J_{\delta(\alpha)-{1\over 2}}\Big(2\pi|\xi|\Big),\qquad \alpha\in\Cx
\end{array}
\eeq
where $\J$ is the Bessel function. See appendix A.

Recall (\ref{Lambda^alpha}).  $\Lambda^\alpha$ is a distribution  in $\R^{n+1}$ defined by analytic continuation from (\ref{Lambda^alpha}).
Consider $\I^\alpha f=f\ast\Lambda^\alpha$ as  (\ref{I_alpha}).
From (\ref{Omega^alpha}), we find
\[
 \begin{array}{lr}\ds
 \I^\alpha f(x,t)~=~\pi^{-\lambda(\alpha)}\Gamma^{-1}\left(1-\lambda(\alpha)\right) \iint_{|u|<|r|} f(x-u,t-r)  \left[{1\over r^2-|u|^2}\right]^{\lambda(\alpha)}dudr
\\\\ \ds~~~~~~~~~~~~~
~=~\iint_{\R^{n+1}} f(x-u,t-r)|r|^{-n}\Omega^\alpha\left({u\over r}\right)|r|^{\left({n+1\over n}\right)\alpha-1}dudr,\qquad\hbox{\small{$\Re\lambda(\alpha)<1$}}.
\end{array}
\]
By the principle of analytic continuation, we must have
 \bel{I_alpha rewrite}
 \begin{array}{lr}\ds
\I^\alpha f(x,t)~=~\pi^{\lambda(\alpha)}\iint_{\R^{n+1}} f(x-u,t-r)\Omega^\alpha_{|r|}(u)| r|^{\left({n+1\over n}\right)\alpha-1}dudr
\qquad
\hbox{\small{$0<\Re\alpha<n$}}
 \end{array}
\eeq 
where $\Omega^\alpha_r, r>0$ is the dilate of $\Omega^\alpha$ in the sense of distribution.

Alternatively, we find
\bel{I_alpha Equi}
\I^\alpha f(x,t)~=~\iint_{\R^{n+1}} e^{2\pi\i x\cdot\xi} \Hat{f}(\xi,t-r) \Hat{\Omega}^\alpha(r\xi) |r|^{\big({n+1\over n}\big)\alpha-1} d\xi dr.
\eeq
Recall $\delta(\alpha)=\Big({n+1\over 2n}\Big)\alpha$. For $0<\delta(\alpha)<1$, the Fourier transform of $\Lambda^\alpha$ agree with the function $\Hat{\Lambda}^\alpha(\xi,\tau)$ in (\ref{Lambda Transform}) whenever $|\xi|\neq|\tau|$. 

Let $0<\Re\alpha<1$. From (\ref{I_alpha rewrite})-(\ref{I_alpha Equi}), we further have
\bel{I_alpha Lambda}
\begin{array}{lr}\ds
\I^\alpha f(x,t)~=~\iint_{\R^{n+1}} e^{2\pi\i \big[x\cdot\xi+t\tau\big]} \Hat{f}(\xi,\tau) \Hat{\Lambda}^\alpha(\xi,\tau)d\xi d\tau,
\\\\ \ds
\Hat{\Lambda}^\alpha(\xi,\tau)~=~\int_\R e^{-2\pi\i r \tau}  \Hat{\Omega}^\alpha(r\xi) |r|^{\big({n+1\over n}\big)\alpha-1} dr
\\\\ \ds~~~~~~~~~~~~~
~=~\int_\R e^{-2\pi\i r \tau}  \left({1\over |r\xi|}\right)^{\delta(\alpha)-{1\over 2}}\J_{\delta(\alpha)-{1\over 2}}\big(2\pi|r\xi|\big) |r|^{\big({n+1\over n}\big)\alpha-1} dr,\qquad \hbox{\small{$|\xi|\neq|\tau|$}}.
\end{array}
\eeq
\v
{\bf Theorem A}~~{\it Let $\I^\alpha$ defined as (\ref{Lambda^alpha})-(\ref{I_alpha}) for $0<\Re\alpha<1$.   We have
\bel{Result A}
\begin{array}{cc}\ds
 \left\| \I^\alpha f\right\|_{\L^q(\R^{n+1})}~\leq~\C_{p~q}~e^{\c|\Im\alpha|}~\left\|f\right\|_{\L^p(\R^{n+1})}
 \\\\ \ds
 \hbox{if}\qquad {\Re\alpha\over n}~=~{1\over p}-{1\over q},
\qquad
 {n-1\over 2n}+{\Re\alpha\over n}~<~{1\over p}~<~{n+1\over 2n}.
 \end{array} 
\eeq}
 
 On the other hand, we find
 \bel{I n}
 \left\|\I^\alpha f\right\|_{\L^\infty(\R^{n+1})}~\leq~\C~e^{\c|\Im\alpha|}~\left\|f\right\|_{\L^1(\R^{n+1})},\qquad\hbox{\small{$\Re\alpha=n$}}.
 \eeq 
 From (\ref{Result A}) and (\ref{I n}), we obtain (\ref{Result One}) for every $0<\alpha<n$   as desired
 by applying Stein interpolation theorem \cite{Stein}.
 
 Let $\varphi\in\mathcal{C}^\infty_o(\R)$ be a smooth bump-function such that $\varphi(r)=1$ if $|r|\leq1$ and $\varphi(r)=0$ if $|r|>2$. Define
 \bel{varphi_j}
 \begin{array}{cc}\ds
 \phi_<(\xi)~=~\varphi(|\xi|),
 \\\\ \ds
  \phi_j(\xi)~=~\varphi\Big( 2^{-j}|\xi|\Big)-\varphi\Big(2^{-j+1}|\xi|\Big),\qquad j\in\Z^+.
  \end{array}
 \eeq
Recall $\Hat{\Omega}^\alpha(\xi)$ from (\ref{Omega^alpha Transform}). Define
 \bel{Omega alpha j}
 \begin{array}{cc}\ds
\Hat{\Omega}^\alpha_<(\xi)~=~ \left({1\over|\xi|}\right)^{\delta(\alpha)-{1\over 2}} \J_{\delta(\alpha)-{1\over 2}}\Big(2\pi|\xi|\Big) \phi_<(\xi),
\\\\ \ds
 \Hat{\Omega}^\alpha_j(\xi)~=~ \left({1\over|\xi|}\right)^{\delta(\alpha)-{1\over 2}} \J_{\delta(\alpha)-{1\over 2}}\Big(2\pi|\xi|\Big) \phi_j(\xi) ,\qquad j>0.
 \end{array}
 \eeq
By using the integral formula of Bessel functions in (\ref{Bessel}), we  have
 \bel{local inte<}
 \Hat{\Omega}^\alpha_< (\xi)~=~ \pi^{\delta(\alpha)-1}\Gamma\left(\delta(\alpha)\right)\phi_<(\xi)\int_{-1}^1 e^{2\pi\i |\xi|s } (1-s^2)^{\delta(\alpha)-1} ds.
   \eeq 
 Observe that $ \Hat{\Omega}^\alpha_<$ has a compact support for $|\xi|<2$. Moreover, a direct computation shows
 \bel{Diff Ineq} 
 \left|\partial_\xi^\gamma \Hat{\Omega}^\alpha_< (\xi)\right|~\leq~\C_\gamma~|\xi|^{-|\gamma|}
 \eeq
 for every multi-index $\gamma$. 
 
 By using the differential inequality in (\ref{Diff Ineq}), we find
 \bel{Omega<size}
 \begin{array}{lr}\ds
 \left|\Omega^\alpha_<(x)\right|~\leq~\C_{\Re\alpha}~\left({1\over 1+|x|}\right)^n
 \\\\ \ds~~~~~~~~~~~~
 ~\leq~\C_{\Re\alpha}~|x|^{\Re\alpha-n},\qquad x\neq0,\qquad\hbox{\small{$0<\Re\alpha<n$}}.
 \end{array}
 \eeq
 Let $0<\a<\N$. Recall a classical result due to Hardy and Littlewood \cite{Hardy-Littlewood} and Sobolev \cite{Sobolev}.
 \v
 {\bf Hardy-Littlewood-Sobolev theorem}~~{Let $0<\a<\N$. We have
 \bel{Hardy-Littlewood-Sobolev inequality}
 \begin{array}{cc}\ds
\left\{\int_{\R^\N} \left\{\int_{\R^\N} |f(x-y)||y|^{\a-\N}dy\right\}^q dx\right\}^{1\over q}~\leq~\C_{p~q}~
\left\| f\right\|_{\L^p(\R^\N)},\qquad 1<p<q<\infty
\\\\ \ds
\hbox{if and only if}\qquad {\a\over\N}~=~{1\over p}-{1\over q}.
\end{array}
\eeq}
 
For $0<\Re\alpha<n$, define
 \bel{I_alpha <}
 \begin{array}{rl}\ds
\I^\alpha_< f(x,t)~=~\pi^{\lambda(\alpha)}\iint_{\R^{n+1}} f(x-u,t-r)|r|^{-n}\Omega^\alpha_<\left({u\over r}\right)| r|^{\left({n+1\over n}\right)\alpha-1}dudr.
 \end{array}
\eeq 
Let ${\Re\alpha\over n}={1\over p}-{1\over q}, 1<p<q<\infty$. We have
 \bel{I<EST}
 \begin{array}{lr}\ds
 \left\|\I^\alpha_< f\right\|_{\L^q(\R^{n+1})}~\lesssim~
\left\{\iint_{\R^{n+1}} \left\{ \iint_{\R^{n+1}} |f(x-u,t-r)| |r|^{-n}\left|{u\over r}\right|^{\Re\alpha-n} |r|^{\big({n+1\over n}\big)\Re\alpha-1} dudr\right\}^q dxdt\right\}^{1\over q}
~~
\hbox{\small{by (\ref{Omega<size})}}
 \\\\ \ds
 \leq~   \left\{\int_\R\left\{ \int_\R\left\{\int_{\R^n} \left\{ \int_{\R^n} |f(x-u,t-r)| |r|^{-n}\left|{u\over r}\right|^{\Re\alpha-n} du\right\}^q dx\right\}^{1\over q} |r|^{\big({n+1\over n}\big)\Re\alpha-1} dr \right\}^q dt\right\}^{1\over q}
 \\ \ds~~~~~~~~~~~~~~~~~~~~~~~~~~~~~~~~~~~~~~~~~~~~~~~~~~~~~~~~~~~~~~~~~~~~~~~~~~~~~~~~~~~~
 \hbox{\small{by Minkowski integral inequality}}
 \\\\ \ds
 \leq~\C_{p~q}  \left\{\int_\R\left\{ \int_\R \left\| f(\cdot,t-r)\right\|_{\L^p(\R^n)}|r|^{{n\over q}-{n\over p}} |r|^{\big({n+1\over n}\big)\Re\alpha-1} dr \right\}^q dt\right\}^{1\over q} 
~~
 \hbox{\small{by (\ref{Hardy-Littlewood-Sobolev inequality}) with $\a=\Re\alpha$ and $\N=n$}}
\\\\ \ds
=~\C_{p~q}  \left\{\int_\R\left\{ \int_\R \left\| f(\cdot,t-r)\right\|_{\L^p(\R^n)} |r|^{{\Re\alpha\over n}-1} dr \right\}^q dt\right\}^{1\over q} 
\\\\ \ds
 \leq~\C_{p~q}~\left\| f\right\|_{\L^p(\R^{n+1})}\qquad  \hbox{\small{by (\ref{Hardy-Littlewood-Sobolev inequality}) with $\a=\Re\alpha$ and $\N=1$}}.
 \end{array}
\eeq 
Let $0<\Re\alpha<1$. Define
  \bel{I_alpha j}
 \begin{array}{rl}\ds
\I^\alpha_j f(x,t)~=~\iint_{\R^{n+1}} e^{2\pi\i x\cdot\xi} \Hat{f}(\xi,t-r) \Hat{\Omega}^\alpha_j(r\xi) |r|^{\big({n+1\over n}\big)\alpha-1} d\xi dr,\qquad j>0.
 \end{array}
\eeq 
Recall the Fourier transform of $\Omega^\alpha_j$ in (\ref{Omega alpha j}). Alternatively, we have
  \bel{I rewrite Transform}
  \begin{array}{cc}\ds
 \I^\alpha_j f(x,t)~=~\iint_{\R^{n+1}} e^{2\pi\i \big[x\cdot\xi+t\tau\big]}\Hat{f}(\xi,\tau)\Lambda^\alpha_j(\xi,\tau)d\xi d\tau,
  \\\\ \ds
\Lambda^\alpha_j(\xi,\tau)~=~\int_\R e^{-2\pi\i r\tau} \left({1\over|r\xi|}\right)^{\delta(\alpha)-{1\over 2}} \J_{\delta(\alpha)-{1\over 2}}\Big(2\pi|r\xi|\Big)\phi_j(r\xi) |r|^{\left({n+1\over n}\right)\alpha-1} dr. 
  \end{array}
 \eeq 
 Our aim is to show
 \bel{Result A j} 
 \begin{array}{cc}\ds
  \left\| \I^\alpha_j f\right\|_{\L^q(\R^{n+1})} ~\leq~\C_{\Re\alpha~p}~e^{\c|\Im\alpha|}~2^{-\ve j}~\left\| f\right\|_{\L^p(\R^{n+1})}, 
  \\\\ \ds
  {\Re\alpha\over n}~=~{1\over p}-{1\over q},\qquad {n-1\over 2n}+{\Re\alpha\over n}~< ~{1\over p}~<~ {n+1\over 2n}  
 \end{array}
  \eeq
for every $j>0$ and some $\ve=\ve(\Re\alpha)>0$. 

In order to prove (\ref{Result A j}), we need to further split $\I^\alpha_j, j>0$ into two  partial operators.

\subsection{A pairing formulation}
First, we introduce an ingenious dyadic decomposition  initially established by Fefferman \cite{Fefferman'}-\cite{Fefferman} and later refined by Seeger, Sogge and Stein \cite{Seeger-Sogge-Stein} in their study of Fourier integral operators. 

Given $j>0$,  $\left\{ \xi^\nu_j\right\}_\nu$ is a collection of points  almost equally distributed on  $\mathds{S}^{n-1}$ having a grid length between $2^{-j/2-1}$ and $2^{-j/2}$:   

{\bf (1)} $|\xi^\mu_j-\xi^\nu_j|\ge 2^{-j/2-1}$ for every $\xi^\mu_j\neq\xi^\nu_j$. 

{\bf (2)} For any $u\in \mathds{S}^{n-1}$, there is a $\xi^\nu_j$ in the open set $\{\xi\in\mathds{S}^{n-1}\colon |\xi-u|<2^{-j/2+1}\}$.

Let $\varphi\in\mathcal{C}^\infty_o(\R)$  be a smooth cut-off function such that $\varphi(t)=1$ if $|t|\leq1$ and $\varphi(t)=0$ for $|t|>2$. Define
\bel{phi^v_j intro}
\varphi^\nu_j(\xi)~=~{\ds\varphi\Bigg[2^{j/2}\left|{\xi\over |\xi|}-\xi^{\nu}_{j}\right|\Bigg]\over \ds \sum_{\nu\colon\xi^\nu_j\in\mathds{S}^{n-1}} \varphi\Bigg[2^{j/2}\left|{\xi\over |\xi|}-\xi^{\nu}_{j}\right|\Bigg]}
\eeq
whose support is contained in the cone
\bel{Gamma_j}
\Gamma^\nu_j~=~\Bigg\{\xi\in\R^n\colon \left|{\xi\over |\xi|}-\xi^{\nu}_{j}\right|<2^{-j/2+1}\Bigg\}.
\eeq
\begin{figure}[h]
\centering
\includegraphics[scale=0.30]{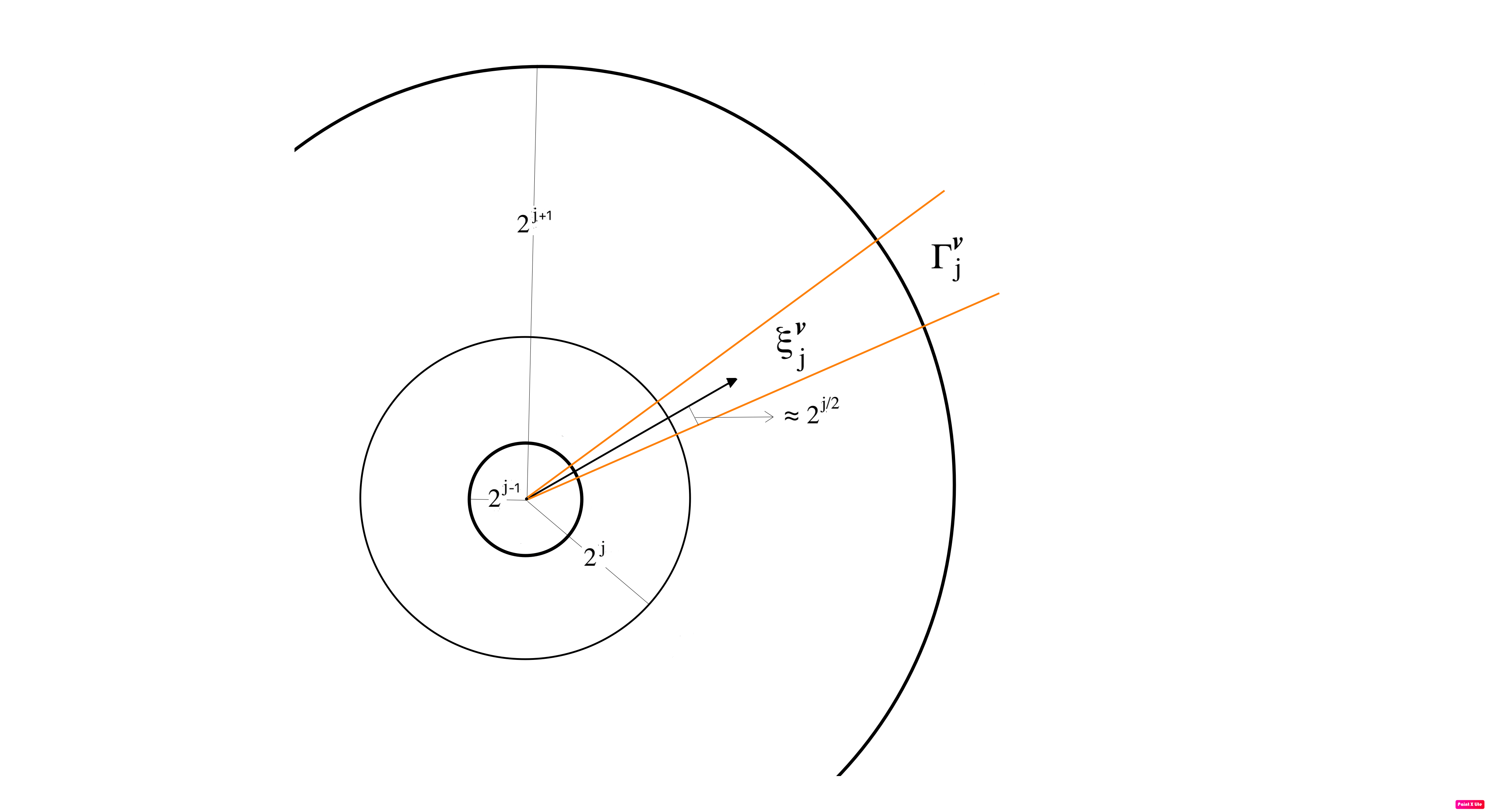}
\end{figure}

Next,  $\L_\nu$ is an $n\times n$-orthogonal matrix with $\det\L_\nu=1$ and 
$\L_\nu^T \xi^\nu_j=(1,0)^T\in\R\times\R^{n-1}$.
Define
 \bel{Psi mu j alpha}
 \begin{array}{cc}\ds
\Psi^{\alpha}_{j}(x)~=~\sum_{\nu\colon\xi^\nu_j\in\mathds{S}^{n-1}}  \Psi^{\alpha~\nu}_{j}(x),
 \\\\ \ds 
 \Psi^{\alpha~\nu}_{j}(x)~=~\varphi\left[ 2^{\big(1-\sigma\big)j}\left|1-\left(\L_\nu^T x\right)_1\right|\right] 
\prod_{i=2}^n \varphi\left[2^{\big[{1\over 2}-\sigma\big]j}\left|\left(\L_\nu^T x\right)_i\right|\right]
\\ \ds~~~~~~~~~~~~~
~+~\varphi\left[ 2^{\big(1-\sigma\big)j}\left|1+\left(\L_\nu^T x\right)_1\right|\right] 
\prod_{i=2}^n \varphi\left[2^{\big[{1\over 2}-\sigma\big]j}\left|\left(\L_\nu^T x\right)_i\right|\right].
\end{array}
\eeq
Each $\Psi^{\alpha~\nu}_{j}$ is supported in the rectangles
\bel{Rectangles}
\begin{array}{cc}\ds
\hbox{R}^\nu_{j}~=~\hbox{R}^{\nu~+}_{j}\cup \hbox{R}^{\nu~-}_{j},
\\\\ \ds
\hbox{R}^{\nu~+}_{j}~=~\Bigg\{x\in\R^n\colon1-2^{\big(\sigma-1\big)j+1}<\left(\L_\nu^T x\right)_1<1+2^{\big(\sigma-1\big)j+1},~\left|\left(\L_\nu^T x\right)_i\right|<2^{\big[\sigma-{1\over 2}\big]j+1},~ i=2,\ldots,n
\Bigg\},
\\\\ \ds
\hbox{R}^{\nu~-}_{j}~=~ \Bigg\{x\in\R^n\colon-1-2^{\big(\sigma-1\big)j+1}<\left(\L_\nu^T x\right)_1<-1+2^{\big(\sigma-1\big)j+1},~\left|\left(\L_\nu^T x\right)_i\right|<2^{\big[\sigma-{1\over 2}\big]j+1},~ i=2,\ldots,n
\Bigg\}.
\end{array}
\eeq
Moreover, $\hbox{R}^\mu_{j}\cap \hbox{R}^\nu_{j}=\emptyset$ if   $\left|\xi^\mu_j-\xi^\nu_j\right|\ge\c2^{-j/2}2^{\sigma j}$ for  $\c$  large.  
\begin{figure}[h]
\centering
\includegraphics[scale=0.52]{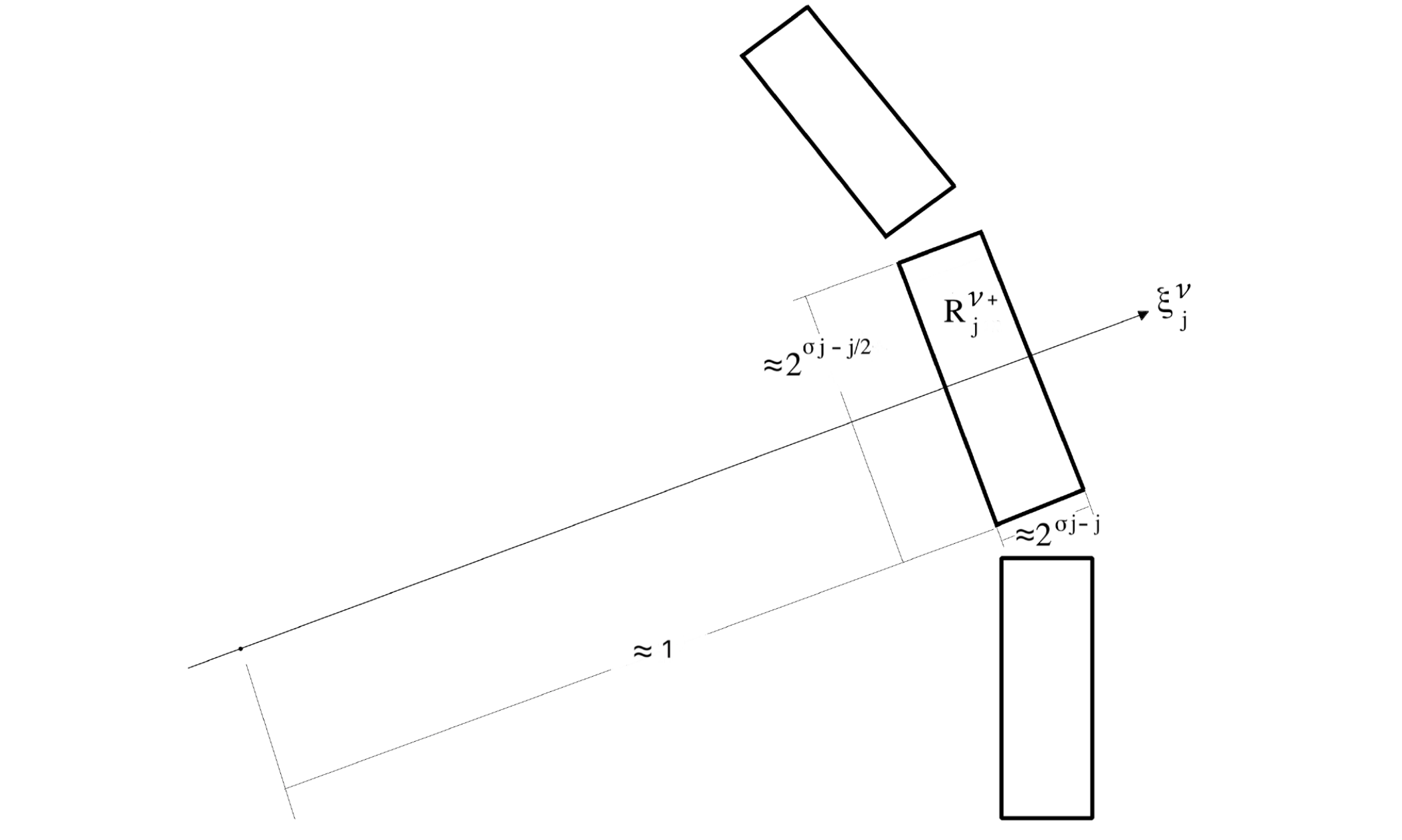}
\end{figure}

Let $\varphi^\nu_j$ defined in (\ref{phi^v_j intro}).
Consider
\bel{P v j}
 \begin{array}{cc}\ds
 \Hat{\Omega}_j^\alpha(\xi)~=~
\sum_{ \xi^\nu_j\in\mathds{S}^{n-1}}  \Hat{\Omega}_j^{\alpha~\nu}(\xi),
\qquad
 \Hat{\Omega}_j^{\alpha~\nu}(\xi)~=~\varphi^\nu_j(\xi) \Hat{\Omega}^\alpha_j (\xi),\qquad j>0.
  \end{array}
\eeq
Let $\Psi^{\alpha}_j$ and $\Psi^{\alpha~\nu}_j$ defined in (\ref{Psi mu j alpha}).
For every $j>0$, we assert
   \bel{II_j U intro}
 \begin{array}{ccc}\ds
 \II^\alpha_jf(x,t)~=~\iint_{\R^{n+1}} f(x-u,t-r) \U^\alpha_{j~|r|}(u) |r|^{\big({n+1\over n}\big)\alpha-1}dudr, 
 \\\\ \ds
 \U^\alpha_j(x) ~=~
 \sum_{\xi^\nu_j\in\mathds{S}^{n-1}}  \U^{\alpha~\nu}_j(x), 
\qquad
  \U^{\alpha~\nu}_j(x)~=~\Omega^{\alpha~\nu}_j(x) \left[1-\Psi^{\alpha~\nu}_j(x)\right], 
    \end{array}
  \eeq 
 \bel{III_j V intro}
 \begin{array}{ccc}\ds
 \III^\alpha_jf(x,t)~=~\iint_{\R^{n+1}} f(x-u,t-r) \V^\alpha_{j~|r|}(u)|r|^{\big({n+1\over n}\big)\alpha-1}dudr,
 \\\\ \ds
 \V^\alpha_j(x) ~=~
 \sum_{\xi^\nu_j\in\mathds{S}^{n-1}}  \V^{\alpha~\nu}_j(x), 
\qquad
  \V^{\alpha~\nu}_j(x)~=~\Omega^{\alpha~\nu}_j(x) \Psi^{\alpha~\nu}_j(x) 
    \end{array}
 \eeq 
 where $\U^\alpha_{j~r}(x)=r^{-n}  \U^\alpha_j\big({x\over r}\big), r>0$ and  $\V^\alpha_{j~r}(x)=r^{-n}  \V^\alpha_j\big({x\over r}\big), r>0$.

From (\ref{II_j U intro})-(\ref{III_j V intro}), we find
\[
\begin{array}{lr}\ds
\U^{\alpha}_j(x)+\V^{\alpha}_j(x)
~=~\sum_{\xi^\nu_j\in\mathds{S}^{n-1}} \Omega^{\alpha~\nu}_j(x) \left[1-\Psi^{\alpha~\nu}_j(x)\right] +
\Omega^{\alpha~\nu}_j(x) \Psi^{\alpha~\nu}_j(x)
\\\\ \ds~~~~~~~~~~~~~~~~~~~~~~~~
~=~\sum_{\xi^\nu_j\in\mathds{S}^{n-1}} \Omega^{\alpha~\nu}_j(x)
~=~\Omega^\alpha_j(x).
\end{array}
\]
In section 4, we will explicitly define
\[
 \I^{\alpha~z}_j f(x,t)~=~\II^{\alpha~z}_j f(x,t)+
 \III^{\alpha~z}_j f(x,t),\qquad \hbox{\small{$0\leq\Re z\leq1$}}
\]
of which $\II^{\alpha~z}_j$ and  $\III^{\alpha~z}_j$
 are two families of analytic operators.
 In particular, we have
 $\II^{\alpha~{1\over n}}_jf(x,t)=\II^{\alpha}_jf(x,t)$ and $\III^{\alpha~{1\over n}}_jf(x,t)=\III^{\alpha}_jf(x,t)$.
 
$\bullet$ For every $j>0$,  there is an $\ve=\ve(\Re\alpha)>0$ and
\bel{II 0 bad}
\begin{array}{cc}\ds
\left\{\int_\R\left\{\int_{\R^n} \left| \II^{\alpha~0+\i\Im z}_j f (x,t)\right|^q dx\right\}^{2\over q} dt \right\}^{1\over 2}
 ~\leq~\C_{\Re\alpha}~e^{\c|\Im\alpha|}~2^{j/2}2^{-\ve j}~\left\| f\right\|_{\L^2(\R^{n+1})}
 \\\\ \ds
\hbox{if}\qquad  {\Re\alpha\over 2n}~=~{1\over 2}-{1\over q},\qquad 2<q<\infty;
 \end{array}
\eeq
\bel{III 0 good}
\begin{array}{cc}\ds
\left\{\int_\R\left\{\int_{\R^n} \left| \III^{\alpha~0+\i\Im z}_j f (x,t)\right|^q dx\right\}^{2\over q} dt \right\}^{1\over 2}
 ~\leq~\C_{\Re\alpha}~e^{\c|\Im\alpha|}~2^{-\ve j}~\left\| f\right\|_{\L^2(\R^{n+1})},
 \\\\ \ds
\hbox{if}\qquad   {\Re\alpha\over 2n}~=~{1\over 2}-{1\over q},\qquad 2<q<\infty.
 \end{array} 
\eeq
Observe that at $\Re z=0$, $\III_j^{\alpha~z}$ is a {\it good} partial operator of $\I^{\alpha~z}_j$ whose $\L^q\L^2$-mixed norm is bound by $\C_{\Re\alpha}e^{\c|\Im\alpha|}2^{-\ve j}$. 
On the other hand,  
$\II^{\alpha~z}_j$ is the {\it bad} one, satisfying (\ref{II 0 bad}). 
\v
$\bullet$ For every $j>0$,  we have
\bel{II 1 good}
\begin{array}{cc}\ds
\left\{\int_\R\left\{\int_{\R^n} \left| \II^{\alpha~1+\i\Im z}_j f (x,t)\right|^q dx\right\}^{s\over q} dt \right\}^{1\over s}
~\leq~\C_{N~\Re\alpha}~e^{\c|\Im\alpha|}~2^{-N j}~\left\| f\right\|_{\L^p(\R^{n+1})},\quad N\ge0
 \\\\ \ds
\hbox{if}\qquad \left({n+1\over 2n}\right)\Re\alpha~=~{1\over p}-{1\over q},\qquad \Re\alpha~=~{1\over p}-{1\over s},\qquad 1<p<q\leq s<\infty;
 \end{array}
\eeq
\bel{III 1 bad}
\begin{array}{cc}\ds
\left\{\int_\R\left\{\int_{\R^n} \left| \III^{\alpha~1+\i\Im z}_j f (x,t)\right|^q dx\right\}^{s\over q} dt \right\}^{1\over s} ~\leq~\C_{\Re\alpha}~e^{\c|\Im\alpha|}~\left\| f\right\|_{\L^p(\R^{n+1})}
 \\\\ \ds
\hbox{if}\qquad \left({n+1\over 2n}\right)\Re\alpha~=~{1\over p}-{1\over q},\qquad \Re\alpha~=~{1\over p}-{1\over s},\qquad 1<p<q\leq s<\infty.
 \end{array} 
 \eeq
At $\Re z=1$, $\III^{\alpha~z}_j$ is a {\it bad} partial operator of $\I^{\alpha~z}_j$ satisfying (\ref{III 1 bad}). On the other hand,  
$\II^{\alpha~z}_j$ becomes the {\it good} one whose $\L^q\L^s$-norm is bounded by $\C_{N~\Re\alpha}e^{\c|\Im\alpha|}2^{-N j}$, $N\ge0$.

Recall $\II^{\alpha~{1\over n}}_j=\II^{\alpha}_j$ and $\III^{\alpha~{1\over n}}_j=\III^{\alpha}_j$. 
Take into account for 
\[
\begin{array}{cc}\ds
{1\over n}\left({n+2\over 2n}\right)\Re\alpha+{n-1\over n} {\Re\alpha\over 2n}
~=~\left[{n+1\over 2n}+{n-1\over 2n}\right]{\Re\alpha\over n}
~=~{\Re\alpha\over n},
\\\\ \ds
{1\over n}\Re\alpha +{n-1\over n} {1\over 2}~=~{n-1\over 2n}+{\Re\alpha\over n},
\qquad
{1\over n}+{n-1\over n}{1\over 2}~=~{n+1\over 2n}.
\end{array}
\]
From (\ref{II 0 bad}) and (\ref{II 1 good}) with $N$ chosen sufficiently large, we have
\bel{II EST}
 \begin{array}{cc}\ds
  \left\| \II^{\alpha}_j f\right\|_{\L^q(\R^{n+1})} ~\leq~\C_{\Re\alpha}~e^{\c|\Im\alpha|}~2^{-\ve\big({n-1\over n}\big) j}~\left\| f\right\|_{\L^p(\R^{n+1})}
  \\\\ \ds
\hbox{if}\qquad {\Re\alpha\over n}~=~{1\over p}-{1\over q},\qquad    {n-1\over 2n}+{\Re\alpha\over n}~< ~{1\over p}~<~ {n+1\over 2n}  
 \end{array}
\eeq
by applying the mixed-norm Stein interpolation theorem due to Benedek and Panzone \cite{Benedek-Panzone}. 
 
From (\ref{III 0 good}) and (\ref{III 1 bad}), we find
\bel{III EST}
 \begin{array}{cc}\ds
  \left\| \III^{\alpha}_j f\right\|_{\L^q(\R^{n+1})} ~\leq~\C_{\Re\alpha}~e^{\c|\Im\alpha|}~2^{-\ve\big({n-1\over n}\big) j}~\left\| f\right\|_{\L^p(\R^{n+1})}
  \\\\ \ds
\hbox{if}\qquad {\Re\alpha\over n}~=~{1\over p}-{1\over q},\qquad   {n-1\over 2n}+{\Re\alpha\over n}~<~ {1\over p}~<~ {n+1\over 2n}
 \end{array}
\eeq
by applying the mixed-norm Stein interpolation theorem again.

By putting together (\ref{II EST})-(\ref{III EST}), we obtain
(\ref{Result A j}) as desired.

\section{Proof of Theorem Two}
\setcounter{equation}{0}
Let $\hbox{U}^\alpha$, $\hbox{V}^\alpha$  be distributions defined by analytic continuation from (\ref{UV}). 
For $0<\Re\delta(\alpha)<1$, we have
\bel{UV=0}
 \hbox{U}^{\alpha\over 2}(\xi,\tau)\hbox{V}^{\alpha\over 2}(\xi,\tau)~=~0
 \eeq
 and
\bel{UV double}
 \begin{array}{lr}\ds 
\Gamma\left(1-\delta(\alpha)\right)\hbox{U}^\alpha(\xi,\tau)
~=~\left[\hbox{\small{$\Gamma\left(1-\delta\left({\alpha\over 2}\right)\right)$}}\hbox{U}^{\alpha\over 2}(\xi,\tau)\right]^2,
\\\\ \ds
\Gamma\left(1-\delta(\alpha)\right)\hbox{V}^\alpha(\xi,\tau)
~=~\left[\hbox{\small{$\Gamma\left(1-\delta\left({\alpha\over 2}\right)\right)$}}\hbox{V}^{\alpha\over 2}(\xi,\tau)\right]^2,\qquad\hbox{\small{$|\xi|\neq|\tau|$}}.
\end{array}
\eeq
By using the formula  in (\ref{Lambda Transform}), we write
\bel{Lambda transform UV}
\begin{array}{lr}
\c_\alpha\Gamma^{-1}(\delta(\alpha))\Hat{\Lambda}^\alpha(\xi,\tau)~=~\Gamma\left(1-\delta(\alpha)\right)
\Bigg\{ \hbox{U}^\alpha(\xi,\tau)-\sin\pi\left(\delta(\alpha)-{1\over 2}\right)\hbox{V}^\alpha(\xi,\tau)\Bigg\},
\\ \ds~~~~~~~~~~~~~~~~~~~~~~~~~~~~~~~~~~~~~~~~~~~~~~~~~~~~~~~~~~~~~~~~~~~~~~~~~~~~
\hbox{\small{$\c_\alpha=\pi^{1+\delta(\alpha)-\lambda(\alpha)}$}},\qquad \hbox{\small{$|\xi|\neq|\tau|$}}.
 \end{array}
\eeq
Note that $\delta(\alpha)=\Big({n+1\over 2n}\Big)\alpha$. We find $0<\Re\delta(\alpha)<1 \Leftrightarrow 0<\Re\alpha<{2n\over n+1}$. 
Let $\I^\alpha f=f\ast\Lambda^\alpha$ defined as (\ref{Lambda^alpha})-(\ref{I_alpha}) for $0<\Re\alpha<1$. The regularity estimate in (\ref{Result A}) shows that $\c_\alpha\Gamma^{-1}(\delta(\alpha))\Hat{\Lambda}^\alpha(\xi,\tau)$ 
is a $\L^p\mt\L^q$-Fourier multiplier, with $p,q,$ and $\Re\alpha$  satisfying 
\bel{Formula A}
{\Re\alpha\over n}~=~{1\over p}-{1\over q},
\qquad
 {n-1\over 2n}+{\Re\alpha\over n}~<~{1\over p}~<~{n+1\over 2n}.
\eeq
Consider 
\bel{double multiplier}
\c_{\alpha/2}^2\hbox{\small{$\Gamma^{-2}\left(\delta\left({\alpha\over 2}\right)\right)$}}~f\ast\Lambda^{\alpha\over 2}\ast\Lambda^{\alpha\over 2},\qquad\hbox{\small{$0<\Re\alpha<1$}}.
\eeq  
From (\ref{UV=0})-(\ref{UV double}) and (\ref{Lambda transform UV}), we find
\bel{Multiplier}
\begin{array}{lr}\ds
\c_{\alpha/2}^2\hbox{\small{$\Gamma^{-2}\left(\delta\left({\alpha\over 2}\right)\right)$}}~\Big[\Hat{\Lambda}^{\alpha\over 2}(\xi,\tau)\Big]^2
 ~=~\Gamma\left(1-\delta(\alpha)\right)
\Bigg\{\hbox{U}^\alpha(\xi,\tau)+\hbox{\small{$\sin^2\pi\left[\delta\left({\alpha\over 2}\right)-{1\over 2}\right]$}}\hbox{V}^\alpha(\xi,\tau)\Bigg\},
\\ \ds~~~~~~~~~~~~~~~~~~~~~~~~~~~~~~~~~~~~~~~~~~~~~~~~~~~~~~~~~~~~~~~~~~~~~~~~~~~~~~~~~~~~~~~~~~~~~~~~~~~~~~~~~~~~~~~~~~~~~~~
 \hbox{\small{$|\xi|\neq|\tau|$}}.
\end{array}
\eeq
Let $0<\Re\alpha<1$. By applying (\ref{Result A})  two times, we have
\bel{mapping1}
\begin{array}{cc}\ds
\left\| f\ast\Lambda^{\alpha\over 2}\right\|_{\L^s(\R^{n+1})}~\leq~\C_{p~s}~e^{\c|\Im\alpha|}~\left\| f\right\|_{\L^p(\R^{n+1})},
\\\\ \ds
 {\Re\alpha\over 2n}~=~{1\over p}-{1\over s},\qquad {n-1\over 2n}+{\Re\alpha\over n}~<~{1\over p}~<~{n+1\over 2n}
\end{array}
\eeq
and
\bel{mapping2}
\begin{array}{cc}\ds
\left\| f\ast\Lambda^{\alpha\over 2}\right\|_{\L^q(\R^{n+1})}~\leq~\C_{s~q}~e^{\c|\Im\alpha|}~\left\| f\right\|_{\L^s(\R^{n+1})},
\\\\ \ds
{\Re\alpha\over 2n}~=~{1\over s}-{1\over q},\qquad{n-1\over 2n}+{\Re\alpha\over 2n}~<~{1\over s}~<~{n+1\over 2n}-{\Re\alpha\over 2n}.
\end{array}
\eeq
From (\ref{mapping1})-(\ref{mapping2}), we conclude
\bel{REGULARITY double}
\begin{array}{cc}\ds
\left\| f\ast\Lambda^{\alpha\over 2}\ast\Lambda^{\alpha\over 2}\right\|_{\L^q(\R^{n+1})}~\leq~\C_{p~q}~e^{\hbox{\small{{\bf c}}}|\Im\alpha|}~\left\| f\right\|_{\L^p(\R^{n+1})},\qquad \hbox{\small{$0<\Re\alpha<1$}},
\\\\ \ds
 {\Re\alpha\over n}~=~{1\over p}-{1\over q},\qquad {n-1\over 2n}+{\Re\alpha\over n}~<~{1\over p}~<~{n+1\over 2n}.
\end{array}
\eeq
Hence that $\c_{\alpha/2}^2\Gamma^{-2}\left(\delta\left({\alpha\over 2}\right)\right)\Big[\Hat{\Lambda}^{\alpha\over 2}(\xi,\tau)\Big]^2$ 
is another $\L^p\mt\L^q$-Fourier multiplier with  $p,q$, $\Re\alpha$ satisfying  (\ref{Formula A}). 

By solving the two  equations  in (\ref{Lambda transform UV}) and (\ref{Multiplier}), we find 
\bel{V rewrite}
\begin{array}{lr}
\c_\alpha\Gamma^{-1}(\delta(\alpha))\Hat{\Lambda}^\alpha(\xi,\tau)- \c_{\alpha/2}^2\Gamma^{-2}\left(\delta\left({\alpha\over 2}\right)\right)\Big[\Hat{\Lambda}^{\alpha\over 2}(\xi,\tau)\Big]^2~=~
\\\\
\Gamma\left(1-\delta(\alpha)\right)
\Bigg\{ -\sin\pi\left(\delta(\alpha)-{1\over 2}\right)- \sin^2\pi\left[\delta\left({\alpha\over 2}\right)-{1\over 2}\right]\Bigg\}\hbox{V}^\alpha(\xi,\tau),\qquad  \hbox{\small{$|\xi|\neq|\tau|$}}.
\end{array}
\eeq
A direct computation shows
\bel{sine computation}
\begin{array}{lr}\ds
\hbox{$\sin\pi\left(\delta(\alpha)-{1\over 2}\right)$}~=~
\\ \ds
-{1\over 2} \left[ e^{-\pi\Im\delta(\alpha)}+e^{\pi\Im\delta(\alpha)}\right]\cos \pi\Re\delta(\alpha)-{\i\over 2} \left[ e^{-\pi\Im\delta(\alpha)}-e^{\pi\Im\delta(\alpha)}\right]\sin\pi\Re\delta(\alpha),
\\\\ \ds
\hbox{$\sin^2\pi\left[\delta\left({\alpha\over 2}\right)-{1\over 2}\right]$}~=~
\\ \ds
{1\over 2}+{1\over 4} \left[ e^{-\pi\Im\delta(\alpha)}+e^{\pi\Im\delta(\alpha)}\right]\cos \pi\Re\delta(\alpha)+{\i\over 4} \left[ e^{-\pi\Im\delta(\alpha)}-e^{\pi\Im\delta(\alpha)}\right]\sin\pi\Re\delta(\alpha).
\end{array}
\eeq
For $0<\Re\delta(\alpha)<1$, (\ref{sine computation}) implies
\bel{computa not 0}
\hbox{$\sin\pi\left(\delta(\alpha)-{1\over 2}\right)$}+\hbox{$\sin^2\pi\left[\delta\left({\alpha\over 2}\right)-{1\over 2}\right]$}~\neq~0.
\eeq
Consequently,  $\hbox{V}^\alpha(\xi,\tau)$ for $0<\Re\alpha<1$  is a $\L^p\mt\L^q$-Fourier multiplier with  $p,q$, $\Re\alpha$ satisfying  (\ref{Formula A}).

Let $\S^{-\alpha}$ defined in (\ref{S^alpha}). Consider
$\W^{-\alpha} f=\Gamma^{-1}\left(\lambda(\alpha)\right)\S^{-\alpha} f$.
We have
\bel{regularity W}
\begin{array}{cc}\ds
\left\| \W^{-\alpha} f\right\|_{\L^q(\R^{n+1})}~\leq~\C_{p~q}~e^{\hbox{\small{{\bf c}}}|\Im\alpha|}~\left\| f\right\|_{\L^p(\R^{n+1})},\qquad 0<\Re\alpha<1,
\\\\ \ds
 {\Re\alpha\over n}~=~{1\over p}-{1\over q},\qquad {n-1\over 2n}+{\Re\alpha\over n}~<~{1\over p}~<~{n+1\over 2n}.
\end{array}
\eeq
The formula derived below (\ref{V inverse transform}) shows
\bel{W inverse Transform <1}
\begin{array}{lr}
\Gamma^{-1}\left(\lambda(\alpha)\right)\hbox{V}^\alpha(x,t)~=~
\pi^{\delta(\alpha)-\lambda(\alpha)-1}\left\{ {\ds\left({1\over t^2-|x|^2}\right)^{\lambda(\alpha)}_-} -\sin\pi\left(\lambda(\alpha)-{1\over 2}\right){\ds\left({1\over t^2-|x|^2}\right)^{\lambda(\alpha)}_+}\right\},
\\\\ \ds~~~~~~~~~~~~~~~~~~~~~~~~~~~~~~~~~~~~~~~~~~~~~~~~~~~~~~~~~~~~~~~~~~~~~~~~~~~~~~~~~~~~~~~~~~~~
\hbox{\small{ $|x|\neq|t|$ }},\qquad\hbox{\small{$0\leq\Re\lambda(\alpha)<1$}}.
\end{array}
\eeq
Clearly,  $\Gamma^{-1}\left(\lambda(\alpha)\right)\hbox{V}^\alpha\in\L^\infty(\R^{n+1})$ at
 $\Re\lambda(\alpha)=0$. This implies
 \bel{Result W n}
\left\| \W^{-\alpha} f\right\|_{\L^\infty(\R^{n+1})}~\leq~\C ~e^{\hbox{\small{{\bf c}}}|\Im\alpha|}~\left\| f\right\|_{\L^1(\R^{n+1})},\qquad \hbox{\small{$\Re\alpha=n$}}.
\eeq
\begin{remark} $\Gamma^{-1}\left(\lambda(\alpha)\right)\hbox{V}^\alpha$ is analytic for $0<\Re\alpha<n$. See appendix B.
\end{remark} 

From (\ref{regularity W}) and (\ref{Result W n}), by applying Stein interpolation theorem \cite{Stein}, we conclude that
$\S^{-\alpha} f=\Gamma(\lambda(\alpha))  \W^{-\alpha} f$ satisfies (\ref{Result Two}) for $0<\alpha<n$.

 \section{End-point estimates regarding $\II^\alpha_j$ and $\III^\alpha_j$}
 \setcounter{equation}{0}
 First, we repeat certain definitions  from Section 1.1.  Recall 
 \[\hbox{$\lambda(\alpha)={n+1\over 2}-\delta(\alpha)$},\qquad \hbox{$\delta(\alpha)=\left({n+1\over 2n}\right)\alpha$},\qquad\hbox{$\alpha\in\Cx$}.\] 
$\Omega^\alpha$ is a distribution  defined in $\R^n$ by analytic continuation from
\bel{Omega^alpha again}
\begin{array}{cc}\ds
\Omega^{\alpha}(x)~=~\pi^{-\lambda(\alpha)}\Gamma^{-1}\left(1-\lambda(\alpha)\right) 
 \left({1 \over1-|x|^2}\right)^{\lambda(\alpha)}_+,\qquad \hbox{\small{$\Re\lambda(\alpha)<1$}}. 
\end{array}
\eeq
Equivalently,  it can be  defined by its Fourier transform
 \bel{Omega^alpha transform} 
\begin{array}{lr}\ds
\Hat{\Omega}^{\alpha}(\xi)~=~\left({1\over|\xi|}\right)^{{n\over 2}-\lambda(\alpha)} \J_{{n\over 2}-\lambda(\alpha)}\Big(2\pi|\xi|\Big)
\\\\ \ds~~~~~~~~~~
~=~\left({1\over|\xi|}\right)^{\delta(\alpha)-{1\over 2}} \J_{\delta(\alpha)-{1\over 2}}\Big(2\pi|\xi|\Big),\qquad \alpha\in\Cx
\end{array}
\eeq
where $\J$ is  Bessel function.

Given $j>0$,  $\left\{ \xi^\nu_j\right\}_\nu$ is a collection of points  almost equally distributed on  $\mathds{S}^{n-1}$ having a grid length between $2^{-j/2-1}$ and $2^{-j/2}$:   

{\bf (1)} $|\xi^\mu_j-\xi^\nu_j|\ge 2^{-j/2-1}$ for every $\xi^\mu_j\neq\xi^\nu_j$. 

{\bf (2)} For any $u\in \mathds{S}^{n-1}$, there is a $\xi^\nu_j$ in the open set $\{\xi\in\mathds{S}^{n-1}\colon |\xi-u|<2^{-j/2+1}\}$.

Let $\varphi\in\mathcal{C}^\infty_o(\R)$  be a smooth cut-off function such that $\varphi(r)=1$ if $|r|\leq1$ and $\varphi(r)=0$ for $|r|>2$. Define
\bel{phi^v_j}
\varphi^\nu_j(\xi)~=~{\ds\varphi\Bigg[2^{j/2}\left|{\xi\over |\xi|}-\xi^{\nu}_{j}\right|\Bigg]\over \ds \sum_{\xi^\nu_j\in\mathds{S}^{n-1}} \varphi\Bigg[2^{j/2}\left|{\xi\over |\xi|}-\xi^{\nu}_{j}\right|\Bigg]}
\eeq
whose support  is contained in the cone
\bel{Gamma}
\Gamma^\nu_j~=~\Bigg\{\xi\in\R^n\colon \left|{\xi\over |\xi|}-\xi^{\nu}_{j}\right|<2^{-j/2+1}\Bigg\}.
\eeq
Given $\xi_j^\nu$, $\L_\nu$ is an $n\times n$-orthogonal matrix with $\det\L_\nu=1$ and 
$\L_\nu^T \xi^\nu_j=(1,0)^T\in\R\times\R^{n-1}$.
Define
 \bel{Psi nu j alpha again}
 \begin{array}{cc}\ds
\Psi^{\alpha}_j(x)~=~\sum_{\nu\colon\xi^\nu_j\in\mathds{S}^{n-1}}  \Psi^{\alpha~\nu}_j(x),
 \\\\ \ds 
 \Psi^{\alpha~\nu}_j(x)~=~\varphi\left[ 2^{\big(1-\sigma\big)j}\left|1-\left(\L_\nu^T x\right)_1\right|\right] 
\prod_{i=2}^n \varphi\left[2^{\big[{1\over 2}-\sigma\big]j}\left|\left(\L_\nu^T x\right)_i\right|\right]
\\ \ds~~~~~~~~~~~~~
~+~\varphi\left[ 2^{\big(1-\sigma\big)j}\left|1+\left(\L_\nu^T x\right)_1\right|\right] 
\prod_{i=2}^n \varphi\left[2^{\big[{1\over 2}-\sigma\big]j}\left|\left(\L_\nu^T x\right)_i\right|\right].
\end{array}
\eeq
Each $\Psi^{\alpha}_j$ has a support inside the rectangles
\bel{Rectangle'}
\begin{array}{cc}\ds
\hbox{R}^\nu_{j}~=~\hbox{R}^{\nu~+}_{j}\cup \hbox{R}^{\nu~-}_{j},
\\\\ \ds
\hbox{R}^{\nu~+}_{j}~=~\Bigg\{x\in\R^n\colon1-2^{\big(\sigma-1\big)j+1}<\left(\L_\nu^T x\right)_1<1+2^{\big(\sigma-1\big)j+1},~\left|\left(\L_\nu^T x\right)_i\right|<2^{\big[\sigma-{1\over 2}\big]j+1},~ i=2,\ldots,n
\Bigg\},
\\\\ \ds
\hbox{R}^{\nu~-}_{j}~=~ \Bigg\{x\in\R^n\colon-1-2^{\big(\sigma-1\big)j+1}<\left(\L_\nu^T x\right)_1<-1+2^{\big(\sigma-1\big)j+1},~\left|\left(\L_\nu^T x\right)_i\right|<2^{\big[\sigma-{1\over 2}\big]j+1},~ i=2,\ldots,n
\Bigg\}.
\end{array}
\eeq
Recall $\phi_j$ and $\Omega^\alpha_j$ defined in (\ref{varphi_j})-(\ref{Omega alpha j}). 
Let $\varphi^\nu_j$ defined in (\ref{phi^v_j}). Consider
\bel{P v j}
 \begin{array}{cc}\ds
 \Hat{\Omega}_j^\alpha(\xi)~=~
\sum_{\nu\colon \xi^\nu_j\in\mathds{S}^{n-1}}  \Hat{\Omega}_j^{\alpha~\nu}(\xi),
\qquad
 \Hat{\Omega}_j^{\alpha~\nu}(\xi)~=~\varphi^\nu_j(\xi) \Hat{\Omega}^\alpha_j (\xi),\qquad j>0. 
  \end{array}
\eeq
Let $\Psi^{\alpha}_j$ and $\Psi^{\alpha~\nu}_j$ defined in (\ref{Psi nu j alpha again}).
For every $j>0$, we assert
   \bel{II_j U}
 \begin{array}{cc}\ds
 \II^\alpha_jf(x,t)~=~\iint_{\R^{n+1}} f(x-u,t-r) \U^\alpha_{j~|r|}(u) |r|^{\big({n+1\over n}\big)\alpha-1}dudr 
 \\\\ \ds
 \U^\alpha_j(x) ~=~
 \sum_{\nu\colon\xi^\nu_j\in\mathds{S}^{n-1}}  \U^{\alpha~\nu}_j(x), 
\qquad
  \U^{\alpha~\nu}_j(x)~=~\Omega^{\alpha~\nu}_j(x) \left[1-\Psi^{\alpha~\nu}_j(x)\right] 
    \end{array}
  \eeq 
 and
 \bel{III_j V}
 \begin{array}{cc}\ds
 \III^\alpha_jf(x,t)~=~\iint_{\R^{n+1}} f(x-u,t-r) \V^\alpha_{j~|r|}(u)|r|^{\big({n+1\over n}\big)\alpha-1}dudr
 \\\\ \ds
 \V^\alpha_j(x) ~=~
 \sum_{\nu\colon\xi^\nu_j\in\mathds{S}^{n-1}}  \V^{\alpha~\nu}_j(x), 
\qquad
  \V^{\alpha~\nu}_j(x)~=~\Omega^{\alpha~\nu}_j(x) \Psi^{\alpha~\nu}_j(x) 
    \end{array}
 \eeq 
 where $\U^\alpha_{j~r}(x)=r^{-n}  \U^\alpha_j\big({x\over r}\big),r>0$ and  $\V^\alpha_{j~r}(x)=r^{-n}  \V^\alpha_j\big({x\over r}\big), r>0$.

\subsection{${^\sharp}\II^{\alpha~\beta}_j$ and ${^\sharp}\III^{\alpha~\beta}_j$}
${^\sharp}\Omega^\alpha$ is a distribution defined in $\R^n$ by analytic continuation from
\bel{Omega sharp}
\begin{array}{cc}\ds
 {^\sharp}\Omega^{\alpha}(x)~=~\pi^{-\lambda(\alpha)-{1\over 2}}\hbox{$\Gamma^{-1}\left({1\over 2}-\lambda(\alpha)\right)$} 
 \left({1 \over1-|x|^2}\right)^{\lambda(\alpha)+{1\over 2}}_+,\qquad \hbox{\small{$\Re\lambda(\alpha)<{1\over 2}$}}. 
 \end{array}
\eeq
Equivalently,  it can be  defined by
 \bel{Omega sharp Transform} 
\begin{array}{lr}\ds
{^\sharp}\Hat{\Omega}^{\alpha}(\xi)~=~
\left({1\over|\xi|}\right)^{{n\over 2}-\big[\lambda(\alpha)+{1\over 2}\big]} \J_{{n\over 2}-\big[\lambda(\alpha)+{1\over 2}\big]}\Big(2\pi|\xi|\Big)
\\\\ \ds~~~~~~~~~~~
~=~\left({1\over|\xi|}\right)^{\delta(\alpha)-1} \J_{\delta(\alpha)-1}\Big(2\pi|\xi|\Big),\qquad \alpha\in\Cx.
\end{array}
\eeq
Recall $\phi_j$ defined in (\ref{varphi_j}). Let
\bel{Omega_j sharp}
 \begin{array}{lr}\ds
 {^\sharp}\Hat{\Omega}^\alpha_j(\xi)~=~{^\sharp}\Hat{\Omega}^{\alpha}(\xi) \phi_j(\xi),\qquad j>0.
  \end{array}
\eeq
For $\alpha,\beta\in\Cx$ and $0<\Re\alpha=\Re\beta<1$, we define
\bel{Lambda sharp}
\begin{array}{cc}\ds
{^\sharp}\I^{\alpha~\beta}_j f(x,t)~=~\iint_{\R^{n+1}}  f(x-u,t-r) {^\sharp}\Lambda^{\alpha~\beta}_j(u,r) dudr,
\\\\ \ds
{^\sharp}\Hat{\Lambda}^{\alpha~\beta}_j(\xi,\tau)~=~\int_\R e^{-2\pi\i \tau r}  {^\sharp}\Hat{\Omega}^\alpha_j(|r|\xi) |r|^{{\beta\over 2}-1}dr.
\end{array}
\eeq
\begin{lemma} We have
\bel{Transform T Est}
\begin{array}{lr}\ds
\sum_{j>0} \left|{^\sharp}\Hat{\Lambda}_j^{\alpha~\beta}(\xi,\tau)\right|~\leq~\C_{\Re \alpha} ~e^{\c|\Im \alpha|}e^{\c|\Im\beta|}~|\xi|^{-{\Re\alpha\over 2}} \left|1-{|\tau|\over |\xi|}\right|^{-{1\over 2}+\ve}  
\end{array}
\eeq
for some $\ve=\ve(\Re\alpha)>0$.
\end{lemma}
Let $\varphi^\nu_j$ defined in (\ref{phi^v_j}). Consider
\bel{P v j sharp}
 \begin{array}{cc}\ds
 {^\sharp}\Hat{\Omega}_j^\alpha(\xi)~=~
\sum_{\nu\colon \xi^\nu_j\in\mathds{S}^{n-1}}  {^\sharp}\Hat{\Omega}_j^{\alpha~\nu}(\xi),
\qquad
 {^\sharp}\Hat{\Omega}_j^{\alpha~\nu}(\xi)~=~\varphi^\nu_j(\xi) {^\sharp}\Hat{\Omega}^\alpha_j (\xi),\qquad j>0. 
  \end{array}
\eeq
Let $\Psi^{\alpha}_j$ and $\Psi^{\alpha~\nu}_j$ defined in (\ref{Psi nu j alpha again}).
For every $j>0$, we assert
   \bel{II_j U sharp}
 \begin{array}{cc}\ds
 {^\sharp}\II^{\alpha~\beta}_jf(x,t)~=~\iint_{\R^{n+1}} f(x-u,t-r) {^\sharp}\U^\alpha_{j~|r|}(u) |r|^{{\beta\over 2}-1}dudr 
 \\\\ \ds
 {^\sharp}\U^\alpha_j(x) ~=~
 \sum_{\nu\colon\xi^\nu_j\in\mathds{S}^{n-1}}  {^\sharp}\U^{\alpha~\nu}_j(x), 
\qquad
  {^\sharp}\U^{\alpha~\nu}_j(x)~=~{^\sharp}\Omega^{\alpha~\nu}_j(x) \left[1-\Psi^{\alpha~\nu}_j(x)\right], 
    \end{array}
  \eeq 
 \bel{III_j V sharp}
 \begin{array}{cc}\ds
{^\sharp}\III^{\alpha~\beta}_jf(x,t)~=~\iint_{\R^{n+1}} f(x-u,t-r) {^\sharp}\V^\alpha_{j~|r|}(u)|r|^{{\beta\over 2}-1}dudr
 \\\\ \ds
 {^\sharp}\V^\alpha_j(x) ~=~
 \sum_{\nu\colon\xi^\nu_j\in\mathds{S}^{n-1}}  {^\sharp}\V^{\alpha~\nu}_j(x), 
\qquad
  {^\sharp}\V^{\alpha~\nu}_j(x)~=~{^\sharp}\Omega^{\alpha~\nu}_j(x) \Psi^{\alpha~\nu}_j(x) 
    \end{array}
 \eeq 
 where ${^\sharp}\U^\alpha_{j~r}(x)=r^{-n}  {^\sharp}\U^\alpha_j\big({x\over r}\big),r>0$ and  ${^\sharp}\V^\alpha_{j~r}(x)=r^{-n}  {^\sharp}\V^\alpha_j\big({x\over r}\big), r>0$.

From (\ref{II_j U sharp})-(\ref{III_j V sharp}), we find
\bel{II sharp transform}
\begin{array}{cc}\ds
{^\sharp}\Hat{\II^{\alpha~\beta}_j f}(\xi,\tau)~=~\Hat{f}(\xi,\tau) {^\sharp}\P^{\alpha~\beta}_j(\xi,\tau),
\qquad
{^\sharp}\P^{\alpha~\beta}_j(\xi,\tau)~=~\int_\R e^{-2\pi\i \tau r} {^\sharp}\Hat{\U}^\alpha_j(|r|\xi) |r|^{{\beta\over 2}-1}dr
\end{array}
\eeq
and
\bel{III sharp transform}
\begin{array}{cc}\ds
{^\sharp}\Hat{\III^{\alpha~\beta}_j f}(\xi,\tau)~=~\Hat{f}(\xi,\tau) {^\sharp}\Q^{\alpha~\beta}_j(\xi,\tau),
\qquad
{^\sharp}\Q^{\alpha~\beta}_j(\xi,\tau)~=~\int_\R e^{-2\pi\i \tau r} {^\sharp}\Hat{\V}^\alpha_j(|r|\xi) |r|^{{\beta\over 2}-1}dr.
\end{array}
\eeq

{\bf Proposition One}~~{\it Let $0< \Re\alpha=\Re\beta<1$. We have
\bel{Result One U}
\left|{^\sharp}\P^{\alpha~\beta}_j(\xi,\tau)
\right|~\leq~\C_{\Re \alpha} ~e^{\c|\Im \alpha|}e^{\c|\Im\beta|}~  |\xi|^{-{\Re\alpha\over 2}} 2^{j/2} 2^{-\ve j},\qquad\xi\neq0
\eeq
and
\bel{Result One V}
\left| {^\sharp}\Q^{\alpha~\beta}_j(\xi,\tau)\right|~\leq~\C_{\Re \alpha} ~e^{\c|\Im \alpha|}e^{\c|\Im\beta|} ~|\xi|^{-{\Re\alpha\over2}}2^{-\ve j},\qquad \xi\neq0
\eeq
for every $j>0$ and some $\ve=\ve(\Re\alpha)>0$.}

\subsection{${^\flat}\II^{\alpha~\beta}_j$ and ${^\flat}\III^{\alpha~\beta}_j$}
${^\flat}\Omega^\alpha$ is a distribution defined in $\R^n$ by analytic continuation from
\bel{Omega flat}
\begin{array}{cc}\ds
 {^\flat}\Omega^{\alpha}(x)~=~\pi^{\delta(\alpha)-1}\Gamma^{-1}\left(\delta(\alpha)\right) 
 \left({1 \over1-|x|^2}\right)^{1-\delta(\alpha)}_+,\qquad \hbox{\small{$\Re\delta(\alpha)>0$}}. 
 \end{array}
\eeq
Equivalently, it can be defined by
 \bel{Omega flat Transform} 
\begin{array}{lr}\ds
{^\flat}\Hat{\Omega}^{\alpha}(\xi)~=~\left({1\over|\xi|}\right)^{{n\over 2}-1+\delta(\alpha)} \J_{{n\over 2}-1+\delta(\alpha)}\Big(2\pi|\xi|\Big)
\\\\ \ds~~~~~~~~~~~
~=~\left({1\over|\xi|}\right)^{{n-1\over 2}+\delta(\alpha)-{1\over2}} \J_{{n-1\over 2}+\delta(\alpha)-{1\over2}}\Big(2\pi|\xi|\Big),\qquad \alpha\in\Cx.
\end{array}
\eeq
\begin{lemma} Let ${^\flat}\Omega^\alpha$ defined in (\ref{Omega flat}) for $0<\Re\delta(\alpha)<1$. We have
\bel{Result 3.2}
\begin{array}{cc}\ds
\left\| f\ast {^\flat}\Omega^\alpha\right\|_{\L^q(\R^n)}~\leq~\C_{\Re\alpha~p}~e^{\c|\Im\alpha|}~\left\| f\right\|_{\L^p(\R^n)},\qquad 1<p<q<\infty
\\\\ \ds
\hbox{if}\qquad \Re\delta(\alpha)~=~{1\over p}-{1\over q}.
\end{array}
\eeq
\end{lemma}
Recall $\phi_j$ defined in (\ref{varphi_j}). Let
\bel{Omega_j flat}
 \begin{array}{lr}\ds
 {^\flat}\Hat{\Omega}^\alpha_j(\xi)~=~{^\flat}\Hat{\Omega}^{\alpha}(\xi) \phi_j(\xi),\qquad j>0.
  \end{array}
\eeq
Note that every $\phi_j, j>0$ is an $\L^p$-multiplier for $1<p<\infty$. Hence, {\bf Lemma 3.2} implies
\bel{Result flat Omega_j}
\begin{array}{cc}\ds
\left\| f\ast {^\flat}\Omega^\alpha_j\right\|_{\L^q(\R^n)}~\leq~\C_{\Re\alpha~p}~e^{\c|\Im\alpha|}~\left\| f\right\|_{\L^p(\R^n)},\qquad j>0
\\\\ \ds
\hbox{for}\qquad \Re\delta(\alpha)~=~{1\over p}-{1\over q},\qquad 1<p<q<\infty.
\end{array}
\eeq
Let $\alpha,\beta\in\Cx$ and $0<\Re\alpha=\Re\beta<1$. We define
\bel{I flat}
{^\flat}\I^{\alpha~\beta}_j f(x,t)~=~\iint_{\R^{n+1}}  f(x-u,t-r) {^\flat}\Omega^\alpha_{j~|r|}(u) |r|^{\big({n+3\over 2}\big)\beta-1}dudr
\eeq
where
${^\flat}\Omega^\alpha_{j~r}(x)=r^{-n} {^\flat}\Omega^\alpha_j\left({x\over r}\right), r>0$.

Recall $\delta(\alpha)=\big({n+1\over 2n}\big)\alpha$. Let 
\[\Re\delta(\alpha)~=~{1\over p}-{1\over q},\qquad \Re\alpha~=~{1\over p}-{1\over s},\qquad 1<p<q\leq s<\infty.\] 
By using (\ref{Result flat Omega_j}), we have
\bel{I j flat EST}
\begin{array}{lr}\ds
\left\{\int_\R \left\{\int_{\R^n} \left| {^\flat}\I^{\alpha~\beta}_j f(x,t)\right|^q dx\right\}^{s\over q}dt\right\}^{1\over s}
\\\\ \ds
\leq~\left\{\int_\R \left\{\int_\R \left\{\int_{\R^n} \left| \int_{\R^n}  f(x-u,t-r) {^\flat}\Omega^\alpha_{j~|r|}(u) du\right|^q dx\right\}^{1\over q} |r|^{\big({n+3\over 2}\big)\Re\alpha-1}dr\right\}^s dt\right\}^{1\over s}
\\ \ds~~~~~~~~~~~~~~~~~~~~~~~~~~~~~~~~~~~~~~~~~~~~~~~~~~~~~~~~~~~~~~~~~~~~~~
\hbox{\small{by Minkowski integral inequality}}
\\\\ \ds
\leq~\C_{\Re\alpha~p}~e^{\c|\Im\alpha|} \left\{\int_\R \left\{\int_\R  \left\| f(\cdot, t-r)\right\|_{\L^p(\R^n)} |r|^{-n\Re\delta(\alpha)} |r|^{\big({n+3\over 2}\big)\Re\alpha-1}dr\right\}^s dt\right\}^{1\over s}\qquad\hbox{\small{by (\ref{Result flat Omega_j})}}
\\\\ \ds
=~\C_{\Re\alpha~p}~e^{\c|\Im\alpha|}\left\{\int_\R \left\{\int_\R  \left\| f(\cdot, t-r)\right\|_{\L^p(\R^n)}  |r|^{\Re\alpha-1}dr\right\}^s dt\right\}^{1\over s}
\\\\ \ds
\leq~\C_{\Re\alpha~p}~e^{\c|\Im\alpha|}\left\| f\right\|_{\L^p(\R^{n+1})}\qquad\hbox{\small{by Hardy-Littlewood-Sobolev inequality \cite{Hardy-Littlewood}-\cite{Sobolev}}}.
\end{array}
\eeq
Let $\varphi^\nu_j$ defined in (\ref{phi^v_j}). Consider 
\bel{P v j flat}
 \begin{array}{cc}\ds
 {^\flat}\Hat{\Omega}_j^\alpha(\xi)~=~
\sum_{ \xi^\nu_j\in\mathds{S}^{n-1}}  {^\flat}\Hat{\Omega}_j^{\alpha~\nu}(\xi),
\qquad
 {^\flat}\Hat{\Omega}_j^{\alpha~\nu}(\xi)~=~\varphi^\nu_j(\xi) {^\flat}\Hat{\Omega}^\alpha_j (\xi),\qquad j>0.
  \end{array}
\eeq
 Let $\Psi^{\alpha}_j$, $\Psi^{\alpha~\nu}_j$ defined in (\ref{Psi nu j alpha again}). For every $j>0$, we assert
 \bel{II_j U flat}
 \begin{array}{cc}\ds
 {^\flat}\II^{\alpha~\beta}_jf(x,t)~=~\iint_{\R^{n+1}} f(x-u,t-r) {^\flat}\U^\alpha_{j~|r|}(u) |r|^{\big({n+3\over 2}\big)\beta-1}dudr 
 \\\\ \ds
 {^\flat}\U^\alpha_j(x) ~=~
 \sum_{\xi^\nu_j\in\mathds{S}^{n-1}}  {^\flat}\U^{\alpha~\nu}_j(x), 
\qquad
  {^\flat}\U^{\alpha~\nu}_j(x)~=~{^\flat}\Omega^{\alpha~\nu}_j(x) \left[1-\Psi^{\alpha~\nu}_j(x)\right] 
    \end{array}
  \eeq 
and
 \bel{III_j V flat}
 \begin{array}{cc}\ds
{^\flat}\III^{\alpha~\beta}_jf(x,t)~=~\iint_{\R^{n+1}} f(x-u,t-r) {^\flat}\V^\alpha_{j~|r|}(u)|r|^{\big({n+3\over 2}\big)\beta-1}dudr
 \\\\ \ds
 {^\flat}\V^\alpha_j(x) ~=~
 \sum_{\xi^\nu_j\in\mathds{S}^{n-1}}  {^\flat}\V^{\alpha~\nu}_j(x), 
\qquad
  {^\flat}\V^{\alpha~\nu}_j(x)~=~{^\flat}\Omega^{\alpha~\nu}_j(x) \Psi^{\alpha~\nu}_j(x) 
    \end{array}
 \eeq 
 where ${^\flat}\U^\alpha_{j~r}(x)=r^{-n}  {^\flat}\U^\alpha_j\big({x\over r}\big),r>0$ and  ${^\flat}\V^\alpha_{j~r}(x)=r^{-n}  {^\flat}\V^\alpha_j\big({x\over r}\big), r>0$. 
\v

{\bf Proposition Two}~~{\it Let $0<\Re\alpha<1$. For every $j>0$, we have
\bel{Result Two U}
\begin{array}{cc}\ds
\left\| f\ast {^\flat}\U^\alpha_j\right\|_{\L^q(\R^n)}~\leq~\C_{\Re\alpha~p} ~e^{\c|\Im\alpha|}~2^{-N j} \left\| f\right\|_{\L^p(\R^n)}, \qquad N\ge0
\\\\ \ds
\hbox{if}\qquad \Re\delta(\alpha)~=~{1\over p}-{1\over q},\qquad 1<p<q<\infty
\end{array}
\eeq
and
\bel{Result Two V}
\begin{array}{cc}\ds
\left\| f\ast {^\flat}\V^\alpha_j\right\|_{\L^q(\R^n)}~\leq~\C_{\Re\alpha~p} ~e^{\c|\Im\alpha|} ~\left\| f\right\|_{\L^p(\R^n)}
\\\\ \ds
\hbox{if}\qquad \Re\delta(\alpha)~=~{1\over p}-{1\over q},\qquad 1<p<q<\infty.
\end{array}
\eeq
}

Let $0<\Re\alpha=\Re\beta<1$. By applying (\ref{Result Two U})-(\ref{Result Two V}) and carrying out an iteration estimate as in (\ref{I j flat EST}), we obtain
\bel{Result Two II}
\begin{array}{cc}\ds
\left\{\int_\R \left\{\int_{\R^n} \left| {^\flat}\II^{\alpha~\beta}_j f(x,t)\right|^q dx\right\}^{s\over q}dt\right\}^{1\over s}
~\leq~\C_{\Re\alpha~p} ~e^{\c|\Im\alpha|}~2^{-N j} \left\| f\right\|_{\L^p(\R^{n+1})}, \qquad N\ge0
\\\\ \ds
\hbox{if}\qquad \Re\delta(\alpha)~=~{1\over p}-{1\over q},\qquad \Re\alpha~=~{1\over p}-{1\over s},\qquad  1<p<q\leq s<\infty
\end{array}
\eeq
and
\bel{Result Two III}
\begin{array}{cc}\ds
\left\{\int_\R \left\{\int_{\R^n} \left| {^\flat}\III^{\alpha~\beta}_j f(x,t)\right|^q dx\right\}^{s\over q}dt\right\}^{1\over s}
~\leq~\C_{\Re\alpha~p} ~e^{\c|\Im\alpha|}~ \left\| f\right\|_{\L^p(\R^{n+1})}
\\\\ \ds
\hbox{if}\qquad \Re\delta(\alpha)~=~{1\over p}-{1\over q},\qquad \Re\alpha~=~{1\over p}-{1\over s},\qquad  1<p<q\leq s<\infty.
\end{array}
\eeq

\section{Proof of Theorem A}
\setcounter{equation}{0}
Recall $\I^\alpha f$ defined in (\ref{I_alpha rewrite}) and $\I^\alpha_<$, $\I^\alpha_j, j>0$  defined in (\ref{I_alpha <}) and (\ref{I_alpha j})-(\ref{I rewrite Transform}). 
We have
\[ \I^\alpha f(x,t)
~=~\I^\alpha_< f(x,t)+\sum_{j>0} \I^\alpha_j f(x,t)\]
of which
\bel{I rewrite again}
 \begin{array}{cc}\ds
 \I^\alpha_j f(x,t)
~=~\iint_{\R^{n+1}} e^{2\pi\i \big[x\cdot\xi+t\tau\big]} \Hat{f}(\xi,\tau) \Lambda^\alpha_j(\xi,\tau) d\xi d\tau,\qquad j>0,
\\\\ \ds
\Hat{\Lambda}^\alpha_j(\xi,\tau)~=~\int_\R e^{-2\pi\i r\tau} \left({1\over|r\xi|}\right)^{\delta(\alpha)-{1\over 2}} \J_{\delta(\alpha)-{1\over 2}}\Big(2\pi|r\xi|\Big)\phi_j(r\xi) |r|^{\left({n+1\over n}\right)\alpha-1} dr.
\end{array}
 \eeq
The $\L^p\mt\L^q$-boundedness of $\I^\alpha_<$ has been discussed in Section 1. 
From now on, we focus on $ \I^\alpha_{j}, j>0$.

For $0\leq\Re z\leq1$, we define
\bel{Omega^a_z}
\begin{array}{lr}\ds
\Hat{\Omega}^{\alpha~z}_j(\xi)~=~
\phi_j(\xi) \left({1\over|\xi|}\right)^{\delta(\alpha)-{1\over 2}+\left({n-1\over 2}\right)z-{1\over 2}(1-z)}
 \J_{\delta(\alpha)-{1\over 2}+\left({n-1\over 2}\right)z-{1\over 2}(1-z)}\left(2\pi|\xi|\right)
 \\\\ \ds~~~~~~~~~~~~
 ~=~\pi^{\delta(\alpha)-1+\left({n-1\over 2}\right)z-{1\over 2}(1-z)}\hbox{$\Gamma^{-1}\left(\delta(\alpha)+\left({n-1\over 2}\right)z-{1\over 2}(1-z)\right)$}
 \\ \ds~~~~~~~~~~~~~~~~~~~
\phi_j(\xi) \int_{-1}^1 e^{2\pi\i |\xi|s} (1-s^2)^{\delta(\alpha)-1+\left({n-1\over 2}\right)z-{1\over 2}(1-z)}ds\qquad\hbox{\small{by (\ref{Bessel})}}
 \end{array}
 \eeq
and
\bel{Lambda^a_z}
\Hat{\Lambda}^{\alpha~z}_j(\xi,\tau)~=~e^{\big[z-{1\over n}\big]^2}\int_{\R} e^{-2\pi\i r\tau} \Hat{\Omega}^{\alpha~z}_j (r\xi) |r|^{\big({n+3\over 2}\big)\Re\alpha z+{\Re\alpha\over 2}(1-z)+\i\Im\alpha-1} dr,\qquad\hbox{\small{$j>0$}}.
\eeq
Observe that  $\Hat{\Lambda}^{\alpha~{1\over n}}_{j}(\xi,\tau)=\Hat{\Lambda}^{\alpha}_{j}(\xi,\tau)$.
 Let $0<\Re\alpha<1$.  Define
 \bel{I^a_z,j}
 \begin{array}{lr}\ds
 \I^{\alpha~z}_{j} f(x,t)~=~e^{\big[z-{1\over n}\big]^2} \iint_{\R^{n+1}} e^{2\pi\i x\cdot\xi}\Hat{f}(\xi,t-r)  \Hat{\Omega}^{\alpha~z}_j (r\xi) |r|^{\big({n+3\over 2}\big)\Re\alpha z+{\Re\alpha\over 2}(1-z)+\i\Im\alpha-1} d\xi dr 
 \\\\ \ds~~~~~~~~~~~~~~~
 ~=~\iint_{\R^{n+1}} e^{2\pi\i \big[x\cdot\xi+t\tau\big]} \Hat{f}(\xi,\tau) \Hat{\Lambda}^{\alpha~z}_j(\xi,\tau)d\xi d\tau,\qquad j>0.
 \end{array}
\eeq
Given $j>0$, $\left\{ \xi^\nu_j\right\}_\nu$ is a collection of points almost equally distributed on  $\mathds{S}^{n-1}$ having  a grid length between $2^{-j/2-1}$ and $2^{-j/2}$. See  {\bf (1)}-{\bf(2)} above (\ref{phi^v_j intro})-(\ref{Gamma_j}).

Let $\varphi^\nu_j$ defined in (\ref{phi^v_j}).
Consider
\bel{Omega z Transform j}
\begin{array}{cc}\ds
\Hat{\Omega}^{\alpha~z}_j(\xi)~=~\sum_{\xi^\nu_j\in\mathds{S}^{n-1}} ~\Hat{\Omega}^{\alpha~\nu~z}_j(\xi),
\\\\ \ds
\Hat{\Omega}^{\alpha~\nu~z}_j(\xi)~=~
\varphi^\nu_j(\xi)\phi_j(\xi) \left({1\over|\xi|}\right)^{\delta(\alpha)-{1\over 2}+\left({n-1\over 2}\right)z-{1\over 2}(1-z)}
 \J_{\delta(\alpha)-{1\over 2}+\left({n-1\over 2}\right)z-{1\over 2}(1-z)}\left(2\pi|\xi|\right)
\end{array}
\eeq
Let $\Psi^{\alpha}_j$ and $\Psi^{\alpha~\nu}_j$ defined in (\ref{Psi nu j alpha again}).
For every $j>0$, we assert
   \bel{II_j U intro z}
 \begin{array}{cc}\ds
 \II^{\alpha~z}_jf(x,t)~=~e^{\big[z-{1\over n}\big]^2}\iint_{\R^{n+1}} f(x-u,t-r) \U^{\alpha~z}_{j~|r|}(u) |r|^{\big({n+3\over 2}\big)\Re\alpha z+{\Re\alpha\over 2}(1-z)+\i\Im\alpha-1} dudr, 
 \\\\ \ds
 \U^{\alpha~z}_j(x) ~=~
 \sum_{\xi^\nu_j\in\mathds{S}^{n-1}}  \U^{\alpha~\nu~z}_j(x), 
\qquad
  \U^{\alpha~\nu~z}_j(x)~=~\Omega^{\alpha~\nu~z}_j(x) \left[1-\Psi^{\alpha~\nu}_j(x)\right], 
    \end{array}
  \eeq 
 \bel{III_j V intro z}
 \begin{array}{cc}\ds
 \III^{\alpha~z}_jf(x,t)~=~e^{\big[z-{1\over n}\big]^2}\iint_{\R^{n+1}} f(x-u,t-r) \V^{\alpha~z}_{j~|r|}(u) |r|^{\big({n+3\over 2}\big)\Re\alpha z+{\Re\alpha\over 2}(1-z)+\i\Im\alpha-1}dudr,
 \\\\ \ds
 \V^{\alpha~z}_j(x) ~=~
 \sum_{\xi^\nu_j\in\mathds{S}^{n-1}}  \V^{\alpha~\nu~z}_j(x), 
\qquad
  \V^{\alpha~\nu~z}_j(x)~=~\Omega^{\alpha~\nu~z}_j(x) \Psi^{\alpha~\nu}_j(x) 
    \end{array}
 \eeq 
 where $\U^{\alpha~z}_{j~r}(x)=r^{-n}\U^{\alpha~z}_{j}\Big({x\over r}\big), r>0$ and  $\V^{\alpha~z}_{j~r}(x)=r^{-n}\V^{\alpha~z}_{j}\Big({x\over r}\big), r>0$.

In particular, at $z={1\over n}$, we find  
\bel{z=1/n}
\Hat{\Omega}^{\alpha~\nu~{1\over n}}_{j}(\xi)~=~\Hat{\Omega}^{\alpha~\nu}_{j}(\xi)\qquad\hbox{and}\qquad  
\II^{\alpha~{1\over n}}_{j} f(x,t)~=~\II^{\alpha}_j f(x,t), \qquad \III^{\alpha~{1\over n}}_{j} f(x,t)~=~\III^{\alpha}_j f(x,t)
\eeq
as defined in (\ref{II_j U})-(\ref{III_j V}). 
 
Recall  the Hardy-Littlewood-Sobolev inequality \cite{Hardy-Littlewood}-\cite{Sobolev} stated in (\ref{Hardy-Littlewood-Sobolev inequality}). Moreover, as a well known formula, we have
\[
\int_{\R^\N}e^{-2\pi\i x\cdot\xi} |x|^{\gamma-\N}dx~=~{\pi^{\N-\gamma\over 2}\Gamma\left({\gamma\over 2}\right)\over \pi^{\gamma\over 2}\Gamma\left({\N-\gamma\over 2}\right)} |\xi|^{-\gamma},\qquad \hbox{\small{$\xi\neq0$}},\qquad\hbox{\small{$0<\Re\gamma<\N$}}.
\]
By applying  (\ref{Result One U})-(\ref{Result One V}) in {\bf Proposition One} together with  (\ref{Hardy-Littlewood-Sobolev inequality}) and Plancherel theorem, we find
\bel{II 0 bad again}
\begin{array}{cc}\ds
\left\{\int_\R\left\{\int_{\R^n} \left| \II^{\alpha~0+\i\Im z}_j f (x,t)\right|^q dx\right\}^{2\over q} dt \right\}^{1\over 2}
 ~\leq~\C_{\Re\alpha}~e^{\c|\Im\alpha|}~2^{j/2}2^{-\ve j}~\left\| f\right\|_{\L^2(\R^{n+1})}
 \\\\ \ds
\hbox{if}\qquad  {\Re\alpha\over 2n}~=~{1\over 2}-{1\over q},\qquad 2<q<\infty
 \end{array}
\eeq
and
\bel{III 0 good again}
\begin{array}{cc}\ds
\left\{\int_\R\left\{\int_{\R^n} \left| \III^{\alpha~0+\i\Im z}_j f (x,t)\right|^q dx\right\}^{2\over q} dt \right\}^{1\over 2}
 ~\leq~\C_{\Re\alpha}~e^{\c|\Im\alpha|}~2^{-\ve j}~\left\| f\right\|_{\L^2(\R^{n+1})}
 \\\\ \ds
\hbox{if}\qquad  {\Re\alpha\over 2n}~=~{1\over 2}-{1\over q},\qquad 2<q<\infty
 \end{array} 
\eeq
for some $\ve=\ve(\Re\alpha)>0$.

By applying  (\ref{Result Two U})-(\ref{Result Two V}) in {\bf Proposition Two} together with  (\ref{Hardy-Littlewood-Sobolev inequality}), we have
\bel{II 1 good again}
\begin{array}{cc}\ds
\left\{\int_\R\left\{\int_{\R^n} \left| \II^{\alpha~1+\i\Im z}_j f (x,t)\right|^q dx\right\}^{s\over q} dt \right\}^{1\over s}
~\leq~\C_{N~\Re\alpha}~e^{\c|\Im\alpha|}~2^{-N j}~\left\| f\right\|_{\L^p(\R^{n+1})},\qquad N\ge0
 \\\\ \ds
\hbox{if}\qquad \left({n+1\over 2n}\right)\Re\alpha~=~{1\over p}-{1\over q},\qquad \Re\alpha~=~{1\over p}-{1\over s},\qquad 1<p<q\leq s<\infty
 \end{array}
\eeq
and
\bel{III 1 bad again}
\begin{array}{cc}\ds
\left\{\int_\R\left\{\int_{\R^n} \left| \III^{\alpha~1+\i\Im z}_j f (x,t)\right|^q dx\right\}^{s\over q} dt \right\}^{1\over s} ~\leq~\C_{\Re\alpha}~e^{\c|\Im\alpha|}~\left\| f\right\|_{\L^p(\R^{n+1})}
 \\\\ \ds
\hbox{if}\qquad \left({n+1\over 2n}\right)\Re\alpha~=~{1\over p}-{1\over q},\qquad \Re\alpha~=~{1\over p}-{1\over s},\qquad 1<p<q\leq s<\infty.
 \end{array} 
 \eeq
Take into account for 
\bel{index computa}
\begin{array}{cc}\ds
{1\over n}\left({n+2\over 2n}\right)\Re\alpha+{n-1\over n} {\Re\alpha\over 2n}~=~\left[{n+1\over 2n}+{n-1\over 2n}\right]{\Re\alpha\over n}
~=~{\Re\alpha\over n},
\\\\ \ds
{1\over n}\Re\alpha +{n-1\over n} {1\over 2}~=~{n-1\over 2n}+{\Re\alpha\over n},
\qquad
{1\over n}+{n-1\over n}{1\over 2}~=~{n+1\over 2n}.
\end{array}
\eeq
From (\ref{II 0 bad}) and (\ref{II 1 good}) with $N$ chosen sufficiently large, we obtain
\bel{II EST again}
 \begin{array}{cc}\ds
  \left\| \II^{\alpha}_j f\right\|_{\L^q(\R^{n+1})} ~\leq~\C_{\Re\alpha}~e^{\c|\Im\alpha|}~2^{-\ve\big({n-1\over n}\big) j}~\left\| f\right\|_{\L^p(\R^{n+1})}
  \\\\ \ds
\hbox{if}\qquad {\Re\alpha\over n}~=~{1\over p}-{1\over q},\qquad    {n-1\over 2n}+{\Re\alpha\over n}~< ~{1\over p}~<~ {n+1\over 2n}  
 \end{array}
\eeq
by applying the mixed-norm Stein interpolation theorem due to Benedek and Panzone \cite{Benedek-Panzone}. 
 
From (\ref{III 0 good}) and (\ref{III 1 bad}),  applying the mixed-norm Stein interpolation theorem again yields
\bel{III EST again}
 \begin{array}{cc}\ds
  \left\| \III^{\alpha}_j f\right\|_{\L^q(\R^{n+1})} ~\leq~\C_{\Re\alpha}~e^{\c|\Im\alpha|}~2^{-\ve\big({n-1\over n}\big) j}~\left\| f\right\|_{\L^p(\R^{n+1})}
  \\\\ \ds
\hbox{if}\qquad {\Re\alpha\over n}~=~{1\over p}-{1\over q},\qquad   {n-1\over 2n}+{\Re\alpha\over n}~<~ {1\over p}~<~ {n+1\over 2n}.
 \end{array}
\eeq
By summing over every $j>0$ and using (\ref{II EST again})-(\ref{III EST again}), we obtain the desired result in (\ref{Result A}).

\section{Proof of Proposition One}
\setcounter{equation}{0}
Recall  (\ref{Omega sharp Transform})-(\ref{Lambda sharp}).
 We find
\[
\begin{array}{cc}\ds
{^\sharp}\Hat{\Lambda}^{\alpha~\beta}_j(\xi,\tau)~=~\int_\R e^{-2\pi\i \tau r}  {^\sharp}\Hat{\Omega}^\alpha_j(|r|\xi) |r|^{{\beta\over 2}-1}dr
\\\\ \ds~~~~~~~~~~~~~~~~~~~~~~~~~~~~~~~~~~~~~~~~~~~~~~~~~~~~~~~~
~=~\int_\R e^{-2\pi\i \tau r} \left({1\over|r\xi|}\right)^{\delta(\alpha)-1} \J_{\delta(\alpha)-1}\Big(2\pi|r\xi|\Big) \phi_j(r\xi) |r|^{{\beta\over 2}-1}dr
\end{array}
\]
where $\phi_j$ is radially symmetric defined in (\ref{varphi_j}).

{\bf Lemma One}  {\it Let $0<\Re\alpha=\Re\beta<1$. For every $j>0$, we have
\bel{Result 1 > j}
\left| {^\sharp}\Hat{\Lambda}_j^{\alpha~\beta}(\xi,\tau)\right|~\leq~\C_{\Re \alpha} ~e^{\c|\Im \alpha|}~ |\xi|^{-{\Re\alpha\over 2}}2^{j/2} 2^{-\big({1\over 2n}\big)\Re\alpha j}\quad\hbox{if}\quad \left|1-{|\tau|\over|\xi|}\right|\leq 2^{-j}
\eeq
and
\bel{Result 1 > xi}
\left| {^\sharp}\Hat{\Lambda}_j^{\alpha~\beta}(\xi,\tau)\right|~\leq~\C_{\Re \alpha} ~e^{\c|\Im \alpha|}e^{\c|\Im\beta|}~ |\xi|^{-{\Re\alpha\over 2}}\left|1- {|\tau|\over |\xi|}\right|^{-{1\over 2}} 2^{-\big({1\over 2n}\big)\Re\alpha j}\quad\hbox{if}\quad \left|1- {|\tau|\over |\xi|}\right|> 2^{-j}.
\eeq
}

Consider
\bel{Sum split}
\begin{array}{lr}\ds
\sum_{j>0} \left|{^\sharp}\Hat{\Lambda}_j^{\alpha~\beta}(\xi,\tau)\right|~\leq~\sum_{j>0\colon \left|1-{|\tau|\over|\xi|}\right|\leq 2^{-j}} \left|{^\sharp}\Hat{\Lambda}_j^{\alpha~\beta}(\xi,\tau)\right|+\sum_{j>0\colon \left|1- {|\tau|\over |\xi|}\right|> 2^{-j}} \left|{^\sharp}\Hat{\Lambda}_j^{\alpha~\beta}(\xi,\tau)\right|.
\end{array}
\eeq
By applying {\bf Lemma One}, we find
\bel{T est.1}
\begin{array}{lr}\ds
\sum_{j>0\colon\left|1-{|\tau|\over|\xi|}\right|\leq 2^{-j}} \left|{^\sharp}\Hat{\Lambda}_j^{\alpha~\beta}(\xi,\tau)\right|
~\leq~\C_{\Re\alpha}~e^{\c\Im\alpha} ~ |\xi|^{-{\Re\alpha\over 2}}\sum_{j>0\colon \left|1-{|\tau|\over|\xi|}\right|\leq 2^{-j}}  2^{j/2} 2^{-\big({1\over 2n}\big)\Re\alpha j}\qquad\hbox{\small{by (\ref{Result 1 > j})}}
\\\\ \ds~~~~~~~~~~~~~~~~~~~~~~~~~~~~~~~~~~~~~~~~
~\leq~\C_{\Re\alpha}~e^{\c\Im\alpha}~ |\xi|^{-{\Re\alpha\over 2}}\left|1- {|\tau|\over |\xi|}\right|^{-{1\over 2}+\big({1\over 4n}\big)\Re\alpha} ~\sum_{j>0} 2^{-\big({1\over 4n}\big)\Re\alpha j}
\\\\ \ds~~~~~~~~~~~~~~~~~~~~~~~~~~~~~~~~~~~~~~~~
~\leq~\C_{\Re\alpha}~e^{\c\Im\alpha} ~ |\xi|^{-{\Re\alpha\over 2}}\left|1- {|\tau|\over |\xi|}\right|^{-{1\over 2}+\big({1\over 4n}\big)\Re\alpha}
\end{array}
\eeq
and
\bel{T est.2}
\begin{array}{lr}\ds
\sum_{j>0\colon \left|1- {|\tau|\over |\xi|}\right|> 2^{-j}} \left|{^\sharp}\Hat{\Lambda}_j^{\alpha~\beta}(\xi,\tau)\right|
~\leq~\C_{\Re\alpha}~e^{\c\Im\alpha}e^{\c\Im\beta} 
 \sum_{j>0\colon \left|1- {|\tau|\over |\xi|}\right|> 2^{-j}}~ |\xi|^{-{\Re\alpha\over 2}}\left|1- {|\tau|\over |\xi|}\right|^{-{1\over 2}}2^{-\big({1\over 2n}\big)\Re\alpha j}
\qquad
 \hbox{\small{by (\ref{Result 1 > xi})}}
\\\\ \ds~~~~~~~~~~~~~~~~~~~~~~~~~~~~~~~~~~~~~~~~
~\leq~\C_{\Re\alpha}~e^{\c\Im\alpha}e^{\c\Im\beta} ~ |\xi|^{-{\Re\alpha\over 2}}\left|1- {|\tau|\over |\xi|}\right|^{-{1\over 2}+\big({1\over 4n}\big)\Re\alpha} ~\sum_{j>0} 2^{-\big({1\over 4n}\big)\Re\alpha j}
\\\\ \ds~~~~~~~~~~~~~~~~~~~~~~~~~~~~~~~~~~~~~~~~
~\leq~\C_{\Re\alpha}~e^{\c\Im\alpha}e^{\c\Im\beta} ~ |\xi|^{-{\Re\alpha\over 2}}\left|1- {|\tau|\over |\xi|}\right|^{-{1\over 2}+\big({1\over 4n}\big)\Re\alpha}.
\end{array}
\eeq
By putting together (\ref{T est.1})-(\ref{T est.2}), we obtain (\ref{Transform T Est}).

\subsection{Proof of Lemma One}

By using the identity of Bessel functions in (\ref{J identity}):
\[
\J_{a-1+\i b} (\rho)~=~2\left[{a+\i b\over \rho}\right] \J_{a+\i b}(\rho)-\J_{a+1+\i b}(\rho),\qquad a,b\in\R,\qquad \rho>0,
\]
we write  ${^\sharp}\Hat{\Lambda}^{\alpha~\beta}_j(\xi,\tau)$  as
\bel{Lambda sharp rewrite}
\begin{array}{lr}\ds
{^\sharp}\Hat{\Lambda}^{\alpha~\beta}_j(\xi,\tau)
~=~\int_\R e^{-2\pi\i \tau r} \left({1\over|r\xi|}\right)^{\delta(\alpha)-1} \J_{\delta(\alpha)-1}\Big(2\pi|r\xi|\Big) \phi_j(r\xi) |r|^{{\beta\over 2}-1}dr
\\\\ \ds~~~~~~~~~~~~~~~~~
~=~{1\over \pi}\delta(\alpha)\int_\R e^{-2\pi\i \tau r} \left({1\over|r\xi|}\right)^{\delta(\alpha)} \J_{\delta(\alpha)}\Big(2\pi|r\xi|\Big) \phi_j(r\xi) |r|^{{\beta\over 2}-1}dr
\\\\ \ds~~~~~~~~~~~~~~~~~
~-\int_\R e^{-2\pi\i \tau r} \left({1\over|r\xi|}\right)^{\delta(\alpha)-1} \J_{\delta(\alpha)+1}\Big(2\pi|r\xi|\Big) \phi_j(r\xi) |r|^{{\beta\over 2}-1}dr.
\end{array}
\eeq
Let $R=r|\xi|$. We have
\bel{Lambda sharp minor}
\begin{array}{lr}\ds
\int_\R e^{-2\pi\i \tau r} \left({1\over|r\xi|}\right)^{\delta(\alpha)} \J_{\delta(\alpha)}\Big(2\pi|r\xi|\Big) \phi_j(r\xi) |r|^{{\beta\over 2}-1}dr
\\\\ \ds
=~|\xi|^{-{\beta\over 2}}\int_\R e^{-2\pi\i \tau R/|\xi|} \left({1\over|R|}\right)^{\delta(\alpha)} \J_{\delta(\alpha)}\Big(2\pi|R|\Big) \phi_j(R) |R|^{{\beta\over 2}-1}dR.
\end{array}
\eeq
Recall $\delta(\alpha)=\Big({n+1\over 2n}\Big)\alpha$. By using the norm estimate of Bessel functions in (\ref{J norm}):
\[
\left| {1\over \rho^{a+\i b}}~\J_{a+\i b}(\rho)\right|~\leq~\C_a~\left({1\over 1+\rho}\right)^{{1\over 2}+a}e^{\c|b|},\qquad a, b\in\R,\qquad \rho>0,
\] 
we find
\bel{Lambda minor norm est}
\begin{array}{lr}\ds
\left|\left({1\over|R|}\right)^{\delta(\alpha)} \J_{\delta(\alpha)}\Big(2\pi|R|\Big)\right|~\leq~\C_{\Re\alpha}~e^{\c\Im\alpha} \left({1\over 1+|R|}\right)^{{1\over 2}+\big({n+1\over 2n}\big)\Re\alpha}.
\end{array}
\eeq
Note that  $\phi_j(\xi)$ is supported in the dyadic shell: $2^{j-1}<|\xi|<2^{j+1}$. Let $0<\Re\alpha=\Re\beta<1$. From (\ref{Lambda sharp minor})-(\ref{Lambda minor norm est}), we have
\bel{Lambda sharp minor EST}
\begin{array}{lr}\ds
\left| {1\over \pi}\delta(\alpha)\int_\R e^{-2\pi\i \tau r} \left({1\over|r\xi|}\right)^{\delta(\alpha)} \J_{\delta(\alpha)}\Big(2\pi|r\xi|\Big) \phi_j(r\xi) |r|^{{\beta\over 2}-1}dr
\right|
\\\\ \ds
~\leq~|\xi|^{-{\Re\alpha\over2}} \int_{2^{j-1}<|R|<2^{j+1}} \left| \left({1\over|R|}\right)^{\delta(\alpha)} \J_{\delta(\alpha)}\Big(2\pi|R|\Big)\right|  |R|^{{\Re\alpha\over 2}-1}dR
\\\\ \ds
~\leq~\C_{\Re\alpha}~e^{\c\Im\alpha}  ~|\xi|^{-{\Re\alpha\over2}} \int_{2^{j-1}<|R|<2^{j+1}} \left({1\over 1+|R|}\right)^{{1\over 2}+\big({n+1\over 2n}\big)\Re\alpha} |R|^{{\Re\alpha\over 2}-1}dR
\\\\ \ds
~\leq~\C_{\Re\alpha}~e^{\c\Im\alpha} ~|\xi|^{-{\Re\alpha\over2}}~2^{-j \big[{1\over 2}+\big({1\over 2n}\big)\Re\alpha\big]},\qquad j>0.
\end{array}
\eeq
Next, by using the asymptotic expansion of Bessel functions in (\ref{J asymptotic}), 
we write 
\bel{Lambda sharp major=S+E}
\begin{array}{lr}\ds
\int_\R e^{-2\pi\i \tau r} \left({1\over|r\xi|}\right)^{\delta(\alpha)-1} \J_{\delta(\alpha)+1}\Big(2\pi|r\xi|\Big) \phi_j(r\xi) |r|^{{\beta\over 2}-1}dr
\\\\ \ds
~=~\int_\R e^{-2\pi\i \tau r}
\Big[{^\sharp}\mathfrak{S}^\alpha\left(|r\xi|\right)+{^\sharp}\mathcal{E}^\alpha\left(|r\xi|\right)\Big]
\phi_j(r\xi) |r|^{{\beta\over 2}-1}dr
\end{array}
\eeq
where
\bel{S principal}
{^\sharp}\mathfrak{S}^\alpha\left(\rho\right)~=~{1\over \pi}  \left({1\over \rho}\right)^{\delta(\alpha)-{1\over 2}} \cos\left[2\pi \rho-{\pi\over 2}\delta(\alpha)+{\pi\over 4}\right],\qquad\rho>0
\eeq
and
\bel{E}
\left| {^\sharp}\mathcal{E}^\alpha\left(\rho\right)\right|~\leq~\C_{\Re\alpha}~e^{\c|\Im\alpha|} \left\{ \begin{array}{lr}\ds \rho^{-\Re\delta(\alpha)+{1\over 2}},\qquad 0<\rho\leq1,
\\ \ds
\rho^{-\Re\delta(\alpha)-{1\over 2}},\qquad~~~ \rho>1.
\end{array}
\right.
\eeq
By using  (\ref{E}), we have  
 \bel{E est}
\begin{array}{lr}\ds
\left| \int_\R e^{-2\pi\i \tau r} {^\sharp}\mathcal{E}^\alpha\left(r|\xi|\right)
|r|^{{\beta\over 2}-1}  dr  \right|
\\\\ \ds
~=~\left||\xi|^{-{\beta\over 2}} \int_\R e^{-2\pi\i \tau R/|\xi|} {^\sharp}\mathcal{E}^\alpha\left(|R|\right)
\phi_j(R)|R|^{{\beta\over 2}-1}  dR  \right|
\\\\ \ds
~\leq~\C_{\Re\alpha}~e^{\c|\Im\alpha|} ~
  |\xi|^{-{\Re\alpha\over2}}\int_{2^{j-1}<|R|<2^{j+1}}  \left({1\over |R|}\right)^{\Re\delta(\alpha)+{1\over 2}} |R|^{{\Re\alpha\over2}-1} dR 
\\\\ \ds
 ~\leq~\C_{\Re\alpha}~e^{\c|\Im\alpha|}~|\xi|^{-{\Re\alpha\over2}}~2^{\big[{\Re\alpha\over 2}-\big({n+1\over 2n}\big)\Re\alpha-{1\over2}\big]j}
 \\\\ \ds
 ~=~\C_{\Re\alpha}~e^{\c|\Im\alpha|}~|\xi|^{-{\Re\alpha\over2}}~2^{-\big[{1\over 2}+\big({1\over 2n}\big)\Re\alpha\big]j}. 
\end{array}
\eeq
Because of Euler's formulae, we replace the cosine function in (\ref{S principal}) with $e^{\pm2\pi\i \rho}$ multiplied by  $e^{\pm \pi\i \big[{\pi\over 4}-{1\over 2}\delta(\alpha)\big]}$ whose amplitude is controlled  by $\C_{\Re\alpha}e^{\c|\Im\alpha|}$. 

Assert
\bel{A_j sharp}
\begin{array}{lr}\ds
{^\sharp}\Hat{\A}_j^{\alpha~\beta}(\xi,\tau)~=~
\int_\R e^{2\pi\i \big[|\xi|-\tau\big] r} \left({1\over |r\xi|}\right)^{\delta(\alpha)-{1\over 2}}\phi_j(r\xi)|r|^{{\beta\over 2}-1}  dr
\\\\ \ds~~~~~~~~~~~~~~~~~
~=~|\xi|^{-{\beta\over 2}} \int_\R e^{2\pi\i \big[1-\tau/|\xi|\big] R} \left({1\over |R|}\right)^{\delta(\alpha)-{1\over 2}}\phi_j(R)|R|^{{\beta\over 2}-1}  dR.
\end{array}
\eeq
Suppose $\left|1-{|\tau|\over|\xi|}\right|\leq2^{-j}$. We have
\bel{A_j sharp Est.1}
\begin{array}{lr}\ds
\left| {^\sharp}\Hat{\A}_j^{\alpha~\beta}(\xi,\tau)\right|~\leq~\C_{\Re\alpha} ~|\xi|^{-{\Re\alpha\over 2}}\int_{2^{j-1}<|R|<2^{j+1}}   |R|^{{\Re\alpha\over2}-\big({n+1\over 2n}\big)\Re\alpha+{1\over 2}-1}  dR
\\\\ \ds~~~~~~~~~~~~~~~~~~~~
~\leq~\C_{\Re\alpha} ~|\xi|^{-{\Re\alpha\over 2}}~2^{j/2} 2^{-\big({1\over 2n}\big)\Re\alpha j},\qquad j>0.
\end{array}
\eeq
On the other hand, suppose $\left|1-{|\tau|\over|\xi|}\right|>2^{-j}$. By integration by parts $w.r.t~R$ in (\ref{A_j sharp}), we have
\bel{A_j sharp by parts}
\begin{array}{lr}\ds
{^\sharp}\Hat{\A}_j^{\alpha~\beta}(\xi,\tau)~=~{1\over 2\pi\i} |\xi|^{-{\beta\over 2}} \left[1-{\tau\over|\xi|}\right]^{-1}\int_\R e^{2\pi\i \big[1-\tau/|\xi|\big]R}\partial_R \Big[ \phi_j(R) |R|^{{\beta\over 2}-\delta(\alpha)+{1\over 2}-1}\Big] dR.
\end{array}
\eeq
By using the Chain-rule, a direct computation shows
\bel{integrand norm}
\begin{array}{lr}\ds
\left| \partial_R \Big[ \phi_j(R) |R|^{{\beta\over 2}-\delta(\alpha)+{1\over 2}-1}\Big]\right|~\leq
\\\\ \ds
\C~ 2^{-j} |R|^{{\Re\alpha\over 2}-\big({n+1\over 2n}\big)\Re\alpha+{1\over 2}-1}+
\C_{\Re\alpha} e^{\c|\Im\alpha|} e^{\c|\Im\beta|} ~|R|^{{\Re\alpha\over 2}-\big({n+1\over 2n}\big)\Re\alpha+{1\over 2}-2}.
\end{array}
\eeq
From (\ref{A_j sharp by parts})-(\ref{integrand norm}), we obtain
\bel{A_j sharp Est.2}
\begin{array}{lr}\ds
\left| {^\sharp}\Hat{\A}_j^{\alpha~\beta}(\xi)\right|~\leq~\C_{\Re\alpha}~e^{\c\Im\alpha}e^{\c\Im\beta} |\xi|^{-{\Re\alpha\over 2}}~\left|1-{\tau\over|\xi|}\right|^{-1} 
\\ \ds~~~~~~~~~~~~~~~~~~~~~~
\int_{2^{j-1}<|R|<2^{j+1}}  2^{-j} |R|^{{\Re\alpha\over 2}-\big({n+1\over 2n}\big)\Re\alpha+{1\over 2}-1}+ |R|^{{\Re\alpha\over 2}-\big({n+1\over 2n}\big)\Re\alpha+{1\over 2}-2} ~dR
\\\\ \ds~~~~~~~~~~~~~~~~~
~\leq~\C_{\Re\alpha}~e^{\c\Im\alpha}e^{\c\Im\beta}~|\xi|^{-{\Re\alpha\over 2}}~\left|1-{|\tau|\over|\xi|}\right|^{-1} 
~2^{\big[{\Re\alpha\over2}-\big({n+1\over 2n}\big)\Re\alpha+{1\over 2}-1\big]j} 
\\\\ \ds~~~~~~~~~~~~~~~~~
~=~\C_{\Re\alpha}~e^{\c\Im\alpha}e^{\c\Im\beta} ~|\xi|^{-{\Re\alpha\over 2}}~\left|1-{|\tau|\over|\xi|}\right|^{-1}  ~2^{-j/2} 2^{-\big({1\over 2n}\big)\Re\alpha j}
\\\\ \ds~~~~~~~~~~~~~~~~~
~\leq~\C_{\Re\alpha}~e^{\c\Im\alpha}e^{\c\Im\beta} ~|\xi|^{-{\Re\alpha\over 2}}~\left|1-{|\tau|\over|\xi|}\right|^{-{1\over2}}  2^{-\big({1\over 2n}\big)\Re\alpha j}.
\end{array}
\eeq
The same estimates in (\ref{A_j sharp Est.1}) and (\ref{A_j sharp Est.2}) hold by replacing $e^{2\pi\i r|\xi|}$ inside (\ref{A_j sharp}) with $e^{-2\pi\i r|\xi|}$. In particular, the term $1-\tau/|\xi|$ will be replaced by $1+\tau/|\xi|$.

Together with (\ref{E est}), we conclude (\ref{Result 1 > j})-(\ref{Result 1 > xi}). This finishes the proof of {\bf Lemma One}.

\subsection{Size of ${^\sharp}\Hat{\P}^{\alpha~\beta}_j(\xi,\tau)$ and ${^\sharp}\Hat{\Q}^{\alpha~\beta}_j(\xi,\tau)$}
Recall {\bf Remark 3.1}. We have $ {^\sharp}\Hat{\G}^{\alpha~\beta}_j(\xi,\tau)={^\sharp}\Hat{\P}^{\alpha~\beta}_j(\xi,\tau)+{^\sharp}\Hat{\Q}^{\alpha~\beta}_j(\xi,\tau)$ for which ${^\sharp}\Hat{\P}^{\alpha~\beta}_j, {^\sharp}\Hat{\Q}^{\alpha~\beta}_j$ are defined in (\ref{II sharp transform})-(\ref{III sharp transform}). 
In order to prove {\bf Proposition One}, it is suffice to show (\ref{Result One V}):
\[
\left| {^\sharp}\Hat{\Q}^{\alpha~\beta}_j(\xi,\tau)\right|~\leq~\C_{\Re \alpha} ~e^{\c|\Im \alpha|}e^{\c|\Im\beta|} ~|\xi|^{-{\Re\alpha\over2}}2^{-\ve j},\qquad \xi\neq0
\]
for every $j>0$ and some $\ve=\ve(\Re\alpha)>0$.

Given $\xi^\nu_j\in\mathds{S}^{n-1}$, 
$\varphi^\nu_j$ and $\Psi^{\alpha~\nu}_j$ are defined in (\ref{phi^v_j}) and (\ref{Psi nu j alpha again}) respectively.  We have
\[
\begin{array}{cc}\ds
{^\sharp}\Hat{\Q}^{\alpha~\beta}_j(\xi,\tau)~=~\int_\R e^{-2\pi\i \tau r} {^\sharp}\Hat{\V}^\alpha_j(|r|\xi) |r|^{{\beta\over 2}-1}dr,
\\\\ \ds
{^\sharp}\V^\alpha_j(x) ~=~
 \sum_{\nu\colon\xi^\nu_j\in\mathds{S}^{n-1}}  {^\sharp}\V^{\alpha~\nu}_j(x), 
\\\\ \ds
  {^\sharp}\V^{\alpha~\nu}_j(x)~=~{^\sharp}\Omega^{\alpha~\nu}_j(x) \Psi^{\alpha~\nu}_j(x),  \qquad
 {^\sharp}\Hat{\Omega}_j^{\alpha~\nu}(\xi)~=~\varphi^\nu_j(\xi) {^\sharp}\Hat{\Omega}^\alpha_j (\xi),\qquad j>0.
\end{array}
\]
The Fourier transform of ${^\sharp}\V^\alpha_j$ equals
\bel{V sharp Transform}
\begin{array}{lr}\ds
{^\sharp}\Hat{\V}^\alpha_j(\xi)~=~\int_{\R^n}e^{-2\pi\i x\cdot\xi}\left\{\sum_{\nu\colon\xi^\nu_j\in\mathds{S}^{n-1}} \Psi^{\alpha~\nu}_j(x){^\sharp}\Omega^{\alpha~\nu}_j(x)\right\}dx
\\\\ \ds~~~~~~~~~~~
~=~\sum_{\nu\colon\xi^\nu_j\in\mathds{S}^{n-1}}\int_{\R^n} \Hat{\Psi}^{\alpha~\nu}_j(\xi-\zeta)
{^\sharp}\Hat{\Omega}^{\alpha~\nu}_j(\zeta) d\zeta.
\end{array}
\eeq
Let $R=r|\xi|$. We write
\bel{Q sharp transform rewrite}
\begin{array}{cc}\ds
{^\sharp}\Hat{\Q}^{\alpha~\beta}_j(\xi,\tau)~=~\int_\R e^{-2\pi\i \tau r} {^\sharp}\Hat{\V}^\alpha_j(|r|\xi) |r|^{{\beta\over 2}-1}dr
\\\\ \ds~~~~~~~~~~~~~~~~~~~~~~~~~~~~~~~~~~~~~~~~~~
~=~|\xi|^{-{\beta\over2}}\int_\R e^{-2\pi\i \big({\tau/|\xi|}\big) R} {^\sharp}\Hat{\V}^\alpha_j\left(|R|{\xi\over|\xi|}\right) |R|^{{\beta\over 2}-1}dR,
\\\\ \ds~~~~~~~~~~~~~~~~~~~~~~~~~~~~~~~~~~~~~~~~~~~~~~~~~~~~~~~~~
~=~|\xi|^{-{\beta\over2}}\sum_{\nu\colon\xi^\nu_j\in\mathds{S}^{n-1}} \int_\R e^{-2\pi\i \big(\tau/|\xi|\big) R} {^\sharp}\Hat{\V}^{\alpha~\nu}_j\left(|R|{\xi\over|\xi|}\right) |R|^{{\beta\over 2}-1}dR,
\\\\ \ds~~~~~~
{^\sharp}\Hat{\V}^{\alpha~\nu}_j\left(|R|{\xi\over|\xi|}\right)~=~\int_{\R^n} \Hat{\Psi}^{\alpha~\nu}_j\left[|R|{\xi\over|\xi|}-\zeta\right]
{^\sharp}\Hat{\Omega}^{\alpha~\nu}_j(\zeta) d\zeta.
\end{array}
\eeq
Recall $\Psi^{\alpha~\nu}_j$  and ${^\sharp}\Hat{\Omega}_j^\alpha$ defined in (\ref{Psi nu j alpha again}) and (\ref{Omega sharp Transform})-(\ref{Omega_j sharp}). We have
\[\begin{array}{lr}\ds
\Psi^{\alpha~\nu}_j(x)~=~\varphi\left[ 2^{\big(1-\sigma\big)j}\left|1-\left(\L_\nu^T x\right)_1\right|\right] 
\prod_{i=2}^n \varphi\left[2^{\big[{1\over 2}-\sigma\big]j}\left|\left(\L_\nu^T x\right)_i\right|\right]
\\ \ds~~~~~~~~~~~~~
~+~\varphi\left[ 2^{\big(1-\sigma\big)j}\left|1+\left(\L_\nu^T x\right)_1\right|\right] 
\prod_{i=2}^n \varphi\left[2^{\big[{1\over 2}-\sigma\big]j}\left|\left(\L_\nu^T x\right)_i\right|\right]
\end{array}\]
supported in the rectangle
\[\begin{array}{cc}\ds
\hbox{R}^\nu_{j}~=~\hbox{R}^{\nu~+}_{j}\cup \hbox{R}^{\nu~-}_{j},
\\\\ \ds
\hbox{R}^{\nu~+}_{j}~=~\Bigg\{x\in\R^n\colon1-2^{\big(\sigma-1\big)j+1}<\left(\L_\nu^T x\right)_1<1+2^{\big(\sigma-1\big)j+1},~\left|\left(\L_\nu^T x\right)_i\right|<2^{\big[\sigma-{1\over 2}\big]j+1},~ i=2,\ldots,n
\Bigg\},
\\\\ \ds
\hbox{R}^{\nu~-}_{j}~=~ \Bigg\{x\in\R^n\colon-1-2^{\big(\sigma-1\big)j+1}<\left(\L_\nu^T x\right)_1<-1+2^{\big(\sigma-1\big)j+1},~\left|\left(\L_\nu^T x\right)_i\right|<2^{\big[\sigma-{1\over 2}\big]j+1},~ i=2,\ldots,n
\Bigg\}.
\end{array}\] 
The volume of $\hbox{R}^\nu_j$ is bounded by a constant multiple of $2^{\big(\sigma-1\big)j}2^{(n-1)\big[\sigma-{1\over 2}\big]j}$.

On the other hand,
\[
 \begin{array}{lr}\ds
 {^\sharp}\Hat{\Omega}^{\alpha~\nu}_j(\xi) ~=~\left({1\over|\xi|}\right)^{\delta(\alpha)-1} \J_{\delta(\alpha)-1}\Big(2\pi|\xi|\Big) \phi_j(\xi)\varphi^\nu_j(\xi) 
 \end{array}
\]
where $\varphi^\nu_j(\xi)$  defined in (\ref{phi^v_j}) is supported in the cone 
\[\Gamma^\nu_j~=~\Bigg\{\xi\in\R^n\colon \left|{\xi\over |\xi|}-\xi^{\nu}_{j}\right|<2^{-j/2+1}\Bigg\}.\]
Recall $\delta(\alpha)=\Big({n+1\over 2n}\Big)\alpha$. By using the norm estimate of Bessel functions in (\ref{J norm}), we find
\bel{J norm Est delta-1}
\begin{array}{lr}\ds
\left|\left({1\over|\xi|}\right)^{\delta(\alpha)-1} \J_{\delta(\alpha)-1}\Big(2\pi|\xi|\Big)\right|~\leq~\C_{\Re\alpha}~e^{\c\Im\alpha} \left({1\over 1+|\xi|}\right)^{\big({n+1\over 2n}\big)\Re\alpha-{1\over 2}}.
\end{array}
\eeq
Consider
\bel{second Int}
\begin{array}{lr}\ds
\int_{|R|\leq2^{j-\h}, |R|\ge2^{j+\h}} e^{-2\pi\i \big(\tau/|\xi|\big) R} {^\sharp}\Hat{\V}^{\alpha~\nu}_j\left(|R|{\xi\over|\xi|}\right) |R|^{{\beta\over 2}-1}dR
\\\\ \ds
=~\int_{|R|\leq2^{j-\h}, |R|\ge2^{j+\h}} e^{-2\pi\i \big(\tau/|\xi|\big) R} \left\{\int_{\R^n} \Hat{\Psi}^{\alpha~\nu}_j\left[|R|{\xi\over|\xi|}-\zeta\right]
{^\sharp}\Hat{\Omega}^{\alpha~\nu}_j(\zeta) d\zeta\right\} |R|^{{\beta\over 2}-1}dR
\\\\ \ds
=~\int_{|R|\leq2^{j-\h}, |R|\ge2^{j+\h}} e^{-2\pi\i \big(\tau/|\xi|\big) R}  \Bigg\{ \int_{\R^n} 
\\ \ds~~~~~
\Bigg\{\int_{\R^n}e^{-2\pi\i x\cdot\big[|R|{\xi\over |\xi|}-\zeta\big]}\varphi\left[ 2^{\big[1-\sigma\big] j}\left|1-\left(\L_\nu^T x\right)_1\right|\right] 
\prod_{i=2}^n \varphi\left[2^{\left[{1\over 2}-\sigma\right]j}\left|\left(\L_\nu^T x\right)_i\right|\right] dx
\\ \ds~~~~~
+\int_{\R^n}e^{-2\pi\i x\cdot\big[|R|{\xi\over |\xi|}-\zeta\big]}\varphi\left[ 2^{\big[1-\sigma\big] j}\left|1+\left(\L_\nu^T x\right)_1\right|\right] 
\prod_{i=2}^n \varphi\left[2^{\left[{1\over 2}-\sigma\right]j}\left|\left(\L_\nu^T x\right)_i\right|\right] dx\Bigg\}
\\ \ds~~~~~
\left({1\over|\zeta|}\right)^{\delta(\alpha)-1} \J_{\delta(\alpha)-1}\Big(2\pi|\zeta|\Big) \phi_j(\zeta)\varphi^\nu_j(\zeta) d\zeta\Bigg\} |R|^{{\beta\over 2}-1}dR,\qquad \h=2+\log_2 n.
\end{array}
\eeq
Suppose $|R|\leq2^{j-\h}$. There is an $i\in\{1,2,\ldots,n\}$ such that $|\zeta_i|>2^{j-1}/n$. We thus have $\Big||R| {\xi_i\over |\xi|}-\zeta_i\Big|\gtrsim 2^j$. By integration by parts $w.r.t$ $x_i$, we find
\bel{second Int Est.1}
\begin{array}{lr}\ds
\left| \Hat{\Psi}^{\alpha~\nu}_j\left[|R|{\xi\over|\xi|}-\zeta\right]\right|~\leq~
\left| \int_{\small{\hbox{R}}^{\nu~+}_j}e^{-2\pi\i x\cdot\big[|R|{\xi\over |\xi|}-\zeta\big]}\varphi\left[ 2^{\big[1-\sigma\big] j}\left|1-\left(\L_\nu^T x\right)_1\right|\right] 
\prod_{i=2}^n \varphi\left[2^{\left[{1\over 2}-\sigma\right]j}\left|\left(\L_\nu^T x\right)_i\right|\right] dx\right|
\\ \ds~~~~~~~~~~~~~~~~~~~~~~~~~~~~~~
~+~\left| \int_{\small{\hbox{R}}^{\nu~-}_j}e^{-2\pi\i x\cdot\big[|R|{\xi\over |\xi|}-\zeta\big]}\varphi\left[ 2^{\big[1-\sigma\big] j}\left|1+\left(\L_\nu^T x\right)_1\right|\right] 
\prod_{i=2}^n \varphi\left[2^{\left[{1\over 2}-\sigma\right]j}\left|\left(\L_\nu^T x\right)_i\right|\right] dx\right|
\\\\ \ds~~~~~~~~~~~~~~~~~~~~~~~~~~~~~
~\leq~\C_N~\left||R| {\xi_i\over |\xi|}-\zeta_i\right|^{-N} 2^{N(1-\sigma)j} ~2^{\big(\sigma-1\big) j} 2^{(n-1)\big[\sigma-{1\over 2}\big]j}
\\\\ \ds~~~~~~~~~~~~~~~~~~~~~~~~~~~~~
~\leq~\C_N~2^{-\sigma N j} 2^{\big(\sigma-1\big) j} 2^{(n-1)\big[\sigma-{1\over 2}\big]j},\qquad N\ge0.
\end{array}
\eeq
Suppose $|R|\ge2^{j+\h}$. There exists an $i\in\{1,2,\ldots,n\}$ such that $|\xi_i|>|\xi|/n$. We must have $\Big||R| {\xi_i\over |\xi|}-\zeta_i\Big|\gtrsim |R|$. By integration by parts $w.r.t$ $x_i$, we find
\bel{second Int Est.2}
\begin{array}{lr}\ds
\left| \Hat{\Psi}^{\alpha~\nu}_j\left[|R|{\xi\over|\xi|}-\zeta\right]\right|~\leq~\left| \int_{\small{\hbox{R}^{\nu~+}_j}}e^{-2\pi\i x\cdot\big[|R|{\xi\over |\xi|}-\zeta\big]}\varphi\left[ 2^{\big[1-\sigma\big] j}\left|1-\left(\L_\nu^T x\right)_1\right|\right] 
\prod_{i=2}^n \varphi\left[2^{\left[{1\over 2}-\sigma\right]j}\left|\left(\L_\nu^T x\right)_i\right|\right] dx\right|
\\ \ds~~~~~~~~~~~~~~~~~~~~~~~~~~~~~~
~+~\left| \int_{\small{\hbox{R}^{\nu~-}_j}}e^{-2\pi\i x\cdot\big[|R|{\xi\over |\xi|}-\zeta\big]}\varphi\left[ 2^{\big[1-\sigma\big] j}\left|1+\left(\L_\nu^T x\right)_1\right|\right] 
\prod_{i=2}^n \varphi\left[2^{\left[{1\over 2}-\sigma\right]j}\left|\left(\L_\nu^T x\right)_i\right|\right] dx\right|
\\\\ \ds~~~~~~~~~~~~~~~~~~~~~~~~~~~~~
~\leq~\C_N~\left||R| {\xi_i\over |\xi|}-\zeta_i\right|^{-N} 2^{N(1-\sigma)j} ~2^{\big(\sigma-1\big) j} 2^{(n-1)\big[\sigma-{1\over 2}\big]j}
\\\\ \ds~~~~~~~~~~~~~~~~~~~~~~~~~~~~~
~\leq~\C_N~|R|^{-\sigma N } 2^{\big(\sigma-1\big) j} 2^{(n-1)\big[\sigma-{1\over 2}\big]j},\qquad N\ge0.
\end{array}
\eeq
Let $0<\Re\alpha=\Re\beta<1$.  From (\ref{second Int}),  we have
\bel{second Int Est}
\begin{array}{lr}\ds
\left| \int_{|R|\leq2^{j-\h},|R|\ge 2^{j+\h}} e^{-2\pi\i \big(\tau/|\xi|\big) R} {^\sharp}\Hat{\V}^{\alpha~\nu}_j\left(|R|{\xi\over|\xi|}\right) |R|^{{\beta\over 2}-1}dR\right|
\\\\ \ds
~\leq~\int_{|R|\leq2^{j-\h}, |R|\ge2^{j+\h}}  \left\{\int_{\Gamma^\nu_j\cap\big\{2^{j-1}<|\zeta|<2^{j+1}\big\}} \left|\Hat{\Psi}^{\alpha~\nu}_j\left[|R|{\xi\over|\xi|}-\zeta\right] \right|
\left|\left({1\over|\zeta|}\right)^{\delta(\alpha)-1} \J_{\delta(\alpha)-1}\Big(2\pi|\zeta|\Big)\right| d\zeta\right\} |R|^{{\Re\alpha\over 2}-1}dR
\\\\ \ds
~\leq~\C_{\Re\alpha} e^{\c|\Im\alpha|}~ \left({1\over 1+2^j}\right)^{\big({n+1\over 2n}\big)\Re\alpha-{1\over 2}} 2^j 2^{\big({n-1\over 2}\big)j}
2^{\big(\sigma-1\big) j} 2^{(n-1)\big[\sigma-{1\over 2}\big]j} 
 2^{-\sigma Nj}\int_{ |R|\leq2^{j-\h}} |R|^{{\Re\alpha\over2}-1} dR
 \\ \ds
~+~
\C_{\Re\alpha} e^{\c|\Im\alpha|}~ \left({1\over 1+2^j}\right)^{\big({n+1\over 2n}\big)\Re\alpha-{1\over 2}} 2^j 2^{\big({n-1\over 2}\big)j}
2^{\big(\sigma-1\big) j} 2^{(n-1)\big[\sigma-{1\over 2}\big]j}  \int_{|R|\ge2^{j+\h}} |R|^{-\sigma N} |R|^{{\Re\alpha\over2}-1}dR
\\ \ds~~~~~~~~~~~~~~~~~~~~~~~~~~~~~~~~~~~~~~~~~~~~~~~~~~~~~~~~~~~~~~~~~~~~~~~~~~~~~~~~~~~~~~~~~~~~~~~~~~~~
\hbox{\small{by (\ref{J norm Est delta-1}) and (\ref{second Int Est.1})-(\ref{second Int Est.2})}}
\\\\ \ds
~\leq~\C_{\Re\alpha} e^{\c|\Im\alpha|}~ 2^{\big({n+1\over 2n}\big)\Re\alpha j} 2^{j/2} 2^{n\sigma j} 2^{\big[{\Re\alpha\over 2}-\sigma N\big]j}
~=~\C_{\Re\alpha} e^{\c|\Im\alpha|}~ 2^{-\big({1\over 2n}\big)\Re\alpha j} 2^{j/2} 2^{n\sigma j} 2^{-\sigma Nj},\qquad\hbox{\small{$N\ge0$}}.
\end{array}
\eeq
There are at most a constant multiple of $2^{\big({n-1\over 2}\big)j}$ many $\xi^\nu_j\in\mathds{S}^{n-1}$.
By using (\ref{second Int Est}) with $N$ chosen large enough, we conclude
\bel{V sharp comple Sum Est} \sum_{\xi^\nu_j\in\mathds{S}^{n-1}}\left| \int_{|R|\leq2^{j-\h},|R|\ge 2^{j+\h}} e^{-2\pi\i \big(\tau/|\xi|\big) R} {^\sharp}\Hat{\V}^{\alpha~\nu}_j\left(|R|{\xi\over|\xi|}\right) |R|^{{\beta\over 2}-1}dR\right|~\leq~\C_{\Re\alpha} e^{\c|\Im\alpha|} 2^{-\big({1\over 2n}\big)\Re\alpha}.
\eeq
Next, we consider
\bel{V v j sharp}
\begin{array}{lr}\ds
\int_{2^{j-\h}<|R|<2^{j+\h}} e^{-2\pi\i \big(\tau/|\xi|\big) R} {^\sharp}\Hat{\V}^{\alpha~\nu}_j\left(|R|{\xi\over|\xi|}\right) |R|^{{\beta\over 2}-1}dR
\\\\ \ds
=~
\int_{2^{j-\h}<|R|<2^{j+\h}} e^{-2\pi\i \big(\tau/|\xi|\big) R} \left\{ \int_{\R^n} \Hat{\Psi}^{\alpha~\nu}_j\left[|R|{\xi\over|\xi|}-\zeta\right]
{^\sharp}\Hat{\Omega}^{\alpha~\nu}_j(\zeta) d\zeta\right\} |R|^{{\beta\over 2}-1}dR
\\\\ \ds
=~\int_{2^{j-\h}<|R|<2^{j+\h}} e^{-2\pi\i \big(\tau/|\xi|\big) R}
 \Bigg\{ \int_{\R^n} \Bigg\{\int_{\R^n}e^{-2\pi\i x\cdot\big[|R|{\xi\over |\xi|}-\zeta\big]}\Psi^{\alpha~\nu}_j(x) dx\Bigg\} 
 \\ \ds~~~~~~~~~~~~~~~~~~~~~~~~
\left({1\over|\zeta|}\right)^{\delta(\alpha)-1} \J_{\delta(\alpha)-1}\Big(2\pi|\zeta|\Big) \phi_j(\zeta)\varphi^\nu_j(\zeta) d\zeta\Bigg\} |R|^{{\beta\over 2}-1}dR
 \\\\ \ds
 =~\int_{2^{j-\h}<|R|<2^{j+\h}} e^{-2\pi\i \big(\tau/|\xi|\big) R}
 \Bigg\{ \int_{\R^n} \Bigg\{\int_{\R^n}e^{-2\pi\i x\cdot\big[{\xi\over |\xi|}-\zeta\big]|R|}\Psi^{\alpha~\nu}_j(x)dx\Bigg\} 
  \\ \ds~~~~~~~~~~~~~~~~~~~~~~~~
 \left({1\over|R\zeta|}\right)^{\delta(\alpha)-1} \J_{\delta(\alpha)-1}\Big(2\pi|R\zeta|\Big) \phi_j(|R|\zeta)\varphi^\nu_j(|R|\zeta) d\zeta\Bigg\} |R|^{{\beta\over 2}+n-1}dR 
\qquad
 \hbox{\small{$\zeta\mt |R|\zeta$}}
 \\\\ \ds
 =~\iint_{\R^n\times\R^n}     {^\sharp}\G^{\alpha~\beta}_j(\xi,\tau,x,\zeta)  \Psi^{\alpha~\nu}_j(x)
 \varphi^\nu_j(\zeta)  dxd \zeta,\qquad\hbox{(\small{$\varphi^\nu_j$ is homogeneous of degree zero})}
\\\\ \ds
{^\sharp}\G^{\alpha~\beta}_j(\xi,\tau,x,\zeta)~=~
\\ \ds
\int_{2^{j-\h}<|R|<2^{j+\h}}  e^{-2\pi\i \big(\tau/|\xi|\big) R} e^{-2\pi\i x\cdot\big[{\xi\over |\xi|}-\zeta\big]|R|}
\left({1\over|R\zeta|}\right)^{\delta(\alpha)-1} \J_{\delta(\alpha)-1}\Big(2\pi|R\zeta|\Big) \phi_j(R\zeta) |R|^{{\beta\over 2}+n-1}dR.
 \end{array}
\eeq

\begin{lemma}  Let $0<\Re\alpha=\Re\beta<1$. Suppose $\left|{\xi\over |\xi|}-\zeta\right|\leq\c 2^{-j/2}$. For every $j>0$ and $x\in\supp\Psi^{\alpha~\nu}_j$, we have
\bel{Result G > j}
\left| {^\sharp}\G^{\alpha~\beta}_j(\xi,\tau,x,\zeta)\right|~\leq~\C_{\Re \alpha} ~e^{\c|\Im \alpha|}~2^{j/2} 2^{-\big({1\over4n}\big)\Re\alpha j} 2^{nj}\quad\hbox{if}\quad \left|{|\tau|\over|\xi|}-|\zeta|\right|\leq 2^{-j}
\eeq
and
\bel{Result G > xi}
\begin{array}{rl}\ds
\left| {^\sharp}\G^{\alpha~\beta}_j(\xi,\tau,x,\zeta)\right|~\leq~\C_{\Re \alpha} ~e^{\c|\Im \alpha|}e^{\c|\Im\beta|}~\left|{|\tau|\over|\xi|}-|\zeta|\right|^{-1+\big({1\over4n}\big)\Re\alpha} 2^{-\big({1\over4n}\big)\Re\alpha j}2^{nj}
\quad
\hbox{if}\quad \left|{|\tau|\over|\xi|}-|\zeta|\right|> 2^{-j}.
 \end{array}
\eeq
\end{lemma}

{\bf Proof}~~ Due to the symmetry reason, we only consider
\bel{G sharp +}
\begin{array}{lr}\ds
{^\sharp}\G^{\alpha~\beta~+}_j(\xi,\tau,x,\zeta)~=~
\\ \ds
\int_{2^{j-\h}<R<2^{j+\h}}  e^{-2\pi\i \big(\tau/|\xi|\big) R}  e^{-2\pi\i x\cdot\big[{\xi\over |\xi|}-\zeta\big]R}
\left({1\over|R\zeta|}\right)^{\delta(\alpha)-1} \J_{\delta(\alpha)-1}\Big(2\pi|R\zeta|\Big) \phi_j(R\zeta) R^{{\beta\over 2}+n-1}dR.
\end{array}
 \eeq
 
By using the identity of Bessel functions in (\ref{J identity}):
\[
\J_{a-1+\i b} (\rho)~=~2\left[{a+\i b\over \rho}\right] \J_{a+\i b}(\rho)-\J_{a+1+\i b}(\rho),\qquad a,b\in\R,\qquad \rho>0,
\]
we write (\ref{G sharp +}) as
\bel{G sharp rewrite}
\begin{array}{lr}\ds
{1\over \pi}\delta(\alpha)\int_{2^{j-\h}<R<2^{j+\h}}  e^{-2\pi\i \big(\tau/|\xi|\big) R} e^{-2\pi\i x\cdot\big[{\xi\over |\xi|}-\zeta\big]R}
\\ \ds~~~~~~~~~~~~~~~~~~~~~~~~~~~~~~
\left({1\over|R\zeta|}\right)^{\delta(\alpha)} \J_{\delta(\alpha)}\Big(2\pi|R\zeta|\Big) \phi_j(R\zeta) R^{{\beta\over 2}+n-1}dR
\\\\ \ds
-\int_{2^{j-\h}<R<2^{j+\h}}  e^{-2\pi\i \big(\tau/|\xi|\big) R}  e^{-2\pi\i x\cdot\big[{\xi\over |\xi|}-\zeta\big]R}
\\ \ds~~~~~~~~~~~~~~~~~~~~~~
\left({1\over|R\zeta|}\right)^{\delta(\alpha)-1} \J_{\delta(\alpha)+1}\Big(2\pi|R\zeta|\Big) \phi_j(R\zeta) R^{{\beta\over 2}+n-1}dR.
\end{array}
\eeq
Recall $\delta(\alpha)=\Big({n+1\over 2n}\Big)\alpha$. By  using the norm estimate of Bessel functions in (\ref{J norm}), we find
\bel{G minor norm est}
\begin{array}{lr}\ds
\left|\left({1\over|R\zeta|}\right)^{\delta(\alpha)} \J_{\delta(\alpha)}\Big(2\pi|R\zeta|\Big)\right|~\leq~\C_{\Re\alpha}~e^{\c\Im\alpha} \left({1\over 1+R}\right)^{{1\over 2}+\big({n+1\over 2n}\big)\Re\alpha}.
\end{array}
\eeq
Let $0<\Re\alpha=\Re\beta<1$. We have
\bel{Lambda sharp minor EST}
\begin{array}{lr}\ds
\Bigg| {1\over \pi}\delta(\alpha)\int_{2^{j-\h}<R<2^{j+\h}}  e^{-2\pi\i \big(\tau/|\xi|\big) R}  e^{-2\pi\i x\cdot\big[{\xi\over |\xi|}-\zeta\big]R}
\\ \ds~~~~~~~~~~~~~~~~~~~~~~~~~~~~~~~~
\left({1\over|R\zeta|}\right)^{\delta(\alpha)} \J_{\delta(\alpha)}\Big(2\pi|R\zeta|\Big) \phi_j(R\zeta) R^{{\beta\over 2}+n-1}dR
\Bigg|
\\\\ \ds
~\leq~\C_{\Re\alpha}~e^{\c\Im\alpha}  ~ \int_{2^{j-\h}<R<2^{j+\h}} \left({1\over 1+R}\right)^{{1\over 2}+\big({n+1\over 2n}\big)\Re\alpha} R^{{\Re\alpha\over 2}+n-1}dR\qquad\hbox{\small{by (\ref{G minor norm est})}}
\\\\ \ds
~\leq~\C_{\Re\alpha}~e^{\c\Im\alpha} ~2^{-j \big[{1\over 2}+\big({1\over 2n}\big)\Re\alpha\big]}2^{n j},\qquad j>0.
\end{array}
\eeq
By using the asymptotic expansion of Bessel functions in (\ref{J asymptotic}), 
we write (\ref{G sharp +}) as
\bel{G sharp major=S+E}
\begin{array}{lr}\ds
\int_{2^{j-\h}<R<2^{j+\h}}  e^{-2\pi\i \big(\tau/|\xi|\big) R}   e^{-2\pi\i x\cdot\big[{\xi\over |\xi|}-\zeta\big]R}
\Bigg[ {^\sharp}\mathfrak{S}^\alpha(|R\zeta|)+{^\sharp}\mathcal{E}^\alpha(|R\zeta|)\Bigg]\phi_j(R\zeta) R^{{\beta\over 2}+n-1}dR
\end{array}
\eeq
where
\[
{^\sharp}\mathfrak{S}^\alpha\left(\rho\right)~=~{1\over \pi}  \left({1\over \rho}\right)^{\delta(\alpha)-{1\over 2}} \cos\left[2\pi \rho-{\pi\over 2}\delta(\alpha)+{\pi\over 4}\right],\qquad\rho>0
\]
and
\[
\left| {^\sharp}\mathcal{E}^\alpha\left(\rho\right)\right|~\leq~\C_{\Re\alpha}~e^{\c|\Im\alpha|} \left\{ \begin{array}{lr}\ds \rho^{-\Re\delta(\alpha)+{1\over 2}},\qquad 0<\rho\leq1,
\\ \ds
\rho^{-\Re\delta(\alpha)-{1\over 2}},\qquad~~~ \rho>1.
\end{array}
\right.
\]
as shown in (\ref{S principal})-(\ref{E}).

We find
 \bel{E est G}
\begin{array}{lr}\ds
\Bigg| \int_{2^{j-\h}<R<2^{j+\h}} e^{-2\pi\i \big(\tau/|\xi|\big) R}   e^{-2\pi\i x\cdot\big[{\xi\over |\xi|}-\zeta\big]R}
{^\sharp}\mathcal{E}^\alpha\left(|R\zeta|\right)
\phi_j(R\zeta) R^{{\beta\over 2}+n-1}dR  \Bigg|
\\\\ \ds
~\leq~\C_{\Re\alpha}~e^{\c|\Im\alpha|} \int_{2^{j-1}<R<2^{j+1}}  \left({1\over R}\right)^{\Re\delta(\alpha)+{1\over 2}} R^{{\Re\alpha\over2}+n-1} dR 
\\\\ \ds
 ~\leq~\C_{\Re\alpha}~e^{\c|\Im\alpha|}~2^{\big[{\Re\alpha\over 2}-\big({n+1\over 2n}\big)\Re\alpha-{1\over2}\big]j} 2^{n j} 
 \\\\ \ds
 ~=~\C_{\Re\alpha}~e^{\c|\Im\alpha|}~2^{-\big[{1\over 2}+\big({1\over 2n}\big)\Re\alpha\big]j}2^{n j} . 
\end{array}
\eeq
Because of Euler's formulae, we replace ${^\sharp}\mathfrak{S}^\alpha\left(\rho\right)$ above with ${1\over \pi}  \left({1\over \rho}\right)^{\delta(\alpha)-{1\over 2}} e^{\pm2\pi\i \rho}$ multiplied by  $e^{\pm \pi\i \big[{\pi\over 4}-{1\over 2}\delta(\alpha)\big]}$ whose amplitude is controlled  by $\C_{\Re\alpha}e^{\c|\Im\alpha|}$. 

Assert
\bel{B_j sharp}
\begin{array}{lr}\ds
{^\sharp}\Hat{\B}_j^{\alpha~\beta}(\xi,\tau,x,\zeta)~=~
\int_{2^{j-\h}<R<2^{j+\h}} e^{-2\pi\i \big[{\tau\over|\xi|}-|\zeta| \big]R} 
e^{-2\pi\i x\cdot\big[{\xi\over |\xi|}-\zeta\big]R}
\left({1\over |R\zeta|}\right)^{\delta(\alpha)-{1\over 2}}
\phi_j(R\zeta) R^{{\beta\over 2}+n-1}dR.
\end{array}
\eeq
Suppose $\left|{|\tau|\over|\xi|}-|\zeta|\right|\leq2^{-j}$. We have
\bel{B_j sharp Est.1}
\begin{array}{lr}\ds
\left| {^\sharp}\Hat{\B}_j^{\alpha~\beta}(\xi,\tau,x,\zeta)\right|~\leq~\C_{\Re\alpha} \int_{2^{j-\h}<R<2^{j+\h}}   R^{-\big({n+1\over 2n}\big)\Re\alpha+{1\over 2}+{\Re\alpha\over2}+n-1}  dR
\\\\ \ds~~~~~~~~~~~~~~~~~~~~~~~~~~~
~\leq~\C_{\Re\alpha} ~2^{j/2} 2^{-\big({1\over 2n}\big)\Re\alpha j} 2^{n j},\qquad j>0.
\end{array}
\eeq
On the other hand, suppose $\left|{|\tau|\over|\xi|}-|\zeta|\right|>2^{-j}$. Let $\left|{\xi\over |\xi|}-\zeta\right|\leq \c2^{j/2}$. 
Recall $\Psi^{\alpha~\nu}_j$   defined in (\ref{Psi nu j alpha again}). For
$x\in\supp\Psi^{\alpha~\nu}_j$, we must have $|x|\approx1$.
By integration by parts $w.r.t~R$, we find 
\bel{B_j sharp by parts}
\begin{array}{lr}\ds
{^\sharp}\Hat{\B}_j^{\alpha~\beta}(\xi,\tau,x,\zeta)~=
\\ \ds
{1\over 2\pi\i} \left[{\tau\over|\xi|}-|\zeta| \right]^{-1}
\int_{2^{j-\h}<R<2^{j+\h}} e^{-2\pi\i \big[{\tau\over|\xi|}-|\zeta| \big]R} 
\partial_R\left\{e^{-2\pi\i x\cdot\big[{\xi\over |\xi|}-\zeta\big]R}\left({1\over |R\zeta|}\right)^{\delta(\alpha)-{1\over 2}}
\phi_j(R\zeta) R^{{\beta\over 2}+n-1}\right\}dR.
\end{array}
\eeq
By using the Chain-rule, a direct computation shows
\bel{integrand norm B}
\begin{array}{lr}\ds
\left| \partial_R \left\{ e^{-2\pi\i x\cdot\big[{\xi\over |\xi|}-\zeta\big]R} \left({1\over |R\zeta|}\right)^{\delta(\alpha)-{1\over 2}}
\phi_j(R\zeta) R^{{\beta\over 2}+n-1}\right\}\right|
\\\\ \ds
\leq~\C~2^{-j/2}R^{{\Re\alpha\over 2}-\big({n+1\over 2n}\big)\Re\alpha+{1\over 2}+n-1} +\C~ 2^{-j}~ R^{{\Re\alpha\over 2}-\big({n+1\over 2n}\big)\Re\alpha+{1\over 2}+n-1} 
\\ \ds~~~~
+\C_{\Re\alpha} e^{\c|\Im\alpha|} e^{\c|\Im\beta|}~ R^{{\Re\alpha\over 2}-\big({n+1\over 2n}\big)\Re\alpha+{1\over 2}+n-2}.
\end{array}
\eeq
From (\ref{B_j sharp by parts})-(\ref{integrand norm B}), we obtain
\bel{B_j sharp Est.2}
\begin{array}{lr}\ds
\left| {^\sharp}\Hat{\B}_j^{\alpha~\beta}(\xi,\tau,x,\zeta)\right|~\leq~
\C_{\Re\alpha}~e^{\c\Im\alpha}e^{\c\Im\beta} \left|{|\tau|\over|\xi|}-|\zeta|\right|^{-1} 
\\ \ds
\int_{2^{j-\h}<R<2^{j+\h}}2^{-j/2}R^{{\Re\alpha\over 2}-\big({n+1\over 2n}\big)\Re\alpha+{1\over 2}+n-1} +  2^{-j} R^{{\Re\alpha\over 2}-\big({n+1\over 2n}\big)\Re\alpha+{1\over 2}+n-1}+ R^{{\Re\alpha\over 2}-\big({n+1\over 2n}\big)\Re\alpha+{1\over 2}+n-2} ~dR
\\\\ \ds
~\leq~\C_{\Re\alpha}~e^{\c\Im\alpha}e^{\c\Im\beta}~\left|{|\tau|\over|\xi|}-|\zeta|\right|^{-1}
~2^{\big[{\Re\alpha\over2}-\big({n+1\over 2n}\big)\Re\alpha\big]j} 2^{n j} 
\\\\ \ds
~=~\C_{\Re\alpha}~e^{\c\Im\alpha}e^{\c\Im\beta} ~\left|{|\tau|\over|\xi|}-|\zeta|\right|^{-1} ~2^{-\big({1\over 2n}\big)\Re\alpha j} 2^{n j} 
\\\\ \ds
~\leq~\C_{\Re\alpha}~e^{\c\Im\alpha}e^{\c\Im\beta} ~\left|{|\tau|\over|\xi|}-|\zeta|\right|^{-1+\big({1\over 4n}\big)\Re\alpha} 2^{-\big({1\over 4n}\big)\Re\alpha j} 2^{n j}.
\end{array}
\eeq
The same estimates in (\ref{B_j sharp Est.1}) and (\ref{B_j sharp Est.2}) hold by replacing $e^{2\pi\i R|\zeta|}$   inside (\ref{B_j sharp}) with $e^{-2\pi\i R|\zeta|}$. In particular, the term ${\tau\over|\xi|}-|\zeta|$ will be replaced by ${\tau\over|\xi|}+|\zeta|$.

Together with (\ref{E est G}), we conclude (\ref{Result G > j})-(\ref{Result G > xi}) for the integral in (\ref{G sharp +}). By carrying out an analogue of (\ref{G sharp +})-(\ref{B_j sharp Est.2}), we obtain the same result for  
\[\begin{array}{lr}\ds
 {^\sharp}\G^{\alpha~\beta~-}_j(\xi,\tau,x,\zeta)~=
 \\ \ds
 \int_{-2^{j+\h}<R<-2^{j-\h}}  e^{-2\pi\i \big(\tau/|\xi|\big) R} e^{2\pi\i \big[{\tau\over|\xi|}-|\zeta| \big]R} 
\left({1\over|R\zeta|}\right)^{\delta(\alpha)-1} \J_{\delta(\alpha)-1}\Big(2\pi|R\zeta|\Big) \phi_j(R\zeta) |R|^{{\beta\over 2}-1}dR.
\end{array}\]
\endproof

Return to (\ref{V v j sharp}). We have
\[
\begin{array}{lr}\ds
\int_{2^{j-\h}<|R|<2^{j+\h}} e^{-2\pi\i \big(\tau/|\xi|\big) R} {^\sharp}\Hat{\V}^{\alpha~\nu}_j\left(|R|{\xi\over|\xi|}\right) |R|^{{\beta\over 2}-1}dR
\\\\ \ds
~=~
\iint_{\R^n\times\R^n}     {^\sharp}\G^{\alpha~\beta}_j(\xi,\tau,x,\zeta)  \Psi^{\alpha~\nu}_j(x)
 \varphi^\nu_j(\zeta)  dxd \zeta.
\end{array}
\]
Recall $\Psi^{\alpha~\nu}_j$  supported in the rectangles  
$\hbox{R}^\nu_j=\hbox{R}^{\nu~+}_j\cup\hbox{R}^{\nu~-}_j$ whose volume  is bounded by a constant multiple of $2^{\big(\sigma-1\big)j}2^{(n-1)\big[\sigma-{1\over 2}\big]j}$.
Moreover,  $\varphi^\nu_j(\xi)$ defined in (\ref{phi^v_j})  has a  support inside the cone $\Gamma^\nu_j=\Bigg\{\xi\in\R^n\colon \left|{\xi\over |\xi|}-\xi^{\nu}_{j}\right|<2^{-j/2+1}\Bigg\}$.
Because $\phi_j(\xi)$  supported  in the dyadic shell: $2^{j-1}<|\xi|<2^{j+1}$, we find $2^{-\h+1}<|\zeta|<2^{\h-1}$ whenever $ {^\sharp}\G^{\alpha~\beta}_j(\xi,\tau,x,\zeta)\neq0$. 

Note that $\h=2+\log_2 n$. Define
\bel{S_k}
\mathcal{S}_k~=~\Bigg\{ \zeta\in\R^n\colon 2^{-\h+1}<|\zeta|<2^{\h-1},~ 2^{-k}<\left| {|\tau|\over |\xi|}-|\zeta|\right|
\leq2^{-k+1}\Bigg\},\qquad k\ge -\log_2 n,~~      k\in\Z.
\eeq
The volume of $\mathcal{S}_k\cap\Gamma^\nu_j$ is bounded by a constant multiple of $2^{-k}2^{-\big({n-1\over 2}\big)j}$. 

Given $\xi\in\R^n$, we must have $\xi\in\Gamma^\mu_j$ for some $\xi^\mu_j\in\mathds{S}^{n-1}$. For every $\xi^\nu_j\in\mathds{S}^{n-1}$ such that $\left|\xi^\mu_j-\xi^\nu_j\right|\leq\c 2^{-j/2}$, we find 
\bel{condition xi zeta}
\left|{\xi\over |\xi|}-\zeta\right|~\leq~\c2^{-j/2},\qquad \zeta\in\Gamma^\nu_j.
\eeq
\begin{remark}
There are only a finitely many $\xi^\nu_j\in\mathds{S}^{n-1}$ of which $\left|\xi^\mu_j-\xi^\nu_j\right|\leq\c 2^{-j/2}$.
\end{remark}
Let $\zeta\in\mathcal{S}_k$ in (\ref{S_k}). Suppose $\left|{\xi\over |\xi|}-\zeta\right|\leq\c2^{-j/2}$ and $x\in\supp\Psi^{\alpha~\nu}_j$. 

For $k>j$, $\left| {|\tau|\over |\xi|}-|\zeta|\right|\leq2^{-k+1}\leq2^{-j}$. We have
\bel{G Est 1 k}
\begin{array}{lr}\ds
\left| {^\sharp}\G^{\alpha~\beta}_j(\xi,\tau,x,\zeta)\right|~\leq~\C_{\Re \alpha} ~e^{\c|\Im \alpha|}~2^{j/2} 2^{-\big({1\over 4n}\big)\Re\alpha j} 2^{nj}\qquad\hbox{\small{by (\ref{Result G > j})}}
\\\\ \ds~~~~~~~~~~~~~~~~~~~~~~~~~~~~~
\leq~\C_{\Re \alpha} ~e^{\c|\Im \alpha|}~2^{k/2} 2^{-\big({1\over 4n}\big)\Re\alpha j} 2^{nj}.
\end{array}
\eeq
For $k\leq j$, $\left| {|\tau|\over |\xi|}-|\zeta|\right|>2^{-k}\ge2^{-j}$. We find
\bel{G Est 2 k}
\begin{array}{lr}\ds
\left| {^\sharp}\G^{\alpha~\beta}_j(\xi,\tau,x,\zeta)\right|~\leq~\C_{\Re \alpha} ~e^{\c|\Im \alpha|}e^{\c|\Im \beta|}~\left| {|\tau|\over |\xi|}-|\zeta|\right|^{-1+\big({1\over 4n}\big)\Re\alpha} 2^{-\big({1\over 4n}\big)\Re\alpha j}2^{nj}\qquad\hbox{\small{by (\ref{Result G > xi})}}
\\\\ \ds~~~~~~~~~~~~~~~~~~~~~~~~~~~~
\leq~\C_{\Re \alpha} ~e^{\c|\Im \alpha|}e^{\c|\Im \beta|}~2^{\big[1-\big({1\over 4n}\big)\Re\alpha\big]k} 2^{-\big({1\over 4n}\big)\Re\alpha j} 2^{nj}.
\end{array}
\eeq
By taking into account (\ref{condition xi zeta}) and using (\ref{G Est 1 k})-(\ref{G Est 2 k}), we obtain
\bel{G Est Sum}
\begin{array}{lr}\ds
 \sum_{\xi^\nu_j\in\mathds{S}^{n-1}\colon~|\xi^\mu_j-\xi^\nu_j|\leq\c 2^{-j/2}}\left|\iint_{\R^n\times\R^n}     {^\sharp}\G^{\alpha~\beta}_j(\xi,\tau,x,\zeta) \Psi^{\alpha~\nu}_j(x)
 \varphi^\nu_j(\zeta)  dx d \zeta\right|
\\\\ \ds
\leq~\C\sum_{k\ge -\log_2 n,~ k\in\Z} \iint_{\hbox{\small{R}}^\nu_j\times \big\{\Gamma^\nu_j\cap\mathcal{S}_k\big\}}    \left| {^\sharp}\G^{\alpha~\beta}_j(\xi,\tau,x,\zeta)\right|   dx d \zeta \qquad\hbox{\small{by {\bf Remark 5.1}}}
 \\\\ \ds
 \leq~\C_{\Re \alpha} ~e^{\c|\Im \alpha|}e^{\c|\Im \beta|} ~\sum_{k\ge -\log_2 n,~k\in\Z}2^{\big[1-\big({1\over 4n}\big)\Re\alpha\big]k} 2^{-\big({1\over 4n}\big)\Re\alpha j} 2^{nj} ~2^{\big(\sigma-1\big)j}2^{(n-1)\big[\sigma-{1\over 2}\big]j}~
  2^{-k} 2^{-\big({n-1\over 2}\big)j}
 \\\\ \ds
 \leq~\C_{\Re \alpha} e^{\c|\Im \alpha|}e^{\c|\Im \beta|} ~2^{-\big({1\over 4n}\big)\Re\alpha j}2^{n\sigma j} \sum_{k\ge -\log_2 n,~ k\in\Z}2^{-\big({1\over 4n}\big)\Re\alpha k} 
 \\\\ \ds
 \leq~\C_{\Re \alpha} e^{\c|\Im \alpha|}e^{\c|\Im \beta|} ~2^{-\big({1\over 8n}\big)\Re\alpha j}\qquad\hbox{\small{by choosing $n\sigma<({1\over 8n})\Re\alpha$ 
}}.
\end{array}
\eeq

From (\ref{V v j sharp}), we also have
\bel{H sharp}
\begin{array}{lr}\ds
\int_{2^{j-\h}<|R|<2^{j+\h}} e^{-2\pi\i \big(\tau/|\xi|\big) R} {^\sharp}\Hat{\V}^{\alpha~\nu}_j\left(|R|{\xi\over|\xi|}\right) |R|^{{\beta\over 2}-1}dR
\\\\ \ds
=~
\int_{2^{j-\h}<|R|<2^{j+\h}} e^{-2\pi\i \big(\tau/|\xi|\big) R} \left\{ \int_{\R^n} \Hat{\Psi}^{\alpha~\nu}_j\left[|R|{\xi\over|\xi|}-\zeta\right]
{^\sharp}\Hat{\Omega}^{\alpha~\nu}_j(\zeta) d\zeta\right\} |R|^{{\beta\over 2}-1}dR
\\\\ \ds
 =~\int_{2^{j-\h}<|R|<2^{j+\h}} e^{-2\pi\i \big(\tau/|\xi|\big) R}
 \Bigg\{ \int_{\R^n} \Bigg\{\int_{\R^n}e^{-2\pi\i x\cdot\big[{\xi\over |\xi|}-\zeta\big]|R|}\Psi^{\alpha~\nu}_j(x)dx\Bigg\} 
  \\ \ds~~~~~~~~~~~~~~~~~~~~~~~~
 \left({1\over|R\zeta|}\right)^{\delta(\alpha)-1} \J_{\delta(\alpha)-1}\Big(2\pi|R\zeta|\Big) \phi_j(|R|\zeta)\varphi^\nu_j(|R|\zeta) d\zeta\Bigg\} |R|^{{\beta\over 2}+n-1}dR 
\\\\ \ds
 =~\int_{2^{j-\h}<|R|<2^{j+\h}} e^{-2\pi\i \big(\tau/|\xi|\big) R}
 \Bigg\{ \int_{\R^n} {^\sharp}\H^{\alpha~\nu}_j(\xi,\zeta,R)
 \\ \ds~~~~~~~~~~~~~~~~~~~~~~~~
 \left({1\over|R\zeta|}\right)^{\delta(\alpha)-1} \J_{\delta(\alpha)-1}\Big(2\pi|R\zeta|\Big) \phi_j(|R|\zeta)\varphi^\nu_j(|R|\zeta) d\zeta\Bigg\} 
 |R|^{{\beta\over 2}+n-1}dR, 
\\\\ \ds
{^\sharp}\H^{\alpha~\nu}_j(\xi,\zeta,R)~=~\int_{\R^n}e^{-2\pi\i x\cdot\big[{\xi\over |\xi|}-\zeta\big]|R|}\Psi^{\alpha~\nu}_j(x)dx
\\\\ \ds~~~~~~~~~~~~~~~~~~~~~
~=~\int_{\R^n}e^{-2\pi\i x\cdot\big[{\xi\over |\xi|}-\zeta\big]|R|}\varphi\left[ 2^{\big[1-\sigma\big] j}\left|1-\left(\L_\nu^T x\right)_1\right|\right] 
\prod_{i=2}^n \varphi\left[2^{\left[{1\over 2}-\sigma\right]j}\left|\left(\L_\nu^T x\right)_i\right|\right]dx
\\ \ds~~~~~~~~~~~~~~~~~~~~~
~+~\int_{\R^n}e^{-2\pi\i x\cdot\big[{\xi\over |\xi|}-\zeta\big]|R|}\varphi\left[ 2^{\big[1-\sigma\big] j}\left|1+\left(\L_\nu^T x\right)_1\right|\right] 
\prod_{i=2}^n \varphi\left[2^{\left[{1\over 2}-\sigma\right]j}\left|\left(\L_\nu^T x\right)_i\right|\right]dx.
\end{array}
\eeq
Recall that $\L_\nu$ is an $n\times n$-orthogonal matrix  with $\det\L_\nu=1$.  Given $\xi^\nu_j\in\mathds{S}^{n-1}$, we require $\L^T_\nu\xi^\nu_j=(1,0)^T\in\R\times\R^{n-1}$. 
Write $\zeta=\L_\nu \eta$ and $x=\L_\nu u$ in (\ref{H sharp}). We find
\bel{H sharp rewrite}
\begin{array}{lr}\ds
{^\sharp}\H^{\alpha~\nu}_j(\xi,\zeta,R)~=~{^\sharp}\H^{\alpha~\nu}_j(\xi,\L_\nu \eta,R)
\\\\ \ds~~~~~~~~~~~~~~~~~~~~~~
~=~\int_{\R^n}e^{-2\pi\i u\cdot\big[{\L^T_\nu\xi\over |\xi|}-\eta\big]|R|}\varphi\left[ 2^{\big[1-\sigma\big] j}\left|1-u_1\right|\right] 
\prod_{i=2}^n \varphi\left[2^{\left[{1\over 2}-\sigma\right]j}\left|u_i\right|\right]du
\\ \ds~~~~~~~~~~~~~~~~~~~~~~
~+~\int_{\R^n}e^{-2\pi\i u\cdot\big[{\L^T_\nu\xi\over |\xi|}-\eta\big]|R|}\varphi\left[ 2^{\big[1-\sigma\big] j}\left|1+u_1\right|\right] 
\prod_{i=2}^n \varphi\left[2^{\left[{1\over 2}-\sigma\right]j}\left|u_i\right|\right]du.
\end{array}
\eeq
\begin{lemma}   Let $\xi\in\Gamma^\mu_j$ and $2^{-\h+1}<|\eta|\leq2^{\h-1}$. Suppose $\left|\xi^\mu_j-\xi^\nu_j\right|\ge \c 2^{-j/2}$ for some $\c>0$ large. Either we have $\left|\Big[\L_\nu^T\xi-\eta\Big]_1\right|\approx2$ or $\left|\Big[\L_\nu^T\xi-\eta\Big]_\imath\right|\gtrsim2^{-j/2}$ for some  $\imath\in\{2,\ldots,n\}$ whenever $\L_\nu \eta\in \Gamma^\nu_j$.
\end{lemma}

{\bf Proof}~~Note that 
\[\L^T_\nu\xi^\nu_j~=~(1,0)^T~\in~\R\times\R^{n-1}.\] 
Consider $\L_\nu \eta\in \Gamma^\nu_j$.
We have 
$\left| \L_\nu^T \xi^\nu_j-(1,0)^T\right|=0$.
 By using the triangle inequality, we find
\bel{eq size est1}
\begin{array}{lr}\ds
{|\eta_i|\over|\eta|}~\leq~\left|{\eta\over|\eta|}-(1,0)^T\right|
\\\\ \ds~~~~~
~\leq~\left|{\eta\over|\eta|}-\L_\nu^T \xi^\nu_j\right|+\left|\L_\nu^T \xi^\nu_j-(1,0)^T\right|
\\\\ \ds~~~~~
~\leq ~2^{-j/2+1},\qquad i=2,\ldots,n.
\end{array}
\eeq
On the other hand, recall $\left|\xi^\mu_j-\xi^\nu_j\right|\ge\c2^{-j/2}$ for some $\c>0$ large. 
Let $\xi\in\Gamma^{\mu}_j$. By using the triangle inequality, we find
\bel{eq size est2}
\begin{array}{lr}\ds
\left| {\L_\nu^T\xi\over |\xi|}- (1,0)^T\right|~=~\left|{\xi\over |\xi|}- \xi^\nu_j\right|
\\\\ \ds~~~~~~~~~~~~~~~~~~~~~~~~
~\ge~\left|\xi^{\mu}_j-\xi^\nu_j\right|-\left| {\xi\over |\xi|}- \xi^{\mu}_j\right|
\\\\ \ds~~~~~~~~~~~~~~~~~~~~~~~~
~\ge~\c 2^{-j/2}-2^{-j/2+1}   
~=~\big[\c-2\big]~ 2^{-j/2}.
\end{array}
\eeq
Suppose $\left| \left( {\L_\nu^T\xi\over |\xi|}\right)_1-1\right|\leq\left[\sum_{i=2}^n \left( {\L_\nu^T\xi\over |\xi|}\right)_i^2\right]^{1\over2}$. The estimate in (\ref{eq size est2})
 further implies
\bel{eq size est3}
\begin{array}{lr}\ds
\left| {\left(\L_\nu^T\xi\right)_\imath\over |\xi|}\right|~\ge~{1\over \sqrt{2}\sqrt{n-1}} \left| {\L_\nu^T\xi\over |\xi|}- (1,0)^T\right|
\\\\ \ds~~~~~~~~~~~~~~
~\ge~{1\over \sqrt{2}\sqrt{n-1}} \big[\c-2\big]~ 2^{-j/2}\qquad \hbox{\small{for some}}\qquad \imath\in\{2,\ldots,n\}.
\end{array}
\eeq
Suppose $\left| \left( {\L_\nu^T\xi\over |\xi|}\right)_1-1\right|>\left[\sum_{i=2}^n \left( {\L_\nu^T\xi\over |\xi|}\right)_i^2\right]^{1\over2}$. As a geometric fact, ${\L_\nu^T\xi\over |\xi|}$ belongs to the semi-sphere opposite to $(1,0)^T$. We either have $\left( {\L_\nu^T\xi\over |\xi|}\right)_1\approx-1$ or $\left| \left( {\L_\nu^T\xi\over |\xi|}\right)_\imath\right|>{1\over \sqrt{2}\sqrt{n-1}} \big[\c-2\big]~ 2^{-j/2}$ for some $\imath\in\{2,\ldots,n\}$ as (\ref{eq size est3}).
By taking into account $2^{-\h+1}<|\eta|\leq2^{\h-1}$ in (\ref{eq size est1}) where $\h=2+\log_2 n$, we find
$|\eta_\imath|\leq2^{\h-1} 2^{-j/2+1}$ .
Consequently, $\left|\left[{\L_\nu^T\xi\over|\xi|}-\eta\right]_\imath\right|\gtrsim2^{-j/2}$  provided that $\c>0$ is  large enough.
\endproof
A direct computation shows
\[
\begin{array}{lr}\ds
\left| \partial_{u_1}^N \varphi\left[ 2^{(1-\sigma) j}\left|1- u_1\right|\right]\right| 
~\leq~\C_N~2^{Nj}2^{-N\sigma j},
\\\\ \ds
\left| \partial_{u_i}^N \varphi\left[2^{\big[{1\over 2}-\sigma\big]j}\left|u_i\right|\right]\right|
~\leq~\C_N~2^{Nj/2}2^{-N\sigma j},\qquad i=2,\ldots,n
\end{array}
\]
for every $N\ge0$.

Let $\xi\in\Gamma^\mu_j$ and  $2^{-\h+1}<|\eta|\leq2^{\h-1}$, $2^{j-1}<|R|<2^{j+1}$. Suppose $\left|\xi^\mu_j-\xi^\nu_j\right|\ge \c 2^{-j/2}$ for $\c>0$ large.
By applying {\bf Lemma 5.2} and carrying out an $N$-fold integration by parts either $w.r.t$ $u_1$ or  $u_\imath, \imath\in\{2,\ldots,n\}$ inside (\ref{H sharp}), we obtain
\bel{H sharp Est}
\begin{array}{lr}\ds
\left| {^\sharp}\H^{\alpha~\nu}_j(\xi,\L_\nu \eta,R)\right|~\leq~\C_N~2^{-N\sigma j} 2^{\big(\sigma-1\big)j}2^{(n-1)\big[\sigma-{1\over 2}\big]j}
\\\\ \ds~~~~~~~~~~~~~~~~~~~~~~~~~~~~~
~=~\C_N~	2^{-\big[N-n\big]\sigma j}  2^{-j}2^{-\big({n-1\over 2}\big)j}.
\end{array}
\eeq
Recall $\delta(\alpha)=\Big({n+1\over 2n}\Big)\alpha$. As in (\ref{J norm Est delta-1}), the norm estimate of Bessel functions in (\ref{J norm}) implies
\[
\begin{array}{lr}\ds
\left|\left({1\over|R\zeta|}\right)^{\delta(\alpha)-1} \J_{\delta(\alpha)-1}\Big(2\pi|R\zeta|\Big)\right|~\leq~\C_{\Re\alpha}~e^{\c\Im\alpha} \left({1\over 1+|R|}\right)^{\big({n+1\over 2n}\big)\Re\alpha-{1\over 2}},\quad\hbox{\small{$2^{-\h+1}<|\zeta|<2^{\h+1}$}}.
\end{array}
\]
\begin{remark}
There are at most a constant multiple of $2^{\big({n-1\over 2}\big)j}$ many $\xi^\nu_j\in\mathds{S}^{n-1}$.
\end{remark}
Let $0<\Re\alpha=\Re\beta<1$. We have
\bel{Sum Est comple}
\begin{array}{lr}\ds
\sum_{\xi^\nu_j\in\mathds{S}^{n-1}\colon~\left|\xi^\mu_j-\xi^\nu_j\right|>\c2^{-j/2}} \left| \int_{2^{j-\h}<|R|<2^{j+\h}} e^{-2\pi\i \big(\tau/|\xi|\big) R} {^\sharp}\Hat{\V}^{\alpha~\nu}_j\left(|R|{\xi\over|\xi|}\right) |R|^{{\beta\over 2}-1}dR\right|
\\\\ \ds
=~\sum_{\xi^\nu_j\in\mathds{S}^{n-1}\colon~\left|\xi^\mu_j-\xi^\nu_j\right|>\c2^{-j/2}} \Bigg|\int_{2^{j-\h}<|R|<2^{j+\h}} e^{-2\pi\i \big(\tau/|\xi|\big) R}
 \Bigg\{ \int_{\R^n} {^\sharp}\H^{\alpha~\nu}_j(\xi,\zeta,R)
 \\ \ds~~~~~~~~~~~~~~~~~~~~~~~~~~~~~~~~~~~
  \left({1\over|R\zeta|}\right)^{\delta(\alpha)-1} \J_{\delta(\alpha)-1}\Big(2\pi|R\zeta|\Big) \phi_j(|R|\zeta)\varphi^\nu_j(\zeta) d\zeta\Bigg\} 
 |R|^{{\beta\over 2}+n-1}dR   
 \\\\ \ds
\leq~\sum_{\xi^\nu_j\in\mathds{S}^{n-1}\colon~\left|\xi^\mu_j-\xi^\nu_j\right|>\c2^{-j/2}}  \int_{2^{j-\h}<|R|<2^{j+\h}} 
 \Bigg\{ \int_{\Gamma_j^\nu} \left|{^\sharp}\H^{\alpha~\nu}_j(\xi,\zeta,R)\right|  \left({1\over 1+|R|}\right)^{\big({n+1\over 2n}\big)\Re\alpha-{1\over 2}}d\zeta\Bigg\} 
 |R|^{{\Re\alpha\over 2}+n-1}dR  
 
\\\\ \ds
\leq~\C_{\Re\alpha~N}~e^{\c\Im\alpha} \sum_{\xi^\nu_j\in\mathds{S}^{n-1}\colon~\left|\xi^\mu_j-\xi^\nu_j\right|>\c2^{-j/2}} ~2^{-\big[N-n\big]\sigma j}  2^{-j}2^{-\big({n-1\over 2}\big)j}~2^{-\big({n+1\over 2n}\big)\Re\alpha j} 2^{j/2}2^{-\big({n-1\over 2}\big)j} ~2^{{\Re\alpha\over 2} j} 2^{nj}
\\\\ \ds
=~\C_{\Re\alpha~N}~e^{\c\Im\alpha} \sum_{\xi^\nu_j\in\mathds{S}^{n-1}\colon~\left|\xi^\mu_j-\xi^\nu_j\right|>\c2^{-j/2}} 

2^{-\big[N-n\big]\sigma j} 2^{-\big({1\over 2n}\big)\Re\alpha j} 2^{j/2}
\\\\ \ds
\leq~\C_{\Re\alpha~N}~e^{\c\Im\alpha}~2^{-\big[N-n\big]\sigma j} 2^{-\big({1\over 2n}\big)\Re\alpha j} 2^{j/2} 2^{\big({n-1\over 2}\big)j}\qquad\hbox{\small{by {\bf Remark 5.2}}}
\\\\ \ds
\leq~\C_{\Re\alpha~\sigma}~e^{\c\Im\alpha}~2^{-\big({1\over 2n}\big)\Re\alpha j}
\end{array}
\eeq
if $N$ is sufficiently large depending on $\sigma=\sigma(\Re\alpha)>0$, i.e: $\big[N-n\big]\sigma\ge{n+1\over 2}$.

\section{Convolution with a kernel having singularity on $\mathds{S}^{n-1}$}
\setcounter{equation}{0}
Define
\bel{I natural}
{^\natural}\I^\delta f (x)~=~\int_{\R^n} f(x-u)\left({1\over 1-|u|^2 }\right)^{1-\delta}du,\qquad \Re\delta>0.
\eeq

\begin{lemma} Let $0<\Re\delta<1$. We have
\bel{result two}
\begin{array}{cc}\ds
\left\| {^\natural}\I^\delta f\right\|_{\L^q(\R^n)}~\leq~\C_{\Re\delta~p}~\left\| f\right\|_{\L^p(\R^n)},
\\\\ \ds
 \Re\delta~=~{1\over p}-{1\over q},\qquad 1<p<q<\infty.
\end{array}
\eeq
\end{lemma}

{\bf Proof}~~For ${^\natural}\I^\delta$ defined in (\ref{I natural}), it is suffice to consider $f\ge0$ and $\delta\in\R\colon 0<\delta<1$. We prove {\bf Lemma 6.1} in the same spirt of Hedberg \cite{Hedberg}. 

Given  $0\leq\ell\in\Z$,   define
 \bel{Partial}
 \begin{array}{cc}\ds
 \Delta_\ell {^\natural}\I^\delta f (x)~=~\int_{\mathcal{S}_\ell} f(x-u)\left({1\over 1-|u|^2 }\right)^{1-\delta}du,
\\\\ \ds
\mathcal{S}_\ell~=~\left\{u\in\R^n~\colon 2^{-\ell-1}\leq 1-|u| < 2^{-\ell}\right\}.
\end{array}
\eeq 
Momentarily, we consider the partial sum operator 
\bel{Partial sum}
{^\natural}\I^\delta_\rho~=~\sum_{0\leq\ell\leq \rho} \Delta_\ell{^\natural}\I^\delta 
\eeq
for $\rho\in\Z$ sufficiently large

Let $\left\{u^\nu_\rho\right\}_\nu$ be
a collection of points  almost equally distributed on  $\mathds{S}^{n-1}$ as before having a grid length between $2^{-\rho-1}$ and $2^{-\rho}$:   

{\bf (1)} $|u^\mu_\rho-u^\nu_\rho|\ge 2^{-\rho-1}$ for every $u^\mu_\rho\neq u^\nu_\rho$. 

{\bf (2)} For any $w\in \mathds{S}^{n-1}$, there is a $u^\nu_\rho$ in the open set $\{u\in\mathds{S}^{n-1}\colon |u-w|<2^{-\rho+1}\}$.

Define
 \bel{cone}
\Gamma^\nu_\rho~=~\left\{u\in\R^n~\colon~\left| {u\over |u|}-u^\nu_\rho\right|~\leq~2^{-\rho}\right\}.
\eeq
Let $\Rec^\nu_{\rho~\ell}( u^\nu_\rho)$ be a rectangle  centered on $u^\nu_\rho\in\mathds{S}^{n-1}$. Assert $u^\nu_\rho$ as an unit vector. 
$\Rec^\nu_{\rho~\ell}( u^\nu_\rho)$ has one side parallel to $u^\nu_\rho$ with side length  $5\times2^{-\ell}$    and others  perpendicular to $u^\nu_\rho$ with side length equal to $5\times2^{-\rho}$. 

Denote
\bel{Rectangle}
\Rec^\nu_{\rho~\ell}=\left\{u\in\R^n\colon u+u^\nu_\rho\in \Rec^\nu_{\rho~\ell }( u^\nu_\rho)\right\}
\eeq
where 
\bel{Rec volumes}
\vol\Bigg\{ \Rec^\nu_{\rho~\ell }\Bigg\}~=~\vol\Bigg\{ \Rec^\p_{\rho~\ell }( u^\nu_\rho)\Bigg\}~=~5^n  2^{-\ell}2^{-\rho(n-1)}.
\eeq
Observe that
\bel{Inclusion R^v_j} 
\begin{array}{cc}\ds
\mathcal{S}_\ell\cap \Gamma^\nu_\rho~~\subset~~\Rec_{\rho~\ell}^\nu ( u^\nu_\rho),\qquad 0\leq\ell\leq \rho.
\end{array}
\eeq
 \begin{figure}[h]
\centering
\includegraphics[scale=0.40]{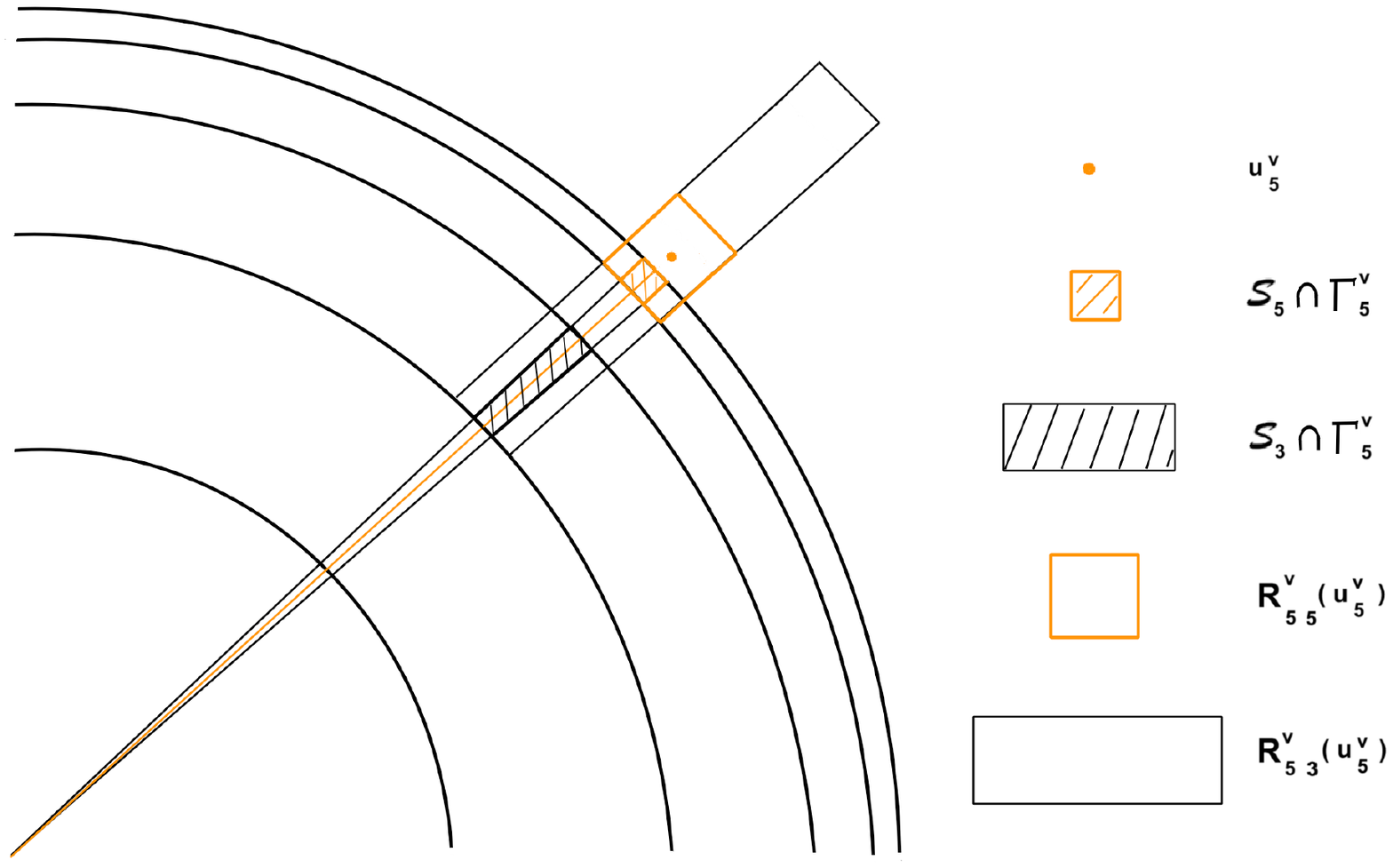}
\end{figure}

Define
\bel{M^v} 
\begin{array}{lr}\ds
 \M^\nu_\rho f(x)~=~\sup_{0\leq\ell\leq \rho} ~{1\over 2^{-\ell}2^{-\rho(n-1)}}\int_{\mathcal{S}_\ell\cap\Gamma^\nu_\rho} f(x-u)du
 \\\\ \ds~~~~~~~~~~~~
 ~\leq~\sup_{0\leq\ell\leq \rho} ~{5^n\over |\Rec^\nu_{\rho~\ell}|}\int_{\Rec_{\rho~\ell}^\nu} f\left(x-u^\nu_\rho-u\right)du.
 \end{array}
 \eeq
From the above, $\M^\nu_\rho$ is a composition of a $u^\nu_\rho$-translation and  a strong maximal operator  defined $w.r.t$ all rectangles  having a same degree centered on the origin. 

For every $\nu$, we have 
\bel{M^v L^p}
\left\| \M^\nu_\rho f\right\|_{\L^p(\R^n)}~\leq~\C_p~\left\| f\right\|_{\L^p(\R^n)},\qquad 1<p<\infty.
\eeq
Define
\bel{M*}
\M_\rho f(x) ~=~2^{-\rho(n-1)}\sum_{\nu\colon u^\nu_\rho\in\mathds{S}^{n-1}} \M^\nu_\rho f(x).
\eeq
There are at most a constant multiple of $2^{\rho\left(n-1\right)}$ elements  
 in $\left\{u^\nu_\rho\right\}_\nu$.  By using (\ref{M^v L^p}) and 
 applying Minkowski inequality, we obtain
\bel{M L^p}
\begin{array}{lr}\ds
\left\| \M_\rho f\right\|_{\L^p(\R^n)}~\leq~2^{-\rho(n-1)}\sum_{\nu\colon u^\nu_\rho\in\mathds{S}^{n-1}}~\left\| \M_\rho^\nu f\right\|_{\L^p(\R^n)}
\\\\ \ds~~~~~~~~~~~~~~~~~~~~~
~\leq~\C_p~\left\| f\right\|_{\L^p(\R^n)},\qquad 1<p<\infty.
\end{array}
\eeq
Given $f\in\L^p(\R^n)$, we claim 
  \bel{Regularity est}
 \begin{array}{cc}\ds
 {^\natural}\I^\delta_\rho f(x)
 ~\leq~\C_{p~q}~  \Big(\M_\rho f\Big)^{p\over q}(x) \left\| f\right\|_{\L^p(\R^n)}^{1-{p\over q}},
 \\\\ \ds
 \delta~=~{1\over p}-{1\over q},  \qquad1<p<q<\infty.
 \end{array}
\eeq
From (\ref{M L^p})-(\ref{Regularity est}), we  find
\bel{Norm Ineq natural eta}
\begin{array}{lr}\ds
\left\|{^\natural}\I^\delta_\rho f\right\|_{\L^q(\R^n)}~\leq~\C_{p~q}~\left\{\int_{\R^n}\Big(\M_\rho f\Big)^{p}(x)dx\right\}^{1\over q} \left\| f\right\|_{\L^p(\R^n)}^{1-{p\over q}}
\\\\ \ds~~~~~~~~~~~~~~~~~~~~
~\leq~\C_{p~q}~\left\|f\right\|_{\L^p(\R^n)}.
\end{array}
\eeq
In particular, (\ref{Norm Ineq natural eta}) is independent from $\rho$.
By allowing $\rho\mt\infty$  and using Fatou's lemma together with the monotone convergence theorem, we conclude (\ref{Result One}).

In order to prove (\ref{Regularity est}), we define $\zeta(x)\in\R$  implicitly by
\bel{zeta}
2^{\zeta(x)/p}~=~{ \M_\rho f(x)\over\left\| f\right\|_{\L^p(\R^n)}}.
\eeq
Let $\mathcal{S}_\ell$ and $\Gamma^\p_\rho$ defined in (\ref{Partial}) and (\ref{cone}). From (\ref{Inclusion R^v_j})-(\ref{M^v}), we have
\bel{Partial Est1 v}
\begin{array}{lr}\ds
\int_{\mathcal{S}_\ell \cap\Gamma^\nu_\rho}f(x-u)\left({1\over 1-|u|^2}\right)^{1-\delta} du
~=~\int_{\mathcal{S}_\ell \cap\Gamma^\nu_\rho}f(x-u)\left({1\over 1+|u|}\right)^{1-\delta} \left({1\over 1-|u|}\right)^{1-\delta}du
\\\\ \ds~~~~~~~~~~~~~~~~~~~~~~~~~~~~~~~~~~~~~~~~~~~~~~~~~~~~
~\leq~\C~   2^{\ell(1-\delta)}\int_{\mathcal{S}_\ell \cap\Gamma^\nu_\rho}f(x-u)du
\\\\ \ds~~~~~~~~~~~~~~~~~~~~~~~~~~~~~~~~~~~~~~~~~~~~~~~~~~~~
~=~\C~2^{\ell\left(1-\delta\right) }2^{-\ell}2^{-\rho(n-1)}\left\{{1\over 2^{-\ell}2^{-\rho(n-1)}}\int_{\mathcal{S}_\ell \cap\Gamma^\nu_\rho}f(x-u)du\right\}
\\\\ \ds~~~~~~~~~~~~~~~~~~~~~~~~~~~~~~~~~~~~~~~~~~~~~~~~~~~~
~=~\C~2^{-\ell\delta }2^{-\rho(n-1)}\left\{{1\over 2^{-\ell}2^{-\rho(n-1)}}\int_{\mathcal{S}_\ell \cap\Gamma^\nu_\rho}f(x-u)du\right\}
\\\\ \ds~~~~~~~~~~~~~~~~~~~~~~~~~~~~~~~~~~~~~~~~~~~~~~~~~~~~
~\leq~\C~2^{-\ell\delta }2^{-\rho(n-1)}\M_\rho^\nu f(x).
\end{array}
\eeq
Note that $\mathcal{S}_\ell\subset \Cup_{ u^\nu_\rho\in\mathds{S}^{n-1}} ~ \mathcal{S}_\ell \cap\Gamma^\p_\rho$.
By using (\ref{Partial Est1 v}) and summing over all $\nu$'s, we find
\bel{Partial Est1}
\begin{array}{lr}\ds
\int_{\mathcal{S}_\ell }f(x-u)\left({1\over 1-|u|^2}\right)^{1-\delta} du
~\leq~\C~2^{-\ell\delta }2^{-\rho(n-1)}\sum_{\nu\colon u^\nu_\rho\in\mathds{S}^{n-1}}\M_\rho^\nu f(x)
\\\\ \ds~~~~~~~~~~~~~~~~~~~~~~~~~~~~~~~~~~~~~~~~~~~~~~~
~=~\C~2^{-\ell\delta }\M_\rho f(x)
\end{array}
\eeq
where $\M_\rho$ is defined in (\ref{M*}).

On the other hand, by applying H\"{o}lder inequality, we have
\bel{Partial Est2 Holder}
\begin{array}{lr}\ds
\int_{\mathcal{S}_\ell }f(x-u)\left({1\over 1-|u|^2}\right)^{1-\delta} du
~\leq~\left\| f\right\|_{\L^p(\R^n)}\left\{\int_{\mathcal{S}_\ell }\left({1\over 1-|u|^2}\right)^{\left(1-\delta\right){p\over p-1}} du\right\}^{p-1\over p}
\end{array}
\eeq
where
\bel{Partial Est2}
\begin{array}{lr}\ds
\left\{\int_{\mathcal{S}_\ell }\left({1\over 1-|u|^2}\right)^{\left(1-\delta\right){p\over p-1}} du\right\}^{p-1\over p}
~=~\left\{\int_{\mathcal{S}_\ell }\left({1\over 1+|u|}\right)^{\left(1-\delta\right){p\over p-1}} \left({1\over 1-|u|}\right)^{\left(1-\delta\right){p\over p-1}} du\right\}^{p-1\over p}
\\\\ \ds~~~~~~~~~~~~~~~~~~~~~~~~~~~~~~~~~~~~~~~~~~~~~~~~~
~\leq~\C~\left[   {1\over 2^{-\ell\left(1-\delta\right){p\over p-1}}} 2^{-\ell}\right]^{p-1\over p}
\\\\ \ds~~~~~~~~~~~~~~~~~~~~~~~~~~~~~~~~~~~~~~~~~~~~~~~~~
~=~\C~2^{\ell\left(1-\delta\right)}2^{-\ell\big({p-1\over p}\big)}
\\\\ \ds~~~~~~~~~~~~~~~~~~~~~~~~~~~~~~~~~~~~~~~~~~~~~~~~~
~=~\C~2^{-\ell \big(\delta-{1\over p}\big)}.
\end{array}
\eeq
By inserting (\ref{zeta}) into (\ref{Partial Est1}), we find
\bel{Est1 partial}
\begin{array}{lr}\ds
\int_{\mathcal{S}_\ell}f(x-u)\left({1\over 1-|u|^2}\right)^{1-\delta} du
~\leq~\C~2^{-\ell\delta}\M_\rho f(x)
\\\\ \ds~~~~~~~~~~~~~~~~~~~~~~~~~~~~~~~~~~~~~~~~~~~~~~
~=~\C~2^{\big(\zeta(x)-\ell\big)\delta} 2^{-\zeta(x)\delta} \M_\rho f(x)
\\\\ \ds~~~~~~~~~~~~~~~~~~~~~~~~~~~~~~~~~~~~~~~~~~~~~~
~=~\C~2^{\big(\zeta(x)-\ell\big)\delta} 2^{-\zeta(x)\big({1\over p}-{1\over q}\big)} \M_\rho f(x)
\\\\ \ds~~~~~~~~~~~~~~~~~~~~~~~~~~~~~~~~~~~~~~~~~~~~~~
~=~\C~2^{\big(\zeta(x)-\ell\big)\delta} \left\{{\left\| f\right\|_{\L^p(\R^n)}\over \M_\rho f(x)}\right\}^{1-{p\over q}}  \M_\rho f(x)
\\\\ \ds~~~~~~~~~~~~~~~~~~~~~~~~~~~~~~~~~~~~~~~~~~~~~~
~=~\C~2^{\big(\zeta(x)-\ell\big)\delta}\Big(\M_\rho f\Big)^{p\over q}(x)\left\| f\right\|_{\L^p(\R^n)}^{1-{p\over q}}.
\end{array}
\eeq
By using (\ref{Est1 partial}) and summing over all $\zeta(x)\leq\ell\leq \rho$, we have
\bel{Est1 Partial}
\begin{array}{lr}\ds
\sum_{\zeta(x)\leq\ell\leq \rho}\int_{\mathcal{S}_\ell}f(x-u)\left({1\over 1-|u|^2}\right)^{1-\delta} du
~\leq~\C\sum_{\zeta(x)\leq\ell\leq \rho} 2^{\big(\zeta(x)-\ell\big)\delta}\Big(\M_\rho f\Big)^{p\over q}(x)\left\| f\right\|_{\L^p(\R^n)}^{1-{p\over q}}
\\\\ \ds~~~~~~~~~~~~~~~~~~~~~~~~~~~~~~~~~~~~~~~~~~~~~~~~~~~~~~~~~~~
~\leq~\C_\delta~\Big(\M_\rho f\Big)^{p\over q}(x)\left\| f\right\|_{\L^p(\R^n)}^{1-{p\over q}}.
\end{array}
\eeq
By inserting (\ref{zeta}) into (\ref{Partial Est2 Holder})-(\ref{Partial Est2}),  we find
\bel{Est2 partial}
\begin{array}{lr}\ds
\int_{\mathcal{S}_\ell }f(x-u)\left({1\over 1-|u|^2}\right)^{1-\delta} du
~\leq~\C~\left\| f\right\|_{\L^p(\R^n)} ~ 2^{-\ell\big(\delta-{1\over p}\big)}
\\\\ \ds~~~~~~~~~~~~~~~~~~~~~~~~~~~~~~~~~~~~~~~~~~~~~~~
~=~\C~2^{\big(\zeta(x)-\ell\big)\big(\delta-{1\over p}\big)} 2^{-\zeta(x)\big(\delta-{1\over p}\big)} \left\| f\right\|_{\L^p(\R^n)} 
\\\\ \ds~~~~~~~~~~~~~~~~~~~~~~~~~~~~~~~~~~~~~~~~~~~~~~~
~=~\C~2^{-\big(\zeta(x)-\ell\big){1\over q}} 2^{\zeta(x){1\over q}} \left\| f\right\|_{\L^p(\R^n)} 
\\\\ \ds~~~~~~~~~~~~~~~~~~~~~~~~~~~~~~~~~~~~~~~~~~~~~~~
~=~\C~2^{\big(\ell-\zeta(x)\big){1\over q}} 
\left\{{\M_\rho f(x)\over\left\| f\right\|_{\L^p(\R^n)}}\right\}^{p\over q}  \left\| f\right\|_{\L^p(\R^n)}
\\\\ \ds~~~~~~~~~~~~~~~~~~~~~~~~~~~~~~~~~~~~~~~~~~~~~~~
~=~\C~2^{\big(\ell-\zeta(x)\big){1\over q}} \Big(\M_\rho f\Big)^{p\over q}(x) \left\| f\right\|_{\L^p(\R^n)}^{1-{p\over q}}.
\end{array}
\eeq
By using (\ref{Est2 partial}) and summing over all $0\leq\ell<\zeta(x)$, we have
\bel{Est2 Partial}
\begin{array}{lr}\ds
\sum_{0\leq\ell<\zeta(x)} \int_{\mathcal{S}_\ell }f(x-u)\left({1\over 1-|u|^2}\right)^{1-\delta} du
~\leq~\C\sum_{0\leq\ell<\zeta(x)} 2^{\big(\ell-\zeta(x)\big){1\over q}} \Big(\M_\rho f\Big)^{p\over q}(x) \left\| f\right\|_{\L^p(\R^n)}^{1-{p\over q}}
\\\\ \ds~~~~~~~~~~~~~~~~~~~~~~~~~~~~~~~~~~~~~~~~~~~~~~~~~~~~~~~~~~~
~\leq~\C_{q}~ \Big(\M_\rho f\Big)^{p\over q}(x) \left\| f\right\|_{\L^p(\R^n)}^{1-{p\over q}}.
\end{array}
\eeq
By putting together (\ref{Est1 Partial}) and (\ref{Est2 Partial}), we obtain (\ref{Regularity est}) as desired. 
\endproof
Recall $\delta(\alpha)=\Big({n+1\over 2n}\Big)\alpha$ and ${^\flat}\Omega^{\alpha}$ defined in (\ref{Omega flat}):
\[
\begin{array}{cc}\ds
 {^\flat}\Omega^{\alpha}(x)~=~\pi^{\delta(\alpha)-1}\Gamma^{-1}\left(\delta(\alpha)\right) 
 \left({1 \over1-|x|^2}\right)^{1-\delta(\alpha)}_+,\qquad \hbox{\small{$\Re\delta(\alpha)>0$}}. 
 \end{array}
\]
Clearly, {\bf Lemma 6.1} implies in (\ref{Result 3.2}) in {\bf Lemma 3.2} which further implies (\ref{Result flat Omega_j}).

Recall ${^\flat}\U^\alpha_j$ and $ {^\flat}\V^\alpha_j$ defined in (\ref{II_j U flat})-(\ref{III_j V flat}). We find
\[
\begin{array}{lr}\ds
{^\flat}\U^{\alpha}_j(x)+{^\flat}\V^{\alpha}_j(x)
~=~\sum_{\xi^\nu_j\in\mathds{S}^{n-1}} {^\flat}\Omega^{\alpha~\nu}_j(x) \left[1-\Psi^{\alpha~\nu}_j(x)\right] +
{^\flat}\Omega^{\alpha~\nu}_j(x) \Psi^{\alpha~\nu}_j(x)
\\\\ \ds~~~~~~~~~~~~~~~~~~~~~~~~~~~
~=~\sum_{\xi^\nu_j\in\mathds{S}^{n-1}} {^\flat}\Omega^{\alpha~\nu}_j(x)
~=~{^\flat}\Omega^\alpha_j(x).
\end{array}
\]
In order to prove {\bf Proposition Two}, it is suffice to obtain (\ref{Result Two U}):
\[
\begin{array}{cc}\ds
\left\| f\ast {^\flat}\U^\alpha_j\right\|_{\L^q(\R^n)}~\leq~\C_{N~\Re\alpha~p} ~e^{\c|\Im\alpha|}~2^{-N j} \left\| f\right\|_{\L^p(\R^n)}, \qquad N\ge0
\\\\ \ds
\hbox{if}\qquad\Re\delta(\alpha)~=~{1\over p}-{1\over q},\qquad 1<p<q<\infty.
\end{array}
\]

\section{Proof of Proposition Two}
\setcounter{equation}{0}
Recall ${^\flat}\Omega^\alpha_j$ defined in (\ref{P v j flat}). We have
\[ \begin{array}{lr}\ds
~~~ {^\flat}\Hat{\Omega}_j^\alpha(\xi)~=~
\sum_{ \xi^\nu_j\in\mathds{S}^{n-1}}  {^\flat}\Hat{\Omega}_j^{\alpha~\nu}(\xi),
\\\\ \ds
{^\flat}\Hat{\Omega}_j^{\alpha~\nu}(\xi)~=~\varphi^\nu_j(\xi) {^\flat}\Hat{\Omega}^\alpha_j (\xi)
 \\\\ \ds~~~~~~~~~~~~~~~
 =~\varphi^\nu_j(\xi)\phi_j(\xi)\left({1\over|\xi|}\right)^{{n-1\over 2}+\delta(\alpha)-{1\over2}} \J_{{n-1\over 2}+\delta(\alpha)-{1\over2}}\Big(2\pi|\xi|\Big) 
 \end{array}
\]
where $\phi_j$ and $\varphi^\nu_j$  are defined in (\ref{varphi_j}) and (\ref{phi^v_j}).

By using
the asymptotic expansion of Bessel functions  in (\ref{J asymptotic})-(\ref{J O}), we write
\bel{Omega flat Sum}
\begin{array}{lr}\ds
{^\flat}\Omega_j^{\alpha~\nu}(x)
~=~\int_{\R^n}e^{2\pi\i x\cdot\xi}  \varphi^\nu_j(\xi)\phi(\xi)\left({1\over|\xi|}\right)^{{n-1\over 2}+\delta(\alpha)-{1\over2}} \J_{{n-1\over 2}+\delta(\alpha)-{1\over2}}\Big(2\pi|\xi|\Big) d\xi
\\\\ \ds~~~~~~~~~~~~~~
~=~{^\flat}\Q^{\alpha~\nu}_j(x)+\sum_{k=1}^L {^\flat}\Q^{\alpha~\nu}_{j~k}(x)+{^\flat}\hbox{\bf R}^{\alpha~\nu}_{j~L}(x).
\end{array}
\eeq
First,  we have
\bel{Q_j v flat}
\begin{array}{lr}\ds
{^\flat}\Q^{\alpha~\nu}_j(x)~=~{1\over \pi}\int_{\R^n}e^{2\pi\i x\cdot\xi} \varphi^\nu_j(\xi)\phi_j(\xi) \left({1\over |\xi|}\right)^{{n-1\over 2}+\delta(\alpha)} \cos\left[2\pi |\xi|-\left({n-1\over 2}+\delta(\alpha)-{1\over2}\right) {\pi\over 2}-{\pi\over 4}\right] d\xi.
\end{array}
\eeq
Next, there are  finitely many 
\bel{Q_j v k}
\begin{array}{lr}\ds
{^\flat}\Q^{\alpha~\nu}_{j~k}(x)~=~
\\\\ \ds
~~~~(2\pi)^{-2k}{\a_k\over \pi}\int_{\R^n}e^{2\pi\i x\cdot\xi} \varphi^\nu_j(\xi)\phi_j(\xi) \left({1\over |\xi|}\right)^{{n-1\over 2}+\delta(\alpha)+2k} \cos\left[2\pi |\xi|-\left({n-1\over 2}+\delta(\alpha)-{1\over2}\right) {\pi\over 2}-{\pi\over 4}\right] d\xi
\\\\ \ds
+~(2\pi)^{-2k+1}{\b_k\over \pi}\int_{\R^n}e^{2\pi\i x\cdot\xi} \varphi^\nu_j(\xi)\phi_j(\xi) \left({1\over |\xi|}\right)^{{n-1\over 2}+\delta(\alpha)+2k-1} \sin\left[2\pi |\xi|-\left({n-1\over 2}+\delta(\alpha)-{1\over2}\right) {\pi\over 2}-{\pi\over 4}\right] d\xi
\end{array}
\eeq
where $\a_k, \b_k$ are constants given in (\ref{J asymptotic}) for $k=1,2,\ldots,L$. 

Finally, the remainder term 
\bel{remainder term}
\begin{array}{cc}\ds
{^\flat}\hbox{\bf R}^{\alpha~\nu}_{j~L}(x)~=~\int_{\R^n}e^{2\pi\i x\cdot\xi} \varphi^\nu_j(\xi)\phi_j(\xi) \H^\alpha\Big(2\pi |\xi|\Big) d\xi,
\\\\ \ds
 \partial_\xi^\gamma \H^\alpha(\rho)~=~\O \left(\rho^{-2L-{1\over 2}}\right),\qquad  \rho\mt\infty
\end{array}
\eeq
where the implied constant  is bounded by $\C_{\gamma~\Re\alpha} e^{\c|\Im\alpha|}$.

Recall $\varphi^\nu_j$  defined in  (\ref{phi^v_j}). We find $\partial_\xi \varphi^\nu_j( \xi)=0$ whenever $\left|{\xi\over |\xi|}-\xi^\nu_j\right|\leq2^{-j/2}$. A direct computation shows
$\left|\partial_\xi^\gamma \left[\varphi_j^{\nu}(\xi)\right]\right|\leq\C_{\gamma} 2^{-|\gamma| j/2}$ 
for every multi-index $\gamma$ whenever $2^{j-1}<|\xi|< 2^{j+1}$.
The support of
$\varphi^\nu_j$ is contained in $\Gamma^\nu_j=\Bigg\{\xi\in\R^n\colon \left|{\xi\over |\xi|}-\xi^{\nu}_{j}\right|<2^{-j/2+1}\Bigg\}$.

 Let $0<\Re\alpha<1$ and $\delta(\alpha)=\big({n+1\over 2n}\big)\alpha$. An $(n+1)$-fold integration by parts $w.r.t~ \xi$ inside (\ref{remainder term}) gives
\bel{R EST}
\begin{array}{cc}\ds
\left| {^\flat}\hbox{\bf R}^{\alpha~\nu}_{j~L}(x)\right|~\leq~\C_{\Re\alpha}~ e^{\c|\Im\alpha|}~\left({1\over 1+|x|}\right)^{n+1}~2^{-\big[2L+{1\over 2}\big]j} 2^j2^{\big({n-1\over 2}\big)j}
\\\\ \ds~~~~~~~~~~~~~~~~~~~~~
~\leq~\C_{\Re\alpha}~ e^{\c|\Im\alpha|}~\left({1\over |x|}\right)^{n-\Re\delta(\alpha)}~2^{-j/2}2^{-\big[2L-{n+1\over 2}\big]j}.
\end{array}
\eeq
\begin{remark}
There are  a constant multiple of $2^{\big({n-1\over 2}\big)j}$ many $\xi^\nu_j\in\mathds{S}^{n-1}$. Moreover, $0\leq \Psi^{\alpha~\nu}_j(x)\leq1$. By using (\ref{R EST}), we obtain
\bel{Result R EST}
\begin{array}{lr}\ds
\sum_{\xi^\nu_j\in\mathds{S}^{n-1}} \left| {^\flat}\hbox{\bf R}^{\alpha~\nu}_{j~L}(x)\right| \left[1-\Psi^{\alpha~\nu}_j(x)\right] 
~\leq~\C_{\Re \alpha} ~e^{\c|\Im \alpha|}~ \left({1\over |x|}\right)^{n-\Re\delta(\alpha)}~2^{-j/2} 2^{-\big[2L-n\big]j}.
\end{array}
\eeq
\end{remark}
Recall $\Psi^{\alpha~\nu}_j$ defined in (\ref{Psi mu j alpha}).
Let $0<\Re\alpha<1$. We aim to show
\bel{Result U flat Q}
\begin{array}{lr}\ds
\sum_{\xi^\nu_j\in\mathds{S}^{n-1}} \left|{^\flat}\Q^{\alpha~\nu}_j(x)\right| \left[1-\Psi^{\alpha~\nu}_j(x)\right]
~\leq~\C_{N~\Re \alpha} ~e^{\c|\Im \alpha|}~ \left({1\over |x|}\right)^{n-\Re\delta(\alpha)}~2^{-N j},\qquad \hbox{\small{$N\ge0$}}
\end{array}
\eeq
for some $\ve=\ve(\Re\alpha)>0$ .
\begin{remark}  Later, it should be clear that (\ref{Result U flat Q}) holds if  ${^\flat}\Q^{\alpha~\nu}_j$  is replaced by ${^\flat}\Q^{\alpha~\nu}_{j~k}, k=1,2,\ldots,L$. Furthermore, an extra  term $2^{-j\big[2k-1\big]}$ appears on the right hand side of the inequality:
\[
\begin{array}{lr}\ds
\sum_{\xi^\nu_j\in\mathds{S}^{n-1}} \left|{^\flat}\Q^{\alpha~\nu}_{j~k}(x)\right| \left[1-\Psi^{\alpha~\nu}_{j~m}(x)\right] 
~\leq~\C_{N~\Re \alpha} ~e^{\c|\Im \alpha|}~\left({1\over |x|}\right)^{n-\Re\delta(\alpha)}~ 2^{-N j} 2^{-j\big[2k-1\big]},\qquad \hbox{\small{$N\ge0$}}.
\end{array}
\]
\end{remark}
From {\bf Remark 7.1} with $L$ chosen large enough, (\ref{Result U flat Q}) and {\bf Remark 7.2}, we obtain (\ref{Result Two U}) by applying the Hardy-Littlewood-Sobolev inequality  stated in (\ref{Hardy-Littlewood-Sobolev inequality}).

Recall ${^\flat}\Q^{\alpha~\nu}_j$ defined in (\ref{Q_j v flat}).
Because of Euler's formulae, we replace the cosine function  with $e^{\pm2\pi\i |\xi|}$  multiplied by  $e^{\pm\pi\i \big[{1\over 4}+{1\over 2}\big({n-1\over 2}+\delta(\alpha)-{1\over2}\big)\big]}$ whose amplitude is controlled by $\C_{\Re\alpha}~e^{\c|\Im\alpha|}$.    This assertion leads us to study
\bel{A flat j nu}
{^\flat}\A^{\alpha}_{j}(x)~=~\sum_{\xi^\nu_j\in\mathds{S}^{n-1}} {^\flat}\A^{\alpha~\nu}_{j}(x),
\qquad
{^\flat}\A^{\alpha~\nu}_{j}(x)~=~\int_{\R^n}e^{2\pi\i \big[x\cdot\xi- |\xi|\big]} \varphi^\nu_j(\xi) \phi_j(\xi)\left({1\over |\xi|}\right)^{{n-1\over 2}+\delta(\alpha)}  d\xi,
\eeq
\bel{B flat j nu}
{^\flat}\B^{\alpha}_{j}(x)~=~\sum_{\xi^\nu_j\in\mathds{S}^{n-1}} {^\flat}\B^{\alpha~\nu}_{j}(x),
\qquad
{^\flat}\B^{\alpha~\nu}_{j}(x)~=~\int_{\R^n}e^{2\pi\i \big[x\cdot\xi+ |\xi|\big]} \varphi^\nu_j(\xi) \phi_j(\xi)\left({1\over |\xi|}\right)^{{n-1\over 2}+\delta(\alpha)}  d\xi.
\eeq

\v
{\bf Lemma Two}~~{\it For every $j>0$, we have
\bel{Result Two A flat}
\sum_{\xi^\nu_j\in\mathds{S}^{n-1}} \left|{^\flat}\A^{\alpha~\nu}_j(x)\right|\left[1-\Psi^{\alpha~\nu}_{j}(x)\right] 
~\leq~\C_{N~\Re\alpha}~e^{\c|\Im\alpha|}~ \left({1\over |x|}\right)^{\Re\delta(\alpha)}~2^{-Nj} 2^{-\ve j},\qquad N\ge0
\eeq
and
\bel{Result Two B flat}
\sum_{\xi^\nu_j\in\mathds{S}^{n-1}} \left|{^\flat}\B^{\alpha~\nu}_{j}(x)\right|\left[1-\Psi^{\alpha~\nu}_{j}(x)\right] 
~\leq~\C_{N~\Re\alpha}~e^{\c|\Im\alpha|}~ \left({1\over |x|}\right)^{\Re\delta(\alpha)}~2^{-Nj} 2^{-\ve j},\qquad N\ge0
\eeq
where $\ve=\ve(\Re\alpha)>0$.}

\subsection{Proof of Lemma Two}
For symmetry reason, we only focus on the norm estimate of ${^\flat}\A^{\alpha~\nu}_{j}(x)$.

Recall  $\Hat{\phi}$ defined in (\ref{varphi_j}) which is radially symmetric, supported in the dyadic shell: $2^{j-1}<|\xi|<2^{j+1}$. 
Moreover, the support of $\varphi^\nu_j$ is contained in
\[\Gamma^\nu_j~=~\Bigg\{\xi\in\R^n\colon \left|{\xi\over |\xi|}-\xi^{\nu}_{j}\right|<2^{-j/2+1}\Bigg\}.\]
Given $\xi^\nu_j\in\mathds{S}^{n-1}$, $\L_\nu$ is an $n\times n$-orthogonal matrix with $\det\L_\nu=1$ and   
$\L_\nu^T \xi^\nu_j=(1,0)^T\in\R\times\R^{n-1}$.
Let $\xi=\L_\nu \eta$ and  $x=\L_\nu u$ in (\ref{A flat j nu}). We have
 \bel{A flat v rewrite}
 \begin{array}{cc}\ds
 {^\flat}\A^{\alpha~\nu}_{j}(\L_\nu u)~=~\int_{\R^n}e^{2\pi\i \big[u\cdot\eta- |\eta|\big]} \varphi^\nu_j(\L_\nu \eta)\phi_j\left(\eta\right) \left({1\over |\eta|}\right)^{{n-1\over 2}+\delta(\alpha)}  d\eta.
 \end{array}
 \eeq 
\begin{remark} Write $\eta=(\eta_1,\eta')\in\R\times\R^{n-1}$. 
 We have 
 $2^{j-1}<\eta_1<2^{j+1}$ and $ |\eta'|\lesssim 2^{j/2}$
  whenever $\varphi^\nu_j(\L_\nu \eta)\phi_j\left(\eta\right) \neq0$. 
\end{remark}
Consider
 \bel{Theta v j}
\begin{array}{lr}\ds
{^\flat}\A^{\alpha~\nu}_{j}(\L_\nu u)~=~\int_{\R^n}e^{2\pi\i  \big[u\cdot\eta- \eta_1\big]}\Theta^{\alpha~\nu}_{j}(\eta) d\eta,
\\\\ \ds
\Theta^{\alpha~\nu}_{j}(\eta)~=~e^{-2\pi\i \big[|\eta|-\eta_1\big]} \varphi^\nu_j(\L_\nu\eta)\phi_j\left(\eta\right) \left({1\over |\eta|}\right)^{{n-1\over 2}+\delta(\alpha)}.
 \end{array}
\eeq
{\bf Remark 7.3} implies
\bel{vol Theta}
\vol ~\Bigg\{\supp \Theta^{\alpha~\nu}_{j}\Bigg\}~\lesssim~2^j 2^{\left({n-1\over 2}\right)j}.
\eeq

We aim to show
\bel{Diff Ineq Theta N}
\left| \left(2^j\partial_{\eta_1}\right)^N \left(2^{j/2} \partial_{\eta'}\right)^\gamma  \Theta^{\alpha~\nu}_{j}(\eta)\right|~\leq~\C_{N~\gamma~\Re\alpha}~ e^{\c|\Im\alpha|}~2^{-\big[{n-1\over 2}+\Re\delta(\alpha)\big]j}
\eeq
for every $N\ge0$ and every multi-index $\gamma$.

{\bf 1.}~~Denote
\[\chi(\eta)~=~ |\eta|-\eta_1,\qquad \eta_1>0.\] 
We claim
\bel{d Est Psi}
\left|\partial_\eta^\gamma\chi(\eta)\right|~\leq~\C_\gamma~2^{-\gamma_1 j}2^{-\big[|\gamma|-\gamma_1\big]j/2}
\eeq
for every multi-index $\gamma$ with $|\gamma|\ge1$ whenever 
\[\L_\nu\eta~\in~\Gamma_j^\nu\cap \Bigg\{2^{j-1}<|\eta|< 2^{j+1}\Bigg\}.\]

Note that $\chi(\eta)$ is homogeneous of degree $1$ in $\eta$: $\chi(\rho\eta)=\rho\chi(\eta),~\rho>0$.  We have
\bel{Psi est1}
\left|\partial_\eta^\gamma \chi(\eta)\right|~\leq~\C_\gamma~|\eta|^{1-|\gamma|}
\eeq
for every multi-index $\gamma$.

Indeed, $\partial_\eta^\gamma \chi(\rho\eta)=\rho^{|\gamma|}\partial_\xi^\gamma\chi(\xi)$ for $\xi=\rho\eta$. On the other hand, $\partial_\eta^\gamma \chi(\rho\eta)=\partial_\eta^\gamma\left[\rho\chi(\eta)\right]=\rho\partial_\eta^\gamma\chi(\eta)$. Choose $\rho=|\eta|^{-1}$. We find $\partial_\eta^\gamma\chi(\eta)=|\eta|^{1-|\gamma|} \partial_\xi^\gamma\chi(\xi)$ where $|\xi|=1$.

Suppose $|\gamma|-\gamma_1\ge2$. From (\ref{Psi est1}), we find
\bel{Psi est2}
\begin{array}{lr}\ds
\left|\partial_\eta^\gamma \chi(\eta)\right|~\leq~\C_\gamma~2^{-\big(|\gamma|-1\big)j}
\\\\ \ds~~~~~~~~~~~~~~
~=~\C_\gamma~2^{-\gamma_1 j}2^{-\big[|\gamma|-\gamma_{1}-1\big]j} 
\\\\ \ds~~~~~~~~~~~~~~
~\leq~\C_\gamma~2^{-\gamma_1 j}2^{-\big[|\gamma|-\gamma_{1}\big]j/2}. 
\end{array}
\eeq
Suppose $|\gamma|-\gamma_{1}=0$ or $1$. Denote  $\eta^1=(\eta_1,0)^T\in\R\times\R^{n-1}$.
A direct computation shows
\[
 \begin{array}{cc}\ds
\left( \nabla_{\eta} \chi\right)(\eta^1)~=~0,
\\\\ \ds
\Big(\partial_{\eta_{1}}^N \nabla_{\eta'}\chi\Big)(\eta^1)
 ~=~ \Big(\nabla_{\eta'}\partial_{\eta_{1}}^N \chi\Big)(\eta^1)~=~0,\qquad \hbox{\small{$N\ge0$}}.
 \end{array}
\]
By writing out  the Taylor expansion of $\partial_{\eta_{1}}^N \chi(\eta)$ and $\nabla_{\eta'}\partial_{\eta_{1}}^N \chi(\eta)$  in the $\eta'$-subspace and using $\Big(\nabla_{\eta'}\partial_{\eta_{1}}^N \chi\Big)(\eta^1)=0$, we have
 \bel{p_eta_1 N norm}
 \begin{array}{lr}\ds
 \partial_{\eta_{1}}^N \chi(\eta)~=~ \O\left(|\eta'|^2|\eta|^{-N-1}\right),
\qquad
 \nabla_{\eta'} \partial_{\eta_{1}}^N \chi(x,\eta)~=~ \O\left(|\eta'||\eta|^{-N-1}\right),\qquad \hbox{\small{$N\ge1$}}.
 \end{array}
  \eeq
From (\ref{p_eta_1 N norm}) and {\bf Remark 7.3}, we conclude
\[
  \left|\partial_{\eta_{1}}^N \chi(\eta)\right|~\leq~\C_N  2^{- Nj},\qquad \left|\nabla_{\eta'}\partial_{\eta_{1}}^N \chi(\eta)\right|~\leq~\C_N  2^{-N j}2^{-j/2},\qquad \hbox{\small{$N\ge1$}}.
\]

 {\bf 2.}~~ Recall $\phi_j$ defined in (\ref{varphi_j}). We easily find
\bel{Diff Ineq phi}
\left| \partial_\eta^\gamma \left\{\phi_j\left(\eta\right) \left({1\over |\eta|}\right)^{{n-1\over 2}+\delta(\alpha)}\right\}\right|~\leq~\C_{\gamma~\Re\alpha} ~e^{\c|\Im\alpha|}~2^{-\big[{n-1\over 2}+\Re\delta(\alpha)\big]j} 2^{-|\gamma|j}
\eeq
for every multi-index $\gamma$.

{\bf 3.}~~Recall $\varphi^\nu_j$  defined in  (\ref{phi^v_j}). We have
\[
\varphi^\nu_j(\L_\nu \eta)~=~{\ds\varphi\Bigg[2^{j/2}\left|{\eta\over |\eta|}-\eta^{\nu}_{j}\right|\Bigg]\over \ds \sum_{\nu\colon\xi^\nu_j\in\mathds{S}^{n-1}} \varphi\Bigg[2^{j/2}\left|{\eta\over |\eta|}-\eta^{\nu}_{j}\right|\Bigg]}.
\]
Observe that $\partial_\eta \varphi^\nu_j(\L_\nu \eta)=0$ whenever $\left|{\eta\over |\eta|}-\eta^\nu_j\right|\leq2^{-j/2}$.

Suppose $\eta_1\approx2^j$ and $|\eta'|\lesssim2^{j/2}$. 
 From direct computation, we find
\bel{deri computation 1}
\begin{array}{lr}\ds
{\partial\over \partial \eta_1}{\eta_1\over |\eta|}~=~{1\over |\eta|}-{\eta_1^2\over |\eta|^3}~=~
{|\eta'|^2\over |\eta|^3},
\\\\ \ds
{\partial\over \partial \eta_1}{\eta_i\over |\eta|}~=~-{\eta_1\eta_i\over |\eta|^3},\qquad i~=~2,\ldots,n
\end{array}
\eeq
and
\bel{deri eta norm 1}
\begin{array}{lr}\ds
\left| {\partial\over \partial \eta_1}{\eta\over |\eta|}\right|~\approx~2^{j/2}|\eta|^{-2}~\approx~2^{-j/2}|\eta|^{-1}~\approx~2^{-j/2}2^{-j},
\\\\ \ds
\left| \left({\partial\over \partial \eta_1}\right)^N{\eta\over |\eta|}\right|~\leq~\C_N~2^{-j/2}2^{-j} 2^{-(N-1)j}~=~\C_N~2^{-j/2} 2^{-Nj},\qquad N\ge1.
\end{array}
\eeq
By carrying out the differentiation $w.r.t$ $(\eta_1,\eta')^T$ and taking into account for $\left|{\eta\over |\eta|}-\eta^\nu_j\right|>2^{-j/2}$, we have
\bel{nu Diff est 1}
\begin{array}{cc}\ds
\left|\partial_{\eta_1}^N \partial_{\eta'}^\gamma \varphi_j^\nu \left(\L_\nu\eta\right)\right|~\leq~\C_{N~\gamma} ~2^{-Nj} 2^{|\gamma| j/2}2^{-|\gamma|j}
\end{array}
\eeq
for every $N\ge0$ and every multi-index $\gamma$.

By putting together the above estimates in {\bf 1}-{\bf 3}, we obtain (\ref{Diff Ineq Theta N}).

Define a differential operator
\bel{D}
\D_j f~=~Id+2^{2j} \partial_{\eta_1}^2 f+2^j \Delta_{\eta'} f.
\eeq
Let ${^\flat}\A^{\alpha~\nu}_{j}(\L_\nu u)$ be given in (\ref{Theta v j}). An $N$-fold integration by parts $w.r.t~\D_j$ shows
\bel{INT by parts}
\begin{array}{lr}\ds
{^\flat}\A^{\alpha~\nu}_{j}(\L_\nu u)~=~
(2\pi\i)^{-2N} \left[1+2^{2j}\big[u_1- 1\big]^2+2^j \sum_{i=2}^n \big(u_i\big)^2\right]^{-N} 
\\\\ \ds~~~~~~~~~~~~~~~~~~~~~~~~~
\int_{\R^n}e^{2\pi\i  \big[u\cdot\eta- \eta_1\big]}\D^N_j \Theta^{\alpha~\nu}_{j}(\eta) d\eta.
\end{array}
\eeq
By using (\ref{vol Theta})-(\ref{Diff Ineq Theta N}), we have
\bel{Int norm}
\begin{array}{lr}\ds
\left|\int_{\R^n}e^{2\pi\i  \big[u\cdot\eta- \eta_1\big]}\D^N_j \Theta^{\alpha~\nu}_{j}(\eta) d\eta\right|
~\leq~\C_{N~\Re\alpha}~ e^{\c|\Im\alpha|}~2^{-\big[{n-1\over 2}+\Re\delta(\alpha)\big]j} 2^j 2^{\left({n-1\over 2}\right)j},\qquad N\ge0.
\end{array}
\eeq

Recall  $\Psi^{\alpha~\nu}_{j}$ defined in (\ref{Psi nu j alpha again}). We find
$0\leq\Psi^{\alpha~\nu}_{j~m}(x)\leq1$.
Define
\bel{Tilde Psi}
\Tilde{\Psi}^{\alpha~\nu}_{j}(x)~=~\left\{\begin{array}{lr}\ds 1,\qquad \Psi^{\alpha~\nu}_{j}(x)=1,
\\ \ds
0,\qquad \Psi^{\alpha~\nu}_{j}(x)\neq1.
\end{array}\right.
\eeq
Note that $\Tilde{\Psi}^{\alpha~\nu}_{j}\left(\L_\nu u\right)$ is supported in the rectangles
\bel{Rec}
\begin{array}{cc}\ds
\Tilde{\hbox{R}}_{j}~=~\Tilde{\hbox{R}}_{j}^{+}\cup \Tilde{\hbox{R}}_{j}^{-},
\\\\ \ds
\Tilde{\hbox{R}}_{j}^{+}~=~
\Bigg\{u\in\R^n\colon 1-2^{\big(\sigma-1\big) j}\leq u_1\leq 1+2^{\big(\sigma-1\big) j}, ~~|u_i|\leq 2^{\big[\sigma-{1\over 2}\big] j},~~i=2,\ldots,n\Bigg\},
\\\\ \ds
\Tilde{\hbox{R}}_{j}^{-}~=~
\Bigg\{u\in\R^n\colon -1-2^{\big(\sigma-1\big) j}\leq u_1\leq -1+2^{\big(\sigma-1\big) j}, ~~|u_i|\leq 2^{\big[\sigma-{1\over 2}\big] j},~~i=2,\ldots,n\Bigg\}.
\end{array}
\eeq
From (\ref{Tilde Psi})-(\ref{Rec}), we have
\bel{INT SUM EST}
\begin{array}{lr}\ds
\sum_{\xi^\nu_j\in\mathds{S}^{n-1}}  \left| {^\flat}\A^{\alpha~\nu}_{j}(x)\right| \left[ 1-\Psi^{\alpha~\nu}_{j}(x)\right]
\\\\ \ds
~\leq~\sum_{\xi^\nu_j\in\mathds{S}^{n-1}} \left| {^\flat}\A^{\alpha~\nu}_{r~j}(x)\right| \left[ 1-\Tilde{\Psi}^{\alpha~\nu}_{j}(x)\right]
\\\\ \ds
~=~\sum_{\xi^\nu_j\in\mathds{S}^{n-1}}\left| {^\flat}\A^{\alpha~\nu}_{r~j}(x)\right| \chi_{\R^n\setminus  \supp \Tilde{\Psi}^{\alpha~\nu}_{j}}(x)
\\\\ \ds
~=~\sum_{\xi^\nu_j\in\mathds{S}^{n-1}}\left| {^\flat}\A^{\alpha~\nu}_{r~j}(\L_\nu u)\right| \chi_{\R^n\setminus\Tilde{\hbox{\small{R}}}_j}(u)
\\\\ \ds
~\leq~\sum_{\xi^\nu_j\in\mathds{S}^{n-1}}\left| {^\flat}\A^{\alpha~\nu}_{r~j}(\L_\nu u)\right| \chi_{\R^n\setminus\Tilde{\hbox{\small{R}}}^+_j}(u)
\end{array}
\eeq
where $\chi$ denotes the indicator function.

Suppose $u_1\notin[1-2^{\big(\sigma-1\big) j}, 1+2^{\big(\sigma-1\big) j}]$.
We find
\bel{lambda r ratio}
\begin{array}{lr}\ds
u_1-1~>~\Big[1+2^{\big(\sigma-1\big) j}\Big]-1,
~=~2^{\big(\sigma-1\big) j},\qquad u_1>1+2^{\big(\sigma-1\big) j}
\\\\ \ds
u_1-1~<~\Big[1-2^{\big(\sigma-1\big) j}\Big]-1~=~-2^{\big(\sigma-1\big) j},\qquad u_1<1-2^{\big(\sigma-1\big) j}.
\end{array}
\eeq
As indicated in {\bf Remark 5.2}, there are  a constant multiple of $ 2^{\left({n-1\over 2}\right)j}$ many $\xi^\nu_j\in\mathds{S}^{n-1}$. 

 Let $0<\Re\alpha<1$ and $\delta(\alpha)=\Big({n+1\over 2n}\Big)\alpha$. Recall (\ref{INT by parts})-(\ref{Int norm}). We have
\bel{INT EST}
\begin{array}{lr}\ds
\sum_{\xi^\nu_j\in\mathds{S}^{n-1}}\left| {^\flat}\A^{\alpha~\nu}_{r~j}(\L_\nu u)\right| \chi_{\R^n\setminus\Tilde{\hbox{\small{R}}}^+_j}(u)
\\\\ \ds

~\leq~\C_{N~\Re\alpha}~ e^{\c|\Im\alpha|} \sum_{\xi^\nu_j\in\mathds{S}^{n-1}}2^{-\big[{n-1\over 2}+\Re\delta(\alpha)\big]j}~ 2^j 2^{\left({n-1\over 2}\right)j}
\left[1+2^{2j}\left[ u_1- 1\right]^2+2^j \sum_{i=2}^n \big(u_i\big)^2\right]^{-N}  \chi_{\R^n\setminus\Tilde{\hbox{\small{R}}}^+_j}(u)
\\\\ \ds

~\leq~\C_{N~\Re\alpha}~ e^{\c|\Im\alpha|} \sum_{\xi^\nu_j\in\mathds{S}^{n-1}}2^{-\big[{n-1\over 2}+\Re\delta(\alpha)\big]j}~ 2^j 2^{\left({n-1\over 2}\right)j}
\left[1+\left[ u_1- 1\right]^2+\sum_{i=2}^n \big(u_i\big)^2\right]^{-\lfloor{n+1\over 2}\rfloor}  
\\ \ds~~~~~~~~~~~~~~~~~~~~~~~~~~~~~~~~~~~~~~~~~
\left[1+2^{2j}\left[ u_1- 1\right]^2+2^j \sum_{i=2}^n \big(u_i\big)^2\right]^{-N+\lfloor{n+1\over 2}\rfloor}  \chi_{\R^n\setminus\Tilde{\hbox{\small{R}}}^+_j}(u)
\\\\ \ds

~\leq~\C_{N~\Re\alpha}~ e^{\c|\Im\alpha|} \sum_{\nu\colon\xi^\nu_j\in\mathds{S}^{n-1}}2^{-\big[{n-1\over 2}+\Re\delta(\alpha)\big]j}~ 2^j 2^{\left({n-1\over 2}\right)j}
\left[1+\left[ u_1- 1\right]^2+\sum_{i=2}^n \big(u_i\big)^2\right]^{-\lfloor{n+1\over 2}\rfloor}  2^{-2\big[N-\lfloor{n+1\over 2}\rfloor\big]\sigma j}
\\\\ \ds
~\leq~\C_{N~\Re\alpha}~ e^{\c|\Im\alpha|} ~2^{\big({n+1\over 2}\big)j} 2^{-\big[{n-1\over 2}+\Re\delta(\alpha)\big]j}~ 2^j 2^{\left({n-1\over 2}\right)j} 2^{-2\big[N-\lfloor{n+1\over 2}\rfloor\big]\sigma j} \left({1\over 1+|u|}\right)^{n+1}
\\\\ \ds
~\leq~\C_{N~\Re\alpha}~ e^{\c|\Im\alpha|} ~2^j2^{\big({n+1\over 2}\big)j} 2^{-\Re\delta(\alpha)j}~  2^{-2\big[N-\lfloor{n+1\over 2}\rfloor\big]\sigma j} \left({1\over |u|}\right)^{n-\Re\delta(\alpha)}.
\end{array}
\eeq
From (\ref{INT SUM EST}) and (\ref{INT EST}) with $N$ sufficiently large depending $\sigma=\sigma(\Re\alpha)>0$, we conclude
\bel{Final A flat}
\sum_{\xi^\nu_j\in\mathds{S}^{n-1}}  \left| {^\flat}\A^{\alpha~\nu}_{j}(x)\right| \left[ 1-\Psi^{\alpha~\nu}_{j}(x)\right] ~\leq~\C_{N~\Re\alpha}~ e^{\c|\Im\alpha|}~\left({1\over |x|}\right)^{n-\Re\delta(\alpha)}~2^{-N j},\qquad\hbox{\small{$N\ge 0$}}
\eeq
as desired.

\appendix 
  \section{Some estimates regarding Bessel functions}
 \setcounter{equation}{0}
 First, we recall certain well known estimates for Bessel functions. More discussion can be found in  the book of Watson \cite{Watson}. 
 
$\bullet$ For  $a>-{1\over 2}, b\in\R$ and $\rho>0$,  a Bessel function $\J_{a+\i b}$ has an integral formula 
\bel{Bessel}
\J_{a+\i b}(\rho)~=~{(\rho/2)^{a+\i b}\over \pi^{1\over 2}\Gamma\left(a+{1\over 2}+\i b\right)} \int_{-1}^1 e^{\i \rho s } (1-s^2)^{a-{1\over 2}+\i b} ds.
\eeq
$\bullet$ For every $a>-{1\over 2}, b\in\R$ and $\rho>0$, we have 
\bel{J asymptotic}
\begin{array}{lr}\ds
\J_{a+\i b}(\rho) ~\sim~\left({\pi \rho\over 2}\right)^{-{1\over 2}} \cos\left[\rho-(a+\i b) {\pi\over 2}-{\pi\over 4}\right]
\\\\ \ds~~~~~~~~~~~~
~+~\left({\pi \rho\over 2}\right)^{-{1\over 2}} \sum_{k=1}^\infty \cos\left[\rho-(a+\i b) {\pi\over 2}-{\pi\over 4}\right] \a_k \rho^{-2k}+\sin\left[\rho-(a+\i b) {\pi\over 2}-{\pi\over 4}\right]~ \b_k \rho^{-2k+1}
\end{array}
\eeq
where
\bel{a_kb_k}
\begin{array}{cc}\ds
\a_k~=~(-1)^k[a+\i b, 2k]2^{-2k},\qquad \b_k~=~(-1)^{k+1}[a+\i b, 2k-1]2^{-2k+1},
\\\\ \ds
[a+\i b,m]~=~{\Gamma\left({1\over 2}+a+\i b+m\right)\over m!\Gamma\left({1\over 2}+a+\i b-m\right)},\qquad  m~=~0,1,2,\ldots
\end{array}
\eeq
in the sense of that
\bel{J O}
\begin{array}{lr}\ds
\left({d\over d\rho}\right)^\ell\left\{ \J_{a+\i b}(\rho) ~-~\left({\pi \rho\over 2}\right)^{-{1\over 2}} \cos\left[\rho-(a+\i b) {\pi\over 2}-{\pi\over 4}\right]\right.
\\ \ds~~~
\left.~-~\left({\pi \rho\over 2}\right)^{-{1\over 2}}\sum_{k=1}^N  \cos\left[\rho-(a+\i b) {\pi\over 2}-{\pi\over 4}\right] \a_k \rho^{-2k}+\sin\left[\rho-(a+\i b) {\pi\over 2}-{\pi\over 4}\right]~ \b_k \rho^{-2k+1}\right\}
\\\\ \ds
~=~\O\left(\rho^{-2N-{1\over 2}}\right)\qquad N\ge0,\qquad \ell\ge0,\qquad \rho\mt\infty.
\end{array}
\eeq
The implied constant is bounded by $\C_a e^{\c|b|}$.

$\bullet$  For every $a,b\in\R$ and $\rho>0$, we have
\bel{J identity}
\J_{a-1+\i b} (\rho)~=~2\left[{a+\i b\over \rho}\right] \J_{a+\i b}(\rho)-\J_{a+1+\i b}(\rho).
\eeq
By putting together (\ref{Bessel}) and (\ref{J asymptotic})-(\ref{J O}) and using the identity in (\ref{J identity}), we obtain the norm estimate in below.

$\bullet$ For  every $a, b\in\R$ and $\rho>0$, 
\bel{J norm}
\left| {1\over \rho^{a+\i b}}~\J_{a+\i b}(\rho)\right|~\leq~\C_a~\left({1\over 1+\rho}\right)^{{1\over 2}+a}e^{\c|b|}.
\eeq

Let $z\in\Cx$.  
$\Omega^z$  is a distribution defined by analytic continuation from
\bel{Omega^z}
\begin{array}{ccc}\ds
\Re z<1,\qquad\Omega^z(x)~=~\pi^{-z}\Gamma^{-1}\left(1-z\right)  \left({1\over 1-|x|^2}\right)^{z}_+.
\end{array}
\eeq
$\bullet$  $\Omega^z$ can be equivalently defined by 
 \bel{Omega^z Transform} 
\begin{array}{lr}\ds
\Hat{\Omega}^z(\xi)~=~\left({1\over|\xi|}\right)^{{n\over 2}-z} \J_{{n\over 2}-z}\Big(2\pi|\xi|\Big),\qquad z\in\Cx.
\end{array}
\eeq
\begin{remark} From (\ref{Bessel}), (\ref{J asymptotic})-(\ref{J O}) and (\ref{J identity}), we conclude that $\Hat{\Omega}^z(\xi)$  for given $\xi\in\R^n$   is an analytic function of $z\in\Cx$.
\end{remark} 
Denote $\omega_{n-2}=2 \pi^{n-1\over 2} \Gamma^{-1}\left({n-1\over 2}\right)$ which is the area of $\mathds{S}^{n-2}$.
From (\ref{Omega^z}), we have
\bel{Omega^lambda Fourier transform}
\begin{array}{lr}\ds
\Hat{\Omega}^z(\xi)~=~\pi^{-z}  \Gamma^{-1}\left(1-z\right)\int_{|x|<1} e^{-2\pi\i x\cdot\xi} \left({1\over 1-|x|^2}\right)^z dx
\\\\ \ds~~~~~~~~~
~=~\pi^{-z}  \Gamma^{-1}\left(1-z\right)\int_0^\pi\left\{\int_0^1 e^{-2\pi\i |\xi|r\cos\vartheta} (1-r^2)^{-z} r^{n-1}dr\right\} \omega_{n-2} \sin^{n-2}\vartheta d\vartheta
\\\\ \ds~~~~~~~~~
~=~\omega_{n-2}\pi^{-z}  \Gamma^{-1}\left(1-z\right)\int_{-1}^1\left\{\int_0^1 e^{2\pi\i|\xi|rs} (1-r^2)^{-z} r^{n-1}dr\right\} (1-s^2)^{n-3\over 2} ds~~~~ \hbox{\small{($-s=\cos\vartheta$)}}
\\\\ \ds~~~~~~~~~
~=~\omega_{n-2}\pi^{-z}  \Gamma^{-1}\left(1-z\right)\int_0^1\left\{\int_{-1}^1 e^{2\pi\i|\xi|rs} (1-s^2)^{n-3\over 2} ds\right\}(1-r^2)^{-z} r^{n-1}dr
\\\\ \ds~~~~~~~~~
~=~\omega_{n-2}\pi^{-z}  \Gamma^{-1}\left(1-z\right)\int_0^1\left\{\int_{-1}^1 \cos(2\pi|\xi|rs) (1-s^2)^{n-3\over 2} ds\right\}(1-r^2)^{-z} r^{n-1}dr.
\end{array}
\eeq
Recall the Beta function identity:
\bel{Beta}
{\Gamma(z)\Gamma(w)\over \Gamma(z+w)}~=~\int_0^1 r^{z-1}(1-r)^{w-1}dr
\eeq
for every $\Re z>0$ and $\Re w>0$.

By writing out the Taylor expansion of the cosine function inside (\ref{Omega^lambda Fourier transform}), we find
\bel{Omega^lambda Fourier transform Sum}
\begin{array}{lr}\ds
\omega_{n-2}\pi^{-z}  \Gamma^{-1}\left(1-z\right)\int_0^1\left\{\int_{-1}^1 \cos(2\pi|\xi|rs) (1-s^2)^{n-3\over 2} ds\right\}(1-r^2)^{-z} r^{n-1}dr
\\\\ \ds
=~\omega_{n-2}\pi^{-z}  \Gamma^{-1}\left(1-z\right)
\\ \ds~~~~~
\sum_{k=0}^\infty (-1)^k {(2\pi|\xi|)^{2k}\over (2k)!}\left\{\int_{-1}^1 s^{2k} (1-s^2)^{n-3\over 2} ds\right\}\left\{\int_0^1 r^{2k+n-1}(1-r^2)^{-z} dr\right\}
\\\\ \ds
=~{1\over 2}\omega_{n-2}\pi^{-z}  \Gamma^{-1}\left(1-z\right)
\\ \ds~~~~~
\sum_{k=0}^\infty (-1)^k {(2\pi|\xi|)^{2k}\over (2k)!}\left\{\int_0^1 t^{k+{1\over 2}-1} (1-t)^{{n-1\over 2}-1} dt\right\}\left\{\int_0^1\rho^{k+{n\over 2}-1}(1-\rho)^{1-z-1} d\rho\right\}
\\ \ds~~~~~~~~~~~~~~~~~~~~~~~~~~~~~~~~~~~~~~~~~~~~~~~~~~~~~~~~~~~~~~~~~~~~~~~~~~~~~~~~~~~~~~~~~~~~~~~~~~~~~~
\hbox{\small{($t=s^2$,~$\rho=r^2$)}}
\\\\ \ds
=~{1\over 2}\omega_{n-2}\pi^{-z}  \Gamma^{-1}\left(1-z\right)\sum_{k=0}^\infty (-1)^k {(2\pi|\xi|)^{2k}\over (2k)!}~{\Gamma\left(k+{1\over 2}\right)\Gamma\left({n-1\over 2}\right)\over\Gamma\left(k+{n\over 2}\right)}~{ \Gamma\left(k+{n\over 2}\right) \Gamma\left(1-z\right)         \over    \Gamma\left(k+{n\over 2}+1-z\right)} 
\\ \ds~~~~~~~~~~~~~~~~~~~~~~~~~~~~~~~~~~~~~~~~~~~~~~~~~~~~~~~~~~~~~~~~~~~~~~~~~~~~~~~~~~~~~~~~~~~~~~~~~~~~~~~~~~~~~~
 \hbox{\small{by (\ref{Beta})}}
\\\\ \ds
=~ \pi^{{n-1\over 2}-z}\sum_{k=0}^\infty (-1)^k {(2\pi|\xi|)^{2k}\over (2k)!}~{\Gamma\left(k+{1\over 2}\right)        \over    \Gamma\left(k+{n\over 2}+1-z\right)}.
\end{array}
\eeq
On the other hand, we have
\bel{Omega cos}
\begin{array}{lr}\ds
\left({1\over|\xi|}\right)^{{n\over 2}-z} \J_{{n\over 2}-z}\Big(2\pi|\xi|\Big)
~=~\hbox{\small{$\pi^{{n-1\over 2}-z}  \Gamma^{-1}\left({n+1\over 2}-z\right)$}}\int_{-1}^1 e^{2\pi\i  |\xi| s} (1-s^2)^{{n-1\over 2}-z} ds\qquad\hbox{\small{by (\ref{Bessel})}}
\\\\ \ds~~~~~~~~~~~~~
~=~\hbox{\small{$\pi^{{n-1\over 2}-z}  \Gamma^{-1}\left({n+1\over 2}-z\right)$}}\int_{-1}^1 \cos\left(2\pi  |\xi| s\right) (1-s^2)^{{n-1\over 2}-z} ds
\\\\ \ds~~~~~~~~~~~~~
~=~\hbox{\small{$\pi^{{n-1\over 2}-z}  \Gamma^{-1}\left({n+1\over 2}-z\right)$}} \sum_{k=0}^\infty (-1)^k {(2\pi|\xi|)^{2k}\over (2k)!}\int_{-1}^1 s^{2k} (1-s^2)^{{n-1\over 2}-z} ds
\\\\ \ds~~~~~~~~~~~~~
~=~\hbox{\small{$\pi^{{n-1\over 2}-z}  \Gamma^{-1}\left({n+1\over 2}-z\right)$}} \sum_{k=0}^\infty (-1)^k {(2\pi|\xi|)^{2k}\over (2k)!}\int_0^1 \rho^{k+{1\over 2}-1} (1-\rho)^{{n+1\over 2}-z-1} d\rho\qquad \hbox{\small{($\rho=s^2$)}}
\\\\ \ds~~~~~~~~~~~~~
~=~\pi^{{n-1\over 2}-z} \sum_{k=0}^\infty (-1)^k {(2\pi|\xi|)^{2k}\over (2k)!}~{\Gamma\left(k+{1\over 2}\right)    \over    \Gamma\left(k+{n\over 2}+1-z\right)}\qquad\hbox{\small{by (\ref{Beta}).}}
\end{array}
\eeq

\section{Fourier transform of $\Lambda^\alpha$}
\setcounter{equation}{0}
We derive the formula of $\Hat{\Lambda}^\alpha(\xi,\tau)$ in (\ref{Lambda Transform})
by following the lines in p. $253-284$,
Chapter III of Gelfand and Shilov \cite{Gelfand-Shilov}.

Recall $\lambda(\alpha)={n+1\over 2}-\delta(\alpha),\delta(\alpha)=\left({n+1\over 2n}\right)\alpha, \alpha\in\Cx$. $\Hat{\hbox{U}}^\alpha$, $\Hat{\hbox{V}}^\alpha$ are distributions defined by analytic continuation from
\bel{UV}
\begin{array}{lr}\ds
\Hat{\hbox{U}}^\alpha(\xi,\tau)~=~\Gamma^{-1}\left(1-\delta(\alpha)\right)
  \left({1\over \tau^2-|\xi|^2}\right)^{\delta(\alpha)}_-,
\\\\ \ds
 \Hat{\hbox{V}}^\alpha(\xi,\tau)~=~\Gamma^{-1}\left(1-\delta(\alpha)\right)
  \left( {1\over \tau^2-|\xi|^2}\right)^{\delta(\alpha)}_+,
\qquad
  \hbox{\small{$\Re\delta(\alpha)<1$}}.  
\end{array}
\eeq
 $\Pi^\alpha$, $\Lambda^\alpha$ are distributions defined by analytic continuation from
\bel{PiLambda}
\begin{array}{lr}\ds
\Pi^\alpha(x,t)~=~\Gamma^{-1}\left(1-\lambda(\alpha)\right)
  \left({1\over t^2-|x|^2}\right)^{\lambda(\alpha)}_-,
\\\\ \ds
 \Lambda^\alpha(x,t)~=~\Gamma^{-1}\left(1-\lambda(\alpha)\right)
  \left({1\over t^2-|x|^2}\right)^{\lambda(\alpha)}_+,
\qquad
 \hbox{\small{$\Re\lambda(\alpha)<1$}}.  
  \end{array}
\eeq 
Let $a>0, b>0$ and
\bel{rho xi,tau}
\begin{array}{cc}\ds
\rho(\xi,\tau,a,b)~=~\Bigg\{ \left[\tau^2-|\xi|^2\right]^2+\left[a \tau^2+b|\xi|^2\right]^2\Bigg\}^{1\over 2},
\\\\ \ds
 \cos\theta(\xi,\tau,a,b)~=~{\tau^2-|\xi|^2\over\rho(\xi,\tau,a,b)}.
 \end{array}
\eeq
Suppose $\Re\alpha<n$ which is equivalent to $\Re\delta(\alpha)<{n+1\over 2}$.   Assert
 \bel{P,R,ab}
 \begin{array}{lr}\ds
\Hat{\hbox{P}}^{\alpha~a~b}(\xi,\tau)~=~\left[ \tau^2-|\xi|^2+\i a \tau^2+\i b |\xi|^2\right]^{-\delta(\alpha)}
\\\\ \ds~~~~~~~~~~~~~~~~~~
~=~\rho(\xi,\tau,a,b)^{-\delta(\alpha)} e^{-\i\alpha \theta(\xi,\tau,a,b)},
\\\\ \ds
\Hat{\hbox{R}}^{\alpha~a~b}(\xi,\tau)~=~\left[ \tau^2-|\xi|^2-\i a \tau^2-\i b |\xi|^2\right]^{-\delta(\alpha)}
\\\\ \ds~~~~~~~~~~~~~~~~~~
~=~\rho(\xi,\tau,a,b)^{-\delta(\alpha)} e^{\i\alpha \theta(\xi,\tau,a,b)},\qquad \hbox{\small{$|\xi|\neq|\tau|$}}.
\end{array}
\eeq
$\Hat{\hbox{P}}^{\alpha~a~b}$ and $\Hat{\hbox{R}}^{\alpha~a~b}$ are distributions defined in $\R^{n+1}$ agree with (\ref{P,R,ab}) whenever $|\xi|\neq|\tau|$. Define
\bel{PR Limit}
 \Hat{\hbox{P}}^\alpha~=~\lim_{a\mt0,~b\mt0} \Hat{\hbox{P}}^{\alpha~a~b},\qquad \Hat{\hbox{R}}^\alpha~=~\lim_{a\mt0,~b\mt0} \Hat{\hbox{R}}^{\alpha~a~b},\qquad \hbox{\small{$\Re\alpha<n$}}.
 \eeq
Suppose $\Re\alpha>0$ which is equivalent to $\lambda(\alpha)<{n+1\over 2}$. Consider
\bel{Phi,Psi,ab}
 \begin{array}{lr}\ds
\Phi^{\alpha~a~b}(x,t)~=~\left[ t^2-|x|^2+\i a t^2+\i b |x|^2\right]^{-\lambda(\alpha)}
\\\\ \ds~~~~~~~~~~~~~~~~~~
~=~\rho(x,t,a,b)^{-\lambda(\alpha)} e^{-\i\lambda(\alpha) \theta(x,t,a,b)},
\\\\ \ds
\Psi^{\alpha~a~b}(x,t)~=~\left[ t^2-|x|^2-\i a t^2-\i b |x|^2\right]^{-\lambda(\alpha)}
\\\\ \ds~~~~~~~~~~~~~~~~~~
~=~\rho(x,t,a,b)^{-\lambda(\alpha)} e^{\i\lambda(\alpha) \theta(x,t,a,b)},\qquad\hbox{\small{$|x|\neq|t|$}}.
\end{array}
\eeq
$\Phi^{\alpha~a~b}$ and $\Psi^{\alpha~a~b}$ are distributions defined in $\R^{n+1}$ agree with (\ref{Phi,Psi,ab}) whenever $|x|\neq|t|$. Define
\bel{PhiPsi Limit}
\Phi^\alpha~=~\lim_{a\mt0,~b\mt0} \Phi^{\alpha~a~b},\qquad \Psi^\alpha~=~\lim_{a\mt0,~b\mt0} \Psi^{\alpha~a~b},\qquad\hbox{\small{$\Re\alpha>0$}}.
\eeq
\begin{remark}
$\Hat{\hbox{P}}^\alpha$, $\Hat{\hbox{R}}^\alpha$ are analytic for 
$\Re\alpha<n$. $\Phi^\alpha$, $\Psi^\alpha$ are analytic for $\Re\alpha>0$.

Regarding details  can be found in $2.2$, 
 Chapter III of  Gelfand and Shilov \cite{Gelfand-Shilov}.
 \end{remark}

Define $\Hat{\hbox{P}}^\alpha$ and $\Hat{\hbox{R}}^\alpha$ by analytic continuation from (\ref{PR Limit}).
Recall (\ref{UV}) and (\ref{rho xi,tau})-(\ref{PR Limit}). We find
\bel{PR}
\begin{array}{lr}\ds
 \Hat{\hbox{P}}^\alpha
 ~=~\Gamma\left(1-\delta(\alpha)\right)\left[e^{-\i\pi \delta(\alpha)}\Hat{\hbox{U}}^\alpha+ \Hat{\hbox{V}}^\alpha\right],
\qquad
\Hat{\hbox{R}}^\alpha~=~\Gamma\left(1-\delta(\alpha)\right)\left[e^{\i\pi\delta(\alpha)}\Hat{\hbox{U}}^\alpha+ \Hat{\hbox{V}}^\alpha\right].
\end{array}
\eeq
Define $\Phi^\alpha$ and $\Psi^\alpha$ by analytic continuation from (\ref{PhiPsi Limit}).
Recall (\ref{PiLambda}), (\ref{rho xi,tau}) and (\ref{Phi,Psi,ab})-(\ref{PhiPsi Limit}).
We have
\bel{PhiPsi} 
\begin{array}{lr}\ds
\Phi^\alpha
~=~\Gamma\left(1-\lambda(\alpha)\right)\left[e^{-\i\pi \lambda(\alpha)} \Pi^\alpha+\Lambda^\alpha\right],
\qquad
\Psi^\alpha
~=~\Gamma\left(1-\lambda(\alpha)\right)\left[e^{\i\pi \lambda(\alpha)} \Pi^\alpha+\Lambda^\alpha\right].
\end{array}
\eeq
For $\Im z\neq0, \Im w\neq0$ and $0<\Re\alpha<n$, assert
\[
Q^{\alpha~z~w}(x,t)~=~\left\{{1\over z |x|^2+ wt^2}\right\}^{\lambda(\alpha)},\qquad \hbox{\small{$(x,t)\neq(0,0)$}}.
\]
Let $a>0, b>0$. From  direct computation, we find
\bel{Q ab Transform}
\begin{array}{lr}\ds
\Hat{Q}^{\alpha~\i a~\i b}(\xi,\tau)~=~(-\i)^{\lambda(\alpha)}\iint_{\R^{n+1}} e^{-2\pi\i \big[x\cdot\xi+t\tau\big]} \left\{{1\over a |x|^2+ b t^2}\right\}^{\lambda(\alpha)} dxdt
\\\\ \ds~~~~~~~~~~~~~~~~~~~
~=~(-\i)^{\lambda(\alpha)}{1\over (\sqrt{a})^n\sqrt{b}}\iint_{\R^{n+1}} e^{-2\pi\i \big[x\cdot\xi/\sqrt{a}+t\tau/\sqrt{b}\big]} \left({1\over  |x|^2+  t^2}\right)^{\lambda(\alpha)} dxdt
\\ \ds~~~~~~~~~~~~~~~~~~~~~~~~~~~~~~~~~~~~~~~~~~ ~~~~~~~~~~~~~~~~~~~~~~~~~~~~~~~~~~~~~~~~~~~~
\hbox{\small{$x\mt x/\sqrt{a}$,~~$t\mt t/\sqrt{b}$}}
\\\\ \ds~~~~~~~~~~~~~~~~~~~
~=~(-\i)^{\lambda(\alpha)}{\pi^{-{n+1\over 2}+2\lambda(\alpha)}\over (\sqrt{a})^n\sqrt{b}}{\Gamma\left({n+1\over 2}-\lambda(\alpha)\right)\over \Gamma\left(\lambda(\alpha)\right)}\left\{{1\over  |\xi|^2/a+  \tau^2/b}\right\}^{{n+1\over 2}-\lambda(\alpha)},
\quad
\hbox{\small{$0<\Re\alpha<n$}}
\end{array}
\eeq
whenever $(\xi,\tau)\neq(0,0)$. 

The last equality in (\ref{Q ab Transform}) is derived by using 
\[
\int_{\R^\N}e^{-2\pi\i x\cdot\xi} |x|^{\gamma-\N}dx~=~{\pi^{\N-\gamma\over 2}\Gamma\left({\gamma\over 2}\right)\over \pi^{\gamma\over 2}\Gamma\left({\N-\gamma\over 2}\right)} |\xi|^{-\gamma},\qquad \hbox{\small{$\xi\neq0$}},\qquad\hbox{\small{$0<\Re\gamma<\N$}}.
\]
Replace $a=-\i z$ and $b=-\i w$ inside (\ref{Q ab Transform}). We have
\bel{Q zv Transform}
\begin{array}{lr}\ds
\Hat{Q}^{\alpha~z~w}(\xi,\tau)~=~\iint_{\R^{n+1}} e^{-2\pi\i \left(x\cdot\xi+t\tau\right)} \left\{{1\over z |x|^2+ w t^2}\right\}^{\lambda(\alpha)} dxdt
\\\\ \ds~~~~~~~~~~~~~~~~~~
~=~{\pi^{\lambda(\alpha)-\delta(\alpha)}\over (\sqrt{ z})^n\sqrt{ w}}{\Gamma\left(\delta(\alpha)\right)\over \Gamma\left(\lambda(\alpha)\right)}\left\{{1\over  |\xi|^2/z+  \tau^2/w}\right\}^{\delta(\alpha)},\qquad \hbox{\small{$(\xi,\tau)\neq(0,0)$}}.
\end{array}
\eeq
\begin{remark}
(\ref{Q zv Transform}) is true  for every $\Im z>0$ and $\Im w>0$.
\end{remark}
Let $w=\i b$ fixed. 
Both sides of (\ref{Q zv Transform}) are analytic whenever $\Im z>0$.  Moreover, they are equal at $z=\i a$ for every $a>0$.  By analytic continuation, this equality holds for $\Im z>0$.
Given $z,~\Im z>0$,  a vice versa argument shows that (\ref{Q zv Transform}) is true for  $\Im w>0$. 

On the other hand, by carrying out the same estimate in (\ref{Q ab Transform}), we find
\bel{Q ab Transform -}
\begin{array}{lr}\ds
\Hat{Q}^{\alpha~-\i a~-\i b}(\xi,\tau)~=~(-\i)^{-\lambda(\alpha)}{\pi^{-{n+1\over 2}+2\lambda(\alpha)}\over (\sqrt{a})^n\sqrt{b}}{\Gamma\left({n+1\over 2}-\lambda(\alpha)\right)\over \Gamma\left(\lambda(\alpha)\right)}\left\{{1\over  |\xi|^2/a+  \tau^2/b}\right\}^{{n+1\over 2}-\lambda(\alpha)},
\\\\ \ds~~~~~~~~~~~~~~~~~~~~~~~~~~~~~~~~~~~~~~~~~~~~~~~~~~~~~~~~~~~~~~~~~~~~~~~~~~~~~~~\hbox{\small{$(\xi,\tau)\neq(0,0)$}},\qquad
\hbox{\small{$0<\Re\alpha<n$}}.
\end{array}
\eeq
Replace $a=\i z$ and $b=\i w$ inside (\ref{Q ab Transform -}),  we obtain (\ref{Q zv Transform}) again. Furthermore, an analogue argument below {\bf Remark B.3} shows (\ref{Q zv Transform}) hold for $\Im z<0, \Im w<0$.
Let
\bel{rho a}
\rho(a)~=~\sqrt{1+a^2},\qquad \cos\theta(a)~=~-1/\rho(a),\qquad \sin\theta(a)~=~a/\rho(a).
\eeq
Note that 
\bel{limit a to zero}
(-1\pm\i a)^{n\over 2}~=~\rho(a)^{n\over 2} e^{\pm\i \theta \left({n\over 2}\right)}~\mt~ e^{\pm\i\pi \left({n\over 2}\right)}\qquad \hbox{\small{as}} \qquad a\mt0.
\eeq
Consider $\Hat{Q}^{\alpha~z~w}(\xi,\tau)\doteq\Hat{Q}^{\alpha~a~b}(\xi,\tau)$ with $z=-1+\i a$ and $w=1+\i b$ inside (\ref{Q zv Transform}).
We have
\bel{Q zv Transform +}
\begin{array}{lr}\ds
\Hat{Q}^{\alpha~a~b}(\xi,\tau)
~=~\iint_{\R^{n+1}} e^{-2\pi\i \left(x\cdot\xi+t\tau\right)} \left\{{1\over (-1+\i a) |x|^2+ (1+\i b) t^2}\right\}^{\lambda(\alpha)} dxdt
\\\\ \ds~~~~~~~~~~~~~~~~~~
~=~{\pi^{\lambda(\alpha)-\delta(\alpha)}\over (-1+\i a)^{n\over2}(1+\i b)^{1\over 2}}{\Gamma\left(\delta(\alpha)\right)\over \Gamma\left(\lambda(\alpha)\right)}\left\{{(-1+\i a)(1+\i b)\over  (1+\i b)|\xi|^2+ (-1+\i a) \tau^2}\right\}^{\delta(\alpha)}
\end{array}
\eeq
for $(\xi,\tau)\neq(0,0)$ and $0<\Re\alpha<n$.

Similarly, for $z=-1-\i a$ and $w=1-\i b$, we find
\bel{Q zv Transform -}
\begin{array}{lr}\ds
\Hat{Q}^{\alpha~a~b}(\xi,\tau)
~=~{\pi^{\lambda(\alpha)-\delta(\alpha)}\over (-1-\i a)^{n\over2}(1-\i b)^{1\over 2}}{\Gamma\left(\delta(\alpha)\right)\over \Gamma\left(\lambda(\alpha)\right)}\left\{{(-1-\i a)(1-\i b)\over  (1-\i b)|\xi|^2+ (-1-\i a) \tau^2}\right\}^{\delta(\alpha)},
\\ \ds~~~~~~~~~~~~~~~~~~~~~~~~~~~~~~~~~~~~~~~~~~~~~~~~~~~~~~~~~~~~~~~~~~~~~~~~~~~~~~
\hbox{\small{$(\xi,\tau)\neq(0,0)$,\qquad $0<\Re\alpha<n$}}.
\end{array}
\eeq
Define $\Hat{Q}^{\alpha~a~b}$, $a>0, b>0,  0<\alpha<n$ as a distribution in $\R^{n+1}$ agree with (\ref{Q zv Transform +})-(\ref{Q zv Transform -}) whenever $|\xi|\neq|\tau|$. 

Recall (\ref{P,R,ab})-(\ref{PhiPsi Limit}). 
By using (\ref{limit a to zero}) and taking $a\mt0, b\mt0$, we find
\bel{Phi, R, Psi, P, Transform}
\begin{array}{lr}\ds
\Hat{\Phi}^\alpha~=~\pi^{\lambda(\alpha)-\delta(\alpha)}e^{-\i\pi\left({ n\over 2}\right)}~{\Gamma\left(\delta(\alpha)\right)\over \Gamma\left(\lambda(\alpha)\right)}~\Hat{\hbox{R}}^\alpha,
\qquad
\Hat{\Psi}^\alpha~=~\pi^{\lambda(\alpha)-\delta(\alpha)}e^{\i\pi\left({ n\over 2}\right)}~{\Gamma\left(\delta(\alpha)\right)\over \Gamma\left(\lambda(\alpha)\right)}~\Hat{\hbox{P}}^\alpha,\qquad \hbox{\small{$0<\Re\alpha<n$}}.
\end{array}
\eeq
Moreover,  $\Hat{\hbox{P}}^\alpha, \Hat{\hbox{R}}^\alpha$ and $\Phi^\alpha, \Psi^\alpha$ can be given   in terms of $\Pi^\alpha, \Lambda^\alpha$ and $\hbox{U}^\alpha, \hbox{V}^\alpha$ as (\ref{PR}) and (\ref{PhiPsi}). 
Let $0<\Re\alpha<n$. A direct computation shows
\bel{Lambda computa}
\begin{array}{lr}\ds
\Hat{\Lambda}^\alpha\left[ e^{i\pi\lambda(\alpha)}-e^{-i\pi\lambda(\alpha)}\right]
~=~\Gamma^{-1}(1-\lambda(\alpha))\Bigg[ e^{\i\pi\lambda(\alpha)} \Hat{\Phi}^\alpha-e^{-\i\pi\lambda(\alpha)}\Hat{\Psi}^\alpha\Bigg]
\\\\ \ds~~~~~~~~~~~~~~~~~~~~~~~~~~~~~~~~~~~
~=~\Gamma^{-1}(1-\lambda(\alpha)) \Gamma^{-1}(\lambda(\alpha)) 
\pi^{\lambda(\alpha)-\delta(\alpha)}\Gamma(\delta(\alpha)) 
\Bigg[ e^{-\i\pi\big[{n\over 2}-\lambda(\alpha)\big]} \Hat{\hbox{R}}^\alpha-e^{\i\pi\big[{n\over 2}-\lambda(\alpha)\big]} \Hat{\hbox{P}}^\alpha\Bigg].
\end{array}
\eeq
From  (\ref{Lambda computa}), by using the identity 
$\Gamma^{-1}(1-z)\Gamma^{-1}(z)={ \sin \pi z\over \pi}, z\in\Cx$, we obtain
\bel{Lambda Transform Formula}
\begin{array}{lr}\ds
\Hat{\Lambda}^\alpha
~=~\pi^{\lambda(\alpha)-\delta(\alpha)-1} \Gamma\left(\delta(\alpha)\right){1\over 2\i} \Bigg\{- e^{\i\pi \big[{n\over 2}-\lambda(\alpha)\big]} \Hat{\hbox{P}}^\alpha+e^{-\i\pi \big[{n\over 2}-\lambda(\alpha)\big] }\Hat{\hbox{R}}^\alpha\Bigg\}
\\\\ \ds~~~~~
~=~\pi^{\lambda(\alpha)-\delta(\alpha)-1} \Gamma\left(\delta(\alpha)\right)\Gamma(1-\delta(\alpha))
\\ \ds~~~~~~~~~~
{1\over 2\i} \Bigg\{- e^{\i\pi \big[{n\over 2}-\lambda(\alpha)\big]} \Big[e^{-\i\pi\delta(\alpha)}\Hat{\hbox{U}}^\alpha+\Hat{\hbox{V}}^\alpha\Big]+e^{-\i\pi \big[{n\over 2}-\lambda(\alpha)\big] }\Big[e^{\i\pi\delta(\alpha)}\Hat{\hbox{U}}^\alpha+\Hat{\hbox{V}}^\alpha\Big]\Bigg\}\qquad \hbox{\small{by (\ref{PR})}}
\\\\ \ds~~~~~
~=~\pi^{\lambda(\alpha)-\delta(\alpha)-1} \Gamma\left(\delta(\alpha)\right)\Gamma(1-\delta(\alpha))
\\ \ds~~~~~~~~~~
{1\over 2\i} \Bigg\{- e^{\i\pi \big[\delta(\alpha)-{1\over2}\big]} \Big[e^{-\i\pi\delta(\alpha)}\Hat{\hbox{U}}^\alpha+\Hat{\hbox{V}}^\alpha\Big]+e^{-\i\pi \big[\delta(\alpha)-{1\over 2}\big] }\Big[e^{\i\pi\delta(\alpha)}\Hat{\hbox{U}}^\alpha+\Hat{\hbox{V}}^\alpha\Big]\Bigg\}
\\\\ ~~~~~
~=~\pi^{\lambda(\alpha)-\delta(\alpha)-1}\Gamma(\delta(\alpha))\Gamma\left(1-\delta(\alpha)\right)
\Bigg\{ \Hat{\hbox{U}}^\alpha-\sin\pi\left(\delta(\alpha)-{1\over 2}\right)\Hat{\hbox{V}}^\alpha\Bigg\}
\end{array}
\eeq
for $\lambda(\alpha)\notin\Z$ and $0<\Re\alpha<n$. 

Recall {\bf Remark B.1}. $\Hat{\hbox{P}}^\alpha$, $\Hat{\hbox{R}}^\alpha$ are analytic for 
$\Re\alpha<n$. By the principle of analytic continuation, (\ref{Lambda Transform Formula}) holds for $0<\Re\alpha<n$.

Let $0<\Re\delta(\alpha)<1$. $\Hat{\hbox{U}}^\alpha$ and $\Hat{\hbox{V}}^\alpha$ agree with  $\Hat{\hbox{U}}^\alpha(\xi,\tau)$ and $\Hat{\hbox{V}}^\alpha(\xi,\tau)$ in (\ref{UV}) whenever $|\xi|\neq|\tau|$. We have
\[
\begin{array}{lr}
  \Hat{\Lambda}^\alpha(\xi,\tau)~=~
\pi^{\lambda(\alpha)-\delta(\alpha)-1}\Gamma\left(\delta(\alpha)\right)
 \left\{ ~{\ds \left({1\over \tau^2-|\xi|^2}\right)^{\delta(\alpha)}_-} - \sin\pi\left(\delta(\alpha)-{1\over 2}\right){\ds \left( {1\over \tau^2- |\xi|^2}\right)^{\delta(\alpha)}_+}~\right\}.
  \end{array}
\]
 An analogue of the above estimates also shows
\bel{V inverse transform formula}
\begin{array}{lr}\ds
\hbox{V}^\alpha~=~\pi^{\delta(\alpha)-\lambda(\alpha)-1}\Gamma(\lambda(\alpha))\Gamma\left(1-\lambda(\alpha)\right)
\Bigg\{\Pi^\alpha-\sin\pi\left(\lambda(\alpha)-{1\over 2}\right)\Lambda^\alpha\Bigg\},\qquad\hbox{\small{$0<\Re\alpha<n$}}. 
\end{array}
\eeq
In particular, we find
\bel{V inverse transform}
\begin{array}{lr}
\hbox{V}^\alpha(x,t)~=~
\pi^{\delta(\alpha)-\lambda(\alpha)-1}\Gamma\left(\lambda(\alpha)\right)
 \left\{ ~{\ds \left({1\over t^2-|x|^2}\right)^{\lambda(\alpha)}_-} - \sin\pi\left(\lambda(\alpha)-{1\over 2}\right){\ds \left( {1\over t^2- |x|^2}\right)^{\lambda(\alpha)}_+}~\right\},
 \\ \ds~~~~~~~~~~~~~~~~~~~~~~~~~~~~~~~~~~~~~~~~~~~~~~~~~~~~~~~~~~~~~~~~~~~~~~~~~~~~~~~~~~~~~~~~~~~~
 \hbox{\small{$|x|\neq|t|$}},\qquad\hbox{\small{$0<\Re\lambda(\alpha)<1$}}.
  \end{array}
\eeq

{\small zipeng.wang@cantab.net}

\end{document}